\documentstyle[12pt,amssymb,amscd]{amsart} 
\newtheorem{theorem}{Theorem}[section]
\newtheorem{lemma}[theorem]{Lemma}
\newtheorem{proposition}[theorem]{Proposition}
\newtheorem{corollary}[theorem]{Corollary}

\newtheorem{thm}{Theorem}[section]
\newtheorem{lem}[theorem]{Lemma}

\newtheorem{cor}[theorem]{Corollary}

\theoremstyle{definition}
\newtheorem{definition}[theorem]{Definition}

\theoremstyle{remark}

\numberwithin{equation}{section}

\newcommand{\QT}{Q_{\theta}}
\newcommand{\TQT}{\tilde{Q}_{\theta}}

\newcommand{\TQN}{\tilde{Q}_{\nu}}
\newcommand{\KTN}{K_{\theta,\, \nu}}

\newcommand{\KT}{K_{\theta}}
\newcommand{\KN}{K_{\nu}}
\newcommand{\LT}{\Lambda_{\theta}}
\newcommand{\LN}{\Lambda_{\nu}}
\newcommand{\FTN}{F_{\theta,\, \nu}}
\newcommand{\BTN}{B_{\theta,\, \nu}}
\newcommand{\TBTN}{\tilde{B}_{\theta,\, \nu}}
\newcommand{\FT}{f_{\theta}}
\newcommand{\FN}{f_{\nu}}
\newcommand{\VT}{\varphi_{\theta}}
\newcommand{\VN}{\varphi_{\nu}}
\newcommand{\PT}{\psi_{\theta}}
\newcommand{\PN}{\psi_{\nu}}
\newcommand{\ZT}{\zeta_{\theta}}
\newcommand{\ZN}{\zeta_{\nu}}
\newcommand{\GT}{\eta_{\theta}}
\newcommand{\GN}{\eta_{\nu}}
\newcommand{\ST}{S_{\theta}}
\newcommand{\SN}{S_{\nu}}
\newcommand{\cheb}{f_{\mbox{\tiny cheb}}}
\newcommand{\gcheb}{\gamma_{\mbox{\tiny cheb}}}
\newcommand{\rat}{\operatorname{\bf{Rat}}_2}

\newcommand{\ov}{\overline}
\newcommand{\io}{\iota}
\newcommand{\dia}{\mbox{diam}}
\renewcommand{\mod}{\operatorname{mod}}
\newcommand{\CCC}{{\mathcal C}}
\newcommand{\PP}{{\Bbb P}}
\newcommand{\CC}{{\Bbb C}}
\newcommand{\RR}{{\Bbb R}}
\newcommand{\TT}{{\Bbb T}}
\newcommand{\ZZ}{{\Bbb Z}}

\newcommand{\DD}{{\Bbb D}}
\newcommand{\iso}{\stackrel{\simeq}{\longrightarrow}}
\newcommand{\sm}{\smallsetminus}
\newcommand{\bp}{{Blaschke product }}
\newcommand{\mate}{\sqcup}
\newcommand{\tmate}{\sqcup_{\cal T}}
\newcommand{\fmate}{\sqcup_{\cal F}}
\newcommand{\infeq}{{\sim}_{\infty}}
\newcommand{\rayeq}{{\sim}_r}

\newcommand{\diam}{\operatorname{diam}}

\newcommand{\dist}{\operatorname{dist}}

\newcommand{\tl}{\tilde}

\newcommand{\eps}{\epsilon}
\newcommand{\es}{\emptyset}

\newcommand{\bd}{\partial}

\newcommand{\cC}{{\cal C}}

\newcommand{\cD}{{\cal D}}

\newcommand{\C}{\Bbb C}

\newcommand{\D}{\Bbb D}
\newcommand{\T}{\Bbb T}

\newcommand{\ang}[2]{\widehat{(#1,#2)}}

\newcommand{\M}{{\cal M}}

\renewcommand{\marginpar}[1]{}
\catcode`\@=12

\def\Empty{}
\newcommand\oplabel[1]{
  \def\OpArg{#1} \ifx \OpArg\Empty {} \else
  	\label{#1}
  \fi}
		
\long\def\realfig#1#2#3#4{
\begin{figure}[htp]
\centerline{\psfig{figure=#2,width=#4}}
\caption[#1]{#3}
\oplabel{#1}
\end{figure}}

\newcommand{\comm}[1]{}
    
\input{psfig}

\newcommand{\thmref}[1]{Theorem~\ref{#1}}
\newcommand{\propref}[1]{Proposition~\ref{#1}}
\newcommand{\secref}[1]{\S\ref{#1}}
\newcommand{\lemref}[1]{Lemma~\ref{#1}}
\newcommand{\corref}[1]{Corollary~\ref{#1}} 
\newcommand{\figref}[1]{Fig.~\ref{#1}}

\begin{document}
\title{Mating Siegel Quadratic Polynomials} 
\author{Michael Yampolsky, Saeed Zakeri}
\thanks{The first author was partially supported
by NSF grant DMS-9804606}


\begin{abstract}
Let $F$ be a quadratic rational map of the sphere which has two fixed Siegel disks with bounded type rotation numbers $\theta$ and $\nu$. 
Using a new degree $3$ Blaschke product model for the dynamics of $F$ and an adaptation of complex a priori bounds for renormalization of 
 critical circle maps, we prove that $F$ can be realized as the mating 
of two Siegel quadratic polynomials with the corresponding rotation numbers $\theta$ and $\nu$.
\end{abstract} 
 
\maketitle

\thispagestyle{empty}
\def\SBIMSMark#1#2#3{
 \font\SBF=cmss10 at 10 true pt
 \font\SBI=cmssi10 at 10 true pt
 \setbox0=\hbox{\SBF Stony Brook IMS Preprint \##1}
 \setbox2=\hbox to \wd0{\hfil \SBI #2}
 \setbox4=\hbox to \wd0{\hfil \SBI #3}
 \setbox6=\hbox to \wd0{\hss
             \vbox{\hsize=\wd0 \parskip=0pt \baselineskip=10 true pt
                   \copy0 \break%
                   \copy2 \break%
                   \copy4 \break}}
 \dimen0=\ht6   \advance\dimen0 by \vsize \advance\dimen0 by 8 true pt
                \advance\dimen0 by -\pagetotal
 \dimen2=\hsize \advance\dimen2 by .25 true in
%
%
  \openin2=publishd.tex
  \ifeof2\setbox0=\hbox to 0pt{}
  \else 
     \setbox0=\hbox to 3.1 true in{
                \vbox to \ht6{\hsize=3 true in \parskip=0pt  \noindent  
                \input publishd.tex 
                \vfill}}
  \fi
  \closein2
  \ht0=0pt \dp0=0pt
 \ht6=0pt \dp6=0pt
 \setbox8=\vbox to \dimen0{\vfill \hbox to \dimen2{\copy0 \hss \copy6}}
 \ht8=0pt \dp8=0pt \wd8=0pt
 \copy8
 \message{*** Stony Brook IMS Preprint #1, #2 ***}
}

\SBIMSMark{1998/8}{July 1998}{}

\tableofcontents
\section{Introduction}
\label{sec:intro}

\subsection{Mating: Definitions and some history.}
Mating quadratic polynomials is a topological construction suggested by Douady and Hubbard \cite{Do1}
to partially parametrize quadratic rational maps of the Riemann sphere by pairs of quadratic
polynomials. Some results on matings of higher degree maps exist, but we will not discuss
them in this paper.
While there exist several, presumably equivalent, ways of describing the construction of mating,
the following approach is perhaps the most standard.
Consider two monic quadratic polynomials $f_1$ and $f_2$ whose filled Julia sets $K(f_i)$ are locally-connected.
  For each $f_i$, let $\Phi_i$ denote the 
conformal isomorphism  between the basin of infinity $\ov{\CC}\sm K(f_i)$
and ${\ov \CC}\sm \ov \DD$, with $\Phi_i(\infty)=\infty$ and $\Phi_i'(\infty)=1$. 
These {\it B\"ottcher maps} conjugate the polynomials to the squaring map:
$$\begin{CD}
\ov \CC\sm K(f_i)@>\Phi_i>>\ov  \CC\sm \ov \DD\\
@VV{f_i}V	@VV{z\mapsto z^2}V\\
\ov \CC\sm K(f_i)@>\Phi_i>>\ov  \CC\sm \ov \DD
\end{CD}
$$

\noindent
 By the Carath\'eodory's Theorem the inverse map $\Phi_i^{-1}$ has a continuous
extension 
$$\Phi_i^{-1}:\bd \DD\to J(f_i),$$ 
where the Julia set $J(f_i)=\bd K(f_i)$ is the topological boundary of the filled Julia set.
 The induced parametrization
$$\gamma_i(t)\equiv \Phi_i^{-1}(e^{2\pi i t}):\TT= \RR/\ZZ\to J(f_i)$$
 is commonly referred to as the {\it Carath\'eodory loop} of $J(f_i)$. Note that by the above commutative diagram, $\gamma_i(2t)=f_i(\gamma_i(t))$. 
Consider the topological space 
$$X=(K(f_1)\sqcup K(f_2))/(\gamma_1(t)\sim\gamma_2(-t))$$
obtained by gluing the two filled Julia sets along their Carath\'eodory loops
in reverse directions. \vspace{0.12in}\\
{\bf Definition I.} Assume that the space $X$ as defined above is homeomorphic to the $2$-sphere $S^2$.
 Then the pair of polynomials $(f_1,f_2)$ is called {\it topologically mateable}. The induced map of $S^2$
$$f_1\tmate f_2 = (f_1|_{K_1}\sqcup f_2|_{K_2})/(\gamma_1(t)\sim\gamma_2(-t))$$
is the {\it topological mating} of $f_1$ and $f_2$. \vspace{0.12in}

It may seem surprising at this point that
topologically mateable quadratics even exist, however, we shall see below that such examples are abundant. For any mateable pair $(f_1,f_2)$, their topological mating is a degree 2 branched covering of the sphere, and it is natural to ask whether it possesses an invariant conformal structure.\vspace{0.12in}\\
{\bf Definition II.} A quadratic rational map $F:\ov \CC\to \ov \CC$ is called a {\it conformal mating}, or simply a {\it mating}, of $f_1$ and $f_2$,
$$F=f_1\mate f_2,$$
if it is conjugate to the topological mating $f_1\tmate f_2$ by a homeomorphism which is conformal in the interiors of $K(f_1)$ and $K(f_2)$ in case there is an interior. If such $F$ is unique up to conjugation by a M\"{o}bius transformation, we refer to it as {\it the} mating of $f_1$ and $f_2$. \vspace{0.12in}

Before proceeding to formulate the known existence results, let us describe another equivalent method of defining a mating.
Let $\copyright$ denote the complex plane $\CC$ compactified by adjoining a circle of directions
at infinity, $\{\infty\cdot e^{2\pi i t}| t\in \TT \}$ with the natural topology. Each $f_i$ extends continuously to a copy of $\copyright_i$, acting as the squaring map $z\mapsto z^2$ on the circle at infinity. 
Gluing the disks $\copyright_i$ together via the equivalence relation $\infeq$ identifying the point $\infty\cdot e^{2\pi i t}\in \copyright_1$ with $\infty\cdot e^{-2\pi i t}\in \copyright_2$,
we obtain a $2$-sphere $(\copyright_1\sqcup \copyright_2)/\infeq$.
The well-defined map $f_1\fmate f_2$ on this sphere given by $f_i$ on $\copyright_i$ is a degree $2$ branched covering of the
sphere with an invariant equator. We shall refer to this map as the {\it formal mating} of $f_1$, $f_2$.

Recall that the {\it external ray of $f_i$ at angle $t$} is the preimage
$$R_i(t)= \Phi_i^{-1}(\{ re^{2\pi i t}|r>1 \})$$
for $t\in \TT$.
Let $\hat R_i(t)$ denote the closure of $R_i(t)$ in $\copyright_i$. The {\it ray equivalence relation} $\rayeq$
on $(\copyright_1\sqcup \copyright_2)/\infeq$ is defined as follows. The points $z$ and $w$ are equivalent,
$z\rayeq w$ if and only if there exists a collection of closed rays $\hat R_j=\hat R_{i}(t_j)$, 
$i \in \{1,2 \}$ and $j=1,\ldots,n$, 
such that $z\in \hat R_1$, $w\in \hat R_n$ and $\hat R_j\cap \hat R_{j+1}\ne\emptyset$ for $j=1, \ldots, n-1$.
It follows immediately from the definition that if $f_1$ and $f_2$ are topologically mateable, then 
 the quotient of $(\copyright_1\sqcup \copyright_2)/\infeq$ modulo $\rayeq$ is again a $2$-sphere, and 
$$(f_1\fmate f_2)/\rayeq \simeq f_1\tmate f_2.$$
Finally, let us formulate another definition of conformal mating, equivalent to the previously given, but more convenient
 for further application:
\vspace{0.12in}\\
{\bf Definition IIa.} Let $f_1$ and $f_2$ be quadratic polynomials with 
locally-connected Julia sets. A quadratic rational map $F$ of the Riemann sphere is called a {\it conformal mating} 
of $f_1$ and $f_2$ if there exist continuous semiconjugacies
$$\varphi_i:K(f_i)\to\ov \CC,\text{ with }\varphi_i\circ f_i=F\circ\varphi_i,$$
conformal in the interiors of the filled Julia sets in case there is an interior, such that $\varphi_1(K(f_1))\cup\varphi_2(K(f_2))=\ov\CC$ and for $i,j=1,2$, $\varphi_i(z)=\varphi_j(w)$ if and only if
 $z\rayeq w$. \vspace{0.12in}

We are now prepared to give an account of known results. The simplest example of a non-mateable pair
is given by quadratic polynomials $f_{c_1}(z)=z^2+c_1$ and $f_{ c_2}(z)=z^2+ c_2$
with locally-connected Julia sets whose parameter values $c_1$ and $c_2$
belong to the conjugate limbs of the Mandelbrot set. In this case the rays $\{ R_1(t_j) \}$
and $\{ R_2 (t_j) \}$ landing at the dividing fixed points $\alpha_1$, $\alpha_2$ of the two polynomials
have opposite angles (see e.g. \cite{Milnor2}). This implies that $\alpha_1\rayeq\alpha_2$, and it is not hard to check that the quotient of $(\copyright_1\sqcup \copyright_2)/\infeq$ modulo $\rayeq$ is not homeomorphic to the $2$-sphere.

Recall that two branched coverings $F$ and $G$ of $S^2$ with finite postcritical sets $P_F$ and $P_G$
are equivalent {\it combinatorially} or {\it in the sense of Thurston} if there exist two orientation
preserving homeomorphisms $\phi, \psi : S^2 \rightarrow S^2$, such that $\phi\circ F=G\circ \psi$,
and $\psi$ is isotopic to $\phi$ rel $P_F$.  Using Thurston's characterization of critically finite
rational maps as branched coverings of the sphere (see \cite{Douady-Hubbard}), Tan Lei \cite{Tan} and Rees \cite{Rees1}
established the following:\vspace{0.12in}\\
{\bf Theorem.} {\it Let $c_1$ and $c_2$ be two parameter values not in conjugate limbs of the Mandelbrot set
such that $f_{c_1}$ and $f_{c_2}$ are postcritically finite. Then the map $F$ is combinatorially equivalent
to a quadratic rational map, where $F$ is either the formal mating $f_{c_1}\fmate f_{c_2}$ or a certain degenerate form of it.}\vspace{0.12in}

Taking this line of investigation further, Rees \cite{Rees2} and Shishikura \cite{Shishikura} demonstrated:\vspace{0.12in}\\
{\bf Theorem.} {\it Under the assumptions of the previous theorem, $f_{c_1}$ and $f_{c_2}$ are
topologically mateable. Moreover, their conformal mating $f_{c_1}\mate f_{c_2}$
exists.}\vspace{0.12in}

The case where the critical points of $f_{c_i}$ are periodic was considered by Rees, the complementary case was done by Shishikura.
 Note, in particular, that when none of the critical points is periodic, the Julia sets are dendrites with no interior,
 which makes the result particularly striking. An example of this phenomenon is analyzed in detail in Milnor's recent
 paper \cite{Milnor4} in which he considers the self-mating $F= f_{c_{1/4}} \mate f_{c_{1/4}}$, where the quadratic polynomial
 $f_{c_{1/4}}$ is the landing point of the $1/4$- external ray of the Mandelbrot set. It is not hard to deduce
 that  $F$ is a Latt\`es map, its Julia set $J(F)=\ov{\CC}$ is obtained by pasting together two copies of the dendrite
  $J(f_{c_{1/4}})$.

The issue of topological mateability is usually settled using the following result of R. L. Moore \cite{Moore}. Recall that an equivalence relation $\sim$ on $S^2$ is {\it closed} if $x_n \to x$, $y_n \to y$ and $x_n \sim y_n$ implies $x\sim y$. \vspace{0.12in}\\
{\bf Theorem (Moore).} {\it Suppose that $\sim$ is a closed equivalence relation on the $2$-sphere $S^2$ such that every equivalence class is a compact connected non-separating proper subset of $S^2$. Then the quotient space $S^2/\sim$ is again homeomorphic to $S^2$.}\vspace{0.12in}

For the application at hand, the theorem is replaced by the following corollary
 (see for example  Proposition 4.4. of \cite{ST}):
\vspace{0.12in}\\
{\bf Corollary.} {\it Let $f_1$ and $f_2$ be two quadratic polynomials with locally-connected Julia sets,
 such that every class 
of the ray equivalence relation  $\rayeq$ is non-separating and
 contains at most $N$ external rays for a fixed $N>0$.
 Then $f_1$ and $f_2$ are topologically mateable.}
\vspace{0.12in}

By means of a standard quasiconformal surgery, the theorem of Rees and Shishikura can be extended
to any pair $f_{c_1}$, $f_{c_2}$ where $c_i$ belong to hyperbolic components $H_1$, $H_2$  of the Mandelbrot
set which do not belong to conjugate limbs.  Mating thus yields an isomorphism between
the product $H_1\times H_2$ and a hyperbolic component in the parameter
space of quadratic rational maps.  This isomorphism, however, does not necessarily extend as a continuous
maps to the product of closures $\ov H_1\times \ov H_2$, as was recently shown by A. Epstein \cite{Epstein}.

So far no example of conformal matings without using Thurston's theorem (that is going beyond postcritically finite/hyperbolic case) has appeared in the literature. However, Jiaqi Luo
in his dissertation \cite{Luo} has outlined a proof of the existence of conformal matings 
of Yoccoz polynomials with star-like polynomials (centers of hyperbolic components attached to the main cardioid of the Mandelbrot set). His approach consists of locating a
candidate rational map for the mating, and then using {\it Yoccoz puzzle partitions} and
complex bounds of Yoccoz to prove that this candidate rational map {\it is} a mating. A somewhat similar philosophy plays a role in this paper.

The question of constructing matings of polynomials
with  connected but non locally-connected Julia sets 
has been completely untouched. While there are definitions of mating which would carry over to 
non locally-connected case (such as approximate matings discussed in \cite{Milnor-remarks}, p. 54)
 no examples of such matings are known.

\subsection{Statement of the results.}
Consider an irrational number $0< \theta <1$ and the quadratic polynomial $z\mapsto e^{2\pi i\theta}z+z^2$ which has an indifferent fixed point with multiplier $e^{2\pi i\theta}$ at the origin. To make this polynomial monic, we conjugate it by an affine map of $\CC$ to put it in the normal form
\begin{equation}
\label{eqn:ft} 
\FT: z \mapsto z^2+c_{\theta}, \ \mbox{with}\ c_{\theta}=\frac{e^{2\pi i\theta}}{2} \left ( 1-\frac{e^{2\pi i\theta}}{2} \right ).
\end{equation}

\realfig{Q1}{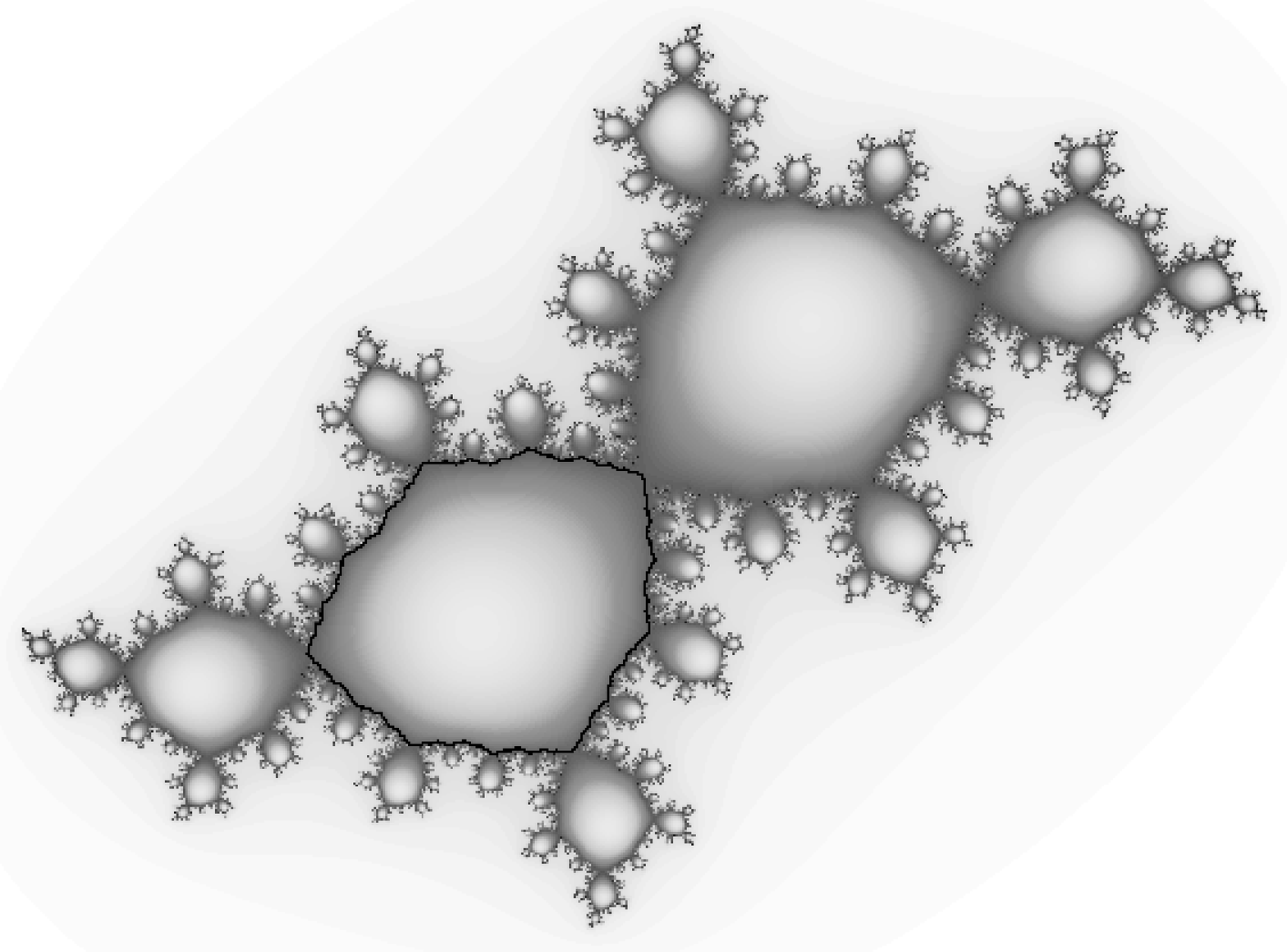}{{\sl Filled Julia set $K(\FT)$ 
for $\theta=(\sqrt{5}-1)/2$.}}{9cm}

\noindent
The corresponding indifferent fixed point of $\FT$ is denoted by $\alpha$.
Assuming $\theta$ is irrational of bounded type, a classical result of Siegel \cite{Carleson-Gamelin} implies that $\FT$ is linearizable near $\alpha$, i.e., there exists an open neighborhood $U$ of $\alpha$ and a conformal isomorphism $\phi: U \iso \DD$ which conjugates $\FT$ on $U$ to the rigid rotation $\varrho_\theta : z \mapsto e^{2\pi i \theta}z$:
$$\phi \circ \FT\circ \phi^{-1}=\varrho_\theta.$$
The maximal such linearization domain is a simply-connected neighborhood of $\alpha$ called the {\it Siegel disk} of $\FT$.
The following result has recently been proved by Petersen \cite{Petersen}:\vspace{0.12in}\\
{\bf Theorem (Petersen).} {\it Let $0< \theta < 1$ be an irrational of bounded type. Then the Julia
set of the quadratic polynomial $\FT$ is locally-connected and has Lebesgue measure zero.}\vspace{0.12in}\\
\figref{Q1} shows the filled Julia set of the quadratic polynomial $\FT$ for the golden mean $\theta=(\sqrt{5}-1)/2$. 

In proving his theorem, Petersen does not work directly with the Julia set of $\FT$, but
instead considers a certain \bp, which is related to $\FT$ via a quasiconformal surgery
procedure. A simplified version of his argument, based on complex {\it a priori} bounds
for renormalization of critical circle maps was presented by one of the authors in \cite{Yampolsky}.
Since the Julia set of $\FT$ is locally-connected, we may pose mateability questions for these polynomials. 
Our main result is the following theorem:\vspace{0.12in}\\
{\bf Main Theorem.} {\it Let $0< \theta, \nu <1$ be two irrationals of bounded type and $\theta\ne 1-\nu$.
Then the polynomials $\FT$ and $\FN$ are topologically mateable. Moreover, there exists a quadratic rational map $F$ such that
$$F=\FT\mate \FN.$$
Any two such rational maps are conjugate by a M\"obius transformation.}

\medskip

\noindent
In other words, one can paste any two filled Julia sets of the type shown in \figref{Q1} along their 
 boundaries to obtain a $2$-sphere, and the actions of the polynomials on their filled Julia sets match up to give an 
action on the sphere which is conjugate to a quadratic rational map with two fixed Siegel disks. \figref{blah} shows the 
result of this pasting in the case $\theta=\nu=(\sqrt{5}-1)/2$. In this picture we normalize the quadratic
 rational map $\FT \mate \FT$ to put the centers of the Siegel disks at zero and infinity. The black and gray regions are
 the images of the copies of the corresponding filled Julia sets in \figref{Q1}. There are, however, some prominent differences
 between these regions and the original filled Julia sets. First, there are infinitely many ``pinch points'' in the ``ends'' of
 the black and gray regions that are not present in the original filled Julia sets. 
An explicit combinatorial description of these pinch points will be presented in \secref{sec:remarks}.
Also, as J. Milnor pointed out to us, an infinite chain of preimages of the Siegel disk  in the filled
 Julia set in \figref{Q1} which lands at an endpoint in $J(\FT)$ maps to a chain in \figref{blah} which appears very stretched out near the end.
This indicates that the continuous semiconjugacies between the filled Julia sets and their
 corresponding regions, although  conformal in the interior of the sets, have a great amount of distortion near the boundary.
\realfig{blah}{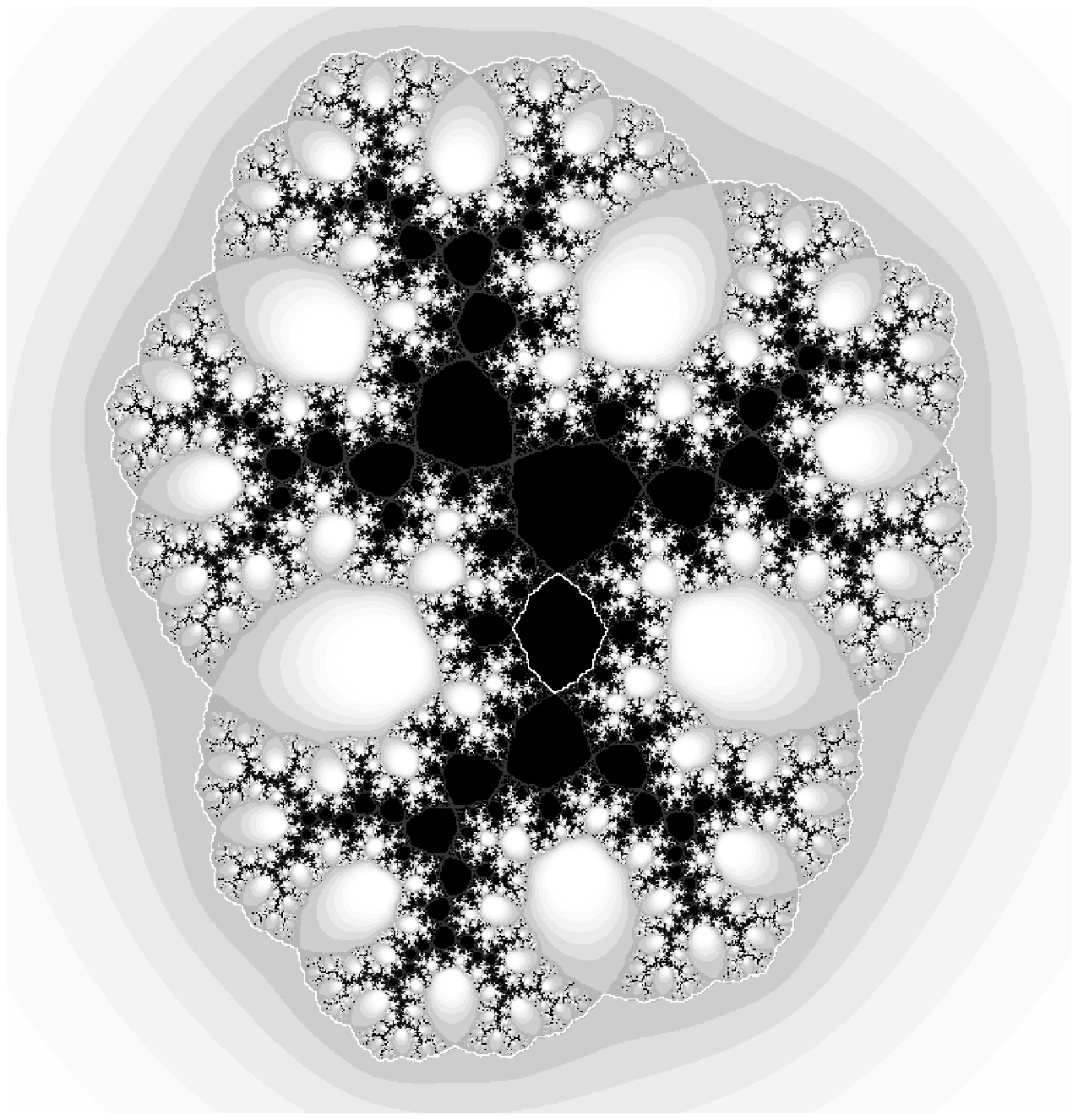}{{\sl The Julia set of the mating $\FT \mate \FT$ for $\theta=(\sqrt{5}-1)/2$.}}{9cm}

In the case $\theta=1-\nu$ the existence of a mating is ruled out for algebraic reasons. In fact, the polynomials are not even
 topologically mateable. Under the assumptions of the theorem, 
the candidate rational map $F$ can be specified algebraically, and the main difficulty
lies in establishing that $F$ is indeed a mating. 
To fix the ideas we may assume that the candidate $F$ has a Siegel disk $\Delta^0$ with rotation number $\theta$
centered at $0$, and another one $\Delta^\infty$  with rotation number $\nu$ centered at $\infty$. 
There is an unambiguous way to 
construct the semiconjugacies of Definition IIa in the interiors of the filled Julia sets,
by mapping the preimages of the Siegel disk of $\FT$ to the corresponding
preimages of $\Delta^0$ and similarly the preimages of the Siegel disk of $\FN$ to the corresponding preimages of $\Delta^\infty$.
To guarantee that these semiconjugacies extend continuously to the filled Julia sets we need
 to demonstrate that the boundaries $\bd \Delta^0$ and $\bd \Delta^\infty$ are Jordan curves each containing a critical point of $F$ and that the Euclidean diameter of the $n$-th preimages of $\Delta^0$ and $\Delta^\infty$ 
goes to zero uniformly in $n$. Proving these properties of the map $F$  directly seems to be quite out of reach. We establish the first property by using a new \bp model for the dynamics of
 $F$ that was discovered by one of the authors when he was working on dynamics of cubic Siegel polynomials
 \cite{Zakeri2}. We then adapt the complex bounds from \cite{Yampolsky} to this model to prove the second property. Further properties of the semiconjugacies of Definition IIa are demonstrated by a combinatorial argument using spines and itineraries.

The symmetry of the construction in the case of a self-mating (i.e., when $\theta=\nu$) has a nice corollary. 
In this case the mating $F=\FT \mate \FT$ given by the Main Theorem commutes with the M\"{o}bius involution $\cal I$
 which interchanges the centers of the two Siegel disks and fixes the third fixed point of $F$. Hence one can pass to 
the quotient Riemann surface $\ov{\CC}/ {\cal I}\simeq \ov{\CC}$ to obtain a new quadratic rational map $G$. It is not
 hard to see that $G$ is the mating of $\FT$ with the Chebyshev quadratic polynomial $\cheb: z\mapsto z^2-2$ whose filled 
Julia set is the interval $[-2,2]$: \vspace{0.12in}\\
{\bf Theorem.} {\it Let $0< \theta <1$ be any irrational of bounded type. Then there exists a quadratic rational map $G$ such that 
$$G= \FT \mate \cheb.$$
Moreover, $G$ is unique up to conjugation with a M\"obius transformation.}\vspace{0.12in}\\
{\bf Acknowledgements.} We would like to express our gratitude to John Milnor for posing the problem 
and encouraging the dynamics group at Stony Brook to look at it.
His picture of the ``presumed mating of golden ratio Siegel disk with itself''
 (\figref{blah} in this paper) posted in the IMS at Stony Brook was the inspiration for this work.
Adam Epstein, who also was enthusiastic about this problem and had learned about our similar ideas, 
brought the two of us together. We are indebted to him because this joint paper would have never existed 
without his persistence. Finally, we gratefully acknowledge the important role that Carsten Petersen's 
ideas in \cite{Petersen} play in our work.

\vspace{0.17in}
\section{Background Material}
\label{sec:prelim}

\subsection{Notations and terminology}
The unit disk in the complex plane will be denoted by $\DD$, its boundary is the unit circle $\TT$.
For a set $X$ in the plane, we use $\ov X$ and $\stackrel{\circ}{X}$ for the closure and the interior of $X$ respectively.
We use $|J|$ for the length of an interval $J$, 
 $\dist$ and $\diam$ for the Euclidean distance and diameter in $\CC$.
We write $[a,b]$  for  the closed interval with endpoints $a$ and $b$ in $\RR$ without
specifying their order. For a hyperbolic Riemann surface $X$,
$\dist_X$ will denote the distance in the hyperbolic metric in $X$.

We call two real numbers $a$ and $b$  {\it $K$-commensurable}
or simply {\it commensurable} if 
$K^{-1}\leq |a|/|b|\leq K$ for some $K>1$ independent of $a,b$. 
Two sets $X$ and $Y$ in $\CC$ are $K$-commensurable, if their
diameters are.
A configuration of points $x_1,\dots\,x_n$ is called {\it $K$-bounded}
if any two intervals $[x_i,x_j]$, and $[x_k,x_l]$ are $K$-commensurable.
For a pair of intervals $I\subset J$ we say that $I$ is
{\it well inside of } $J$ if there exists a universal constant $K>0$,
such that for each component $L$ of $J\setminus I$
we have $|L|\geq K|I|$.

For two points $a$, $b$ on the  circle which are not diagonally opposite $[a,b]$ will denote, unless otherwise
specified, the shorter
of the two closed arcs connecting them. 
When working with a homeomorphism $f$ of the unit circle, which extends beyond the  circle,
we will reserve the notation $f^{-i}(z)$ for the $i$-th preimage of $z\in\TT$
contained in the circle $\TT$.

\subsection{Quadratic rational maps}
\label{sec:rational-maps}
The reader may find a detailed discussion of the dynamics of quadratic rational maps
in Milnor's paper \cite{Milnor-remarks}. Below we give a brief summary of some relevant
facts. A quadratic rational map of the Riemann sphere $\ov\CC$ may be expressed as
a ratio
$$F(z)=\frac{a_0z^2+a_1z+a_2}{b_0z^2+b_1z+b_2}$$
with one of the coefficients $a_0$, $b_0$ different from $0$. The six-tuple
$(a_0:a_1:a_2:b_0:b_1:b_2)$ may be viewed as a point in the complex projective space
$ \CC\PP^5$. The space of all quadratic rational maps $\rat$ is identified
in this way with a Zariski open subset of $\CC\PP^5$ (see \cite{Milnor-remarks}
for a description of the topology of this set).
From the point of view of complex dynamics the quadratic rational
maps which are conjugate by a conformal isomorphism of the Riemann sphere are identified. 
That is, we consider the quotient space of $\rat$ by the action of the group   {\bf{M\"{o}b}} $\simeq PSL_2(\CC)$ of M\"obius transformations. This {\it moduli space}
of quadratic rational maps will be denoted ${\M}_2$. The action of {\bf{M\"{o}b}}
on $\rat$ is locally free, and the quotient space has the structure of a 2-dimensional
complex orbifold branched over a set ${\cal S}\subset {\M}_2$. This {\it symmetry locus} $\cal S$
consists of maps possessing a nontrivial automorphism group.

A more useful parametrization of the moduli space ${\M}_2$ comes from the following considerations.
Every map $F\in\rat$ has three not necessarily distinct fixed points. Let $\mu_1$, $\mu_2$, $\mu_3$
denote the multipliers of the fixed points. (By definition, the multiplier of $F$ at a fixed point $p$ is simply the derivative $F'(p)$ with appropriate modification if $p=\infty$.) Let
$$\sigma_1=\mu_1+\mu_2+\mu_3,\;\sigma_2=\mu_1\mu_2+\mu_1\mu_3+\mu_2\mu_3,\;\sigma_3=\mu_1\mu_2\mu_3$$
be the elementary symmetric functions of these multipliers. \vspace{0.12 in}\\
{\bf Proposition }(\cite{Milnor-remarks}, Lemma 3.1). {\it The numbers $\sigma_1$, $\sigma_2$, $\sigma_3$ determine 
$F$ up to a M\"{o}bius conjugacy, and are subject only to the restriction that
$$\sigma_3=\sigma_1-2.$$
Hence the moduli space ${\M}_2$ is canonically isomorphic to $\CC^2$, with coordinates $\sigma_1$ and $\sigma_2$.}
\vspace{0.12 in}

Note that for any choice of $\mu_1$, $\mu_2$ with $\mu_1\mu_2\ne 1$
there exists a quadratic rational map $F$, unique up to a M\"obius conjugacy, which has distinct fixed points
with these multipliers. The third multiplier can be computed as $\mu_3=(2-\mu_1-\mu_2)/(1-\mu_1\mu_2)$.

As a special case, let $F$ be a quadratic rational map which has two  Siegel disks centered at two fixed points of multipliers $e^{2 \pi i \theta}$ and $e^{2 \pi i \nu}$, where $0 < \theta, \nu <1$. Note that we necessarily have $\theta \neq 1-\nu$. By conjugating $F$ with a M\"{o}bius transformation which sends the two centers to $0$ and $\infty$ and the third fixed point to $1$, we obtain a quadratic rational map which fixes $0, 1, \infty$ and has multipliers $e^{2 \pi i \theta}$ at $0$ and $e^{2 \pi i \nu}$ at $\infty$. It is easy to see that these conditions determine the map uniquely. In fact, we obtain the normal form
\begin{equation}
\label{eqn:normal form} 
\FTN:z\mapsto z\ \frac{(1-e^{2 \pi i \theta})z+e^{2 \pi i \theta}(1-e^{2 \pi i \nu})}{(1-e^{2 \pi i \theta})e^{2 \pi i \nu} z+(1-e^{2 \pi i \nu})}.
\end{equation} 

\subsection{Critical circle maps.} 
\label{subsection:circle maps}
Throughout this paper, we shall identify the unit circle $\TT=\{ z\in \CC: |z|=1\}$ with the affine manifold $\RR/ \ZZ$ using the canonical projection from the real line given by $x\mapsto e^{2\pi i x}$. By definition, a {\it critical circle map } is an
 orientation-preserving homeomorphism of the circle $\TT$
 of class $C^3$ with a single critical point
 $c$. 
We further
assume that the critical point is of cubic type. This means that for a lift $\hat f:\RR \to \RR$
of $f$ with critical points at integer translates of $\hat c$,
$$\hat f(x)-\hat f(\hat c)=(x-\hat c)^3(\operatorname{const}+O(x-\hat c)).$$
The standard examples of analytic critical circle maps are provided by the projections to $\TT$ of
homeomorphisms in the {\it Arnold family}:
$$A^t: x \mapsto x+t-\frac{1}{2\pi}\sin 2\pi x.$$
Another group of examples, more relevant for our considerations, is given by the family of degree $3$ Blaschke products
$$Q^t: z \mapsto e^{2\pi i t}z^2 \left ( \frac{z-3}{1-3z} \right ).$$
The restriction of $Q^t$ to the unit circle $\TT$
is a real-analytic homeomorphism. Every $Q^t$ has a critical point of cubic type at $1\in\TT$ and no 
other critical points in $\TT$, thus  $Q^t|_\TT$ is a critical circle map.

The quantity 
$$\rho(f)= \lim_{n\to\infty}\frac{\hat f^{\circ n}(x)}{n}\ (\mod 1)$$
 is independent both of the choice of $x\in \RR$ and the lift
$\hat f$ of a critical circle map $f$, and is referred to as the {\it rotation number} of $f$.
The rotation number is rational of the form $\rho(f)=p/q$ if and only if $f$ has an orbit of period $q$. 
To further illustrate the connection between the number-theoretic properties of $\rho(f)$ and
the dynamics of $f$, let us introduce the notion of a closest return of the critical point $c$.
The iterate $f^{\circ n}(c)$ is a {\it closest return}, or equivalently, $n$ is a {\it closest return moment}, if the interior of the arc $[f^{\circ n}(c),c]$ contains no iterates $f^{\circ j}(c)$ with $j<n$. Consider the representation of $\rho(f)$ as a (possibly finite) continued fraction
$$\rho(f)=\cfrac{1}{a_1+\cfrac{1}{a_2+\cfrac{1}{a_3+\dotsb}}}\ ,$$
with the $a_i$ being positive integers.
For convenience we will write $\rho(f)=[a_1,a_2,a_3,\ldots]$.
The $n$-{\it th convergent} of the continued fraction of $\rho(f)$ is the rational number
$$\frac{p_n}{q_n}=[a_1,a_2,\ldots,a_n]$$
written in the reduced form. We set $p_0=0$, $q_0=1$. One easily verifies the recursive relations
$$p_n=a_np_{n-1}+p_{n-2},$$
$$q_n=a_nq_{n-1}+q_{n-2},$$
for $n\geq 2$. In this notation, the iterates $\{f^{\circ q_n}(c)\}$ are the 
consecutive closest returns of the critical point $c$ (see for example \cite{deMelo-vanStrien}).

The rotation number $\rho(f)$ is said to be of {\it bounded type} if $\operatorname{sup}a_i<\infty$.
We will make use of two linearization theorems for critical circle maps. Let us denote
by $\varrho_\theta$ the rigid rotation $x \mapsto x+\theta\ (\mod  \ZZ)$.
Yoccoz \cite{Yoccoz1} has shown:\vspace{0.12 in}\\
{\bf Theorem.} {\it Let $f$ be a critical circle map with irrational rotation number $\theta$. Then
there exists a homeomorphic change of coordinates $h:\TT \to \TT$ such that}
$$h \circ f\circ h^{-1}=\varrho_{\theta}.$$

In general the homeomorphism $h$ may not be regular at all, even if the map $f$ is real-analytic. 
However, some regularity for $h$ may be gained at the  expense of extra assumptions on the rotation number $\rho(f)$.
The following theorem of Herman \cite{Herman2} provides us with a sharp result which will be useful further in performing a
quasiconformal surgery. Recall that a homeomorphism $h:{\Bbb R} \rightarrow {\Bbb R}$ is called $K$-{\it quasisymmetric} if 
$$0< K^{-1} \leq \frac{|h(x+t)-h(x)|}{|h(x)-h(x-t)|}\leq K < +\infty$$
for all $x$ and all $t>0$. A homeomorphism $h:\TT \rightarrow \TT$ is $K$-quasisymmetric if its lift to ${\Bbb R}$ is such a homeomorphism. We simply call $h$ quasisymmetric if it is $K$-quasisymmetric for some $K$.\vspace{0.12 in}\\
{\bf Theorem.}{ \it A critical circle map $f$ is conjugate to a rigid rotation by a quasisymmetric
homeomorphism $h$ if and only if the rotation number $\rho(f)$ is irrational of bounded type.
}
\vspace{0.12 in}

The above result is based on the following a priori estimates called the
{\it \'Swia\c\negthinspace tek-Herman real a priori bounds} (see \cite{Swiatek},\cite{dFdM}):\vspace{0.12 in}\\
{\bf Theorem.} {\it Let $f$ be a critical circle map with irrational rotation number. Let $I_n$ denote the {\it $n$-th closest return interval} $[c,f^{\circ q_n}(c)]$. Then there exists $N=N(f)>0$ such that 
$$K^{-1}|I_n|\leq |I_{n+1}|\leq K|I_n|$$
for $n\geq N$ and a universal constant $K>1$. Moreover, let $\alpha_n: \RR \to \RR$ denote the affine map which restricts to a map $I_{n-1} \to [0,1]$ sending $c$ to $0$,
and set $q(z)=z^3$. Then, there exists a $C^2$-compact family $\cal F$ of 
$C^3$ 
diffeomorphisms of the interval $[0,1]$ into $\RR$ such that for $n>N$, 
$$\alpha_n \circ f^{\circ q_n}\circ \alpha_n ^{-1}|_{[0,1]}=H_n\circ q\circ h_n,$$
where $H_n\in \cal F$ and $h_n$ is a $C^3$ diffeomorphism of $[0,1]$ with $h_n \to \operatorname{id}$ in $C^2$-topology. }

\medskip
We conclude this section with a useful observation on the combinatorics of closest returns. Let the continued
fraction expansion $[a_1,a_2,\ldots]$ of the rotation number $\rho(f)$ of a critical circle map $f$ contain 
at least $n+1$ terms. Then (see \cite{deMelo-vanStrien}) for any $i\leq n$, 
the consecutive closest returns $f^{\circ q_i}(c)$ and $f^{\circ q_{i+1}}(c)$ occur on different sides of the critical point $c$,
that is $[f^{\circ q_i}(c),f^{\circ q_{i+1}}(c)]\ni c$. Let us list some of the points in the forward orbit of $c$ in the order
they are encountered when going from $f^{\circ q_{i-1}}(c)$ to $f^{q_i}(c)$:
$$f^{\circ q_{i-1}}(c),f^{\circ q_{i-1}+q_i}(c),f^{\circ q_{i-1}+2q_i}(c),\ldots,
f^{\circ q_{i-1}+a_{i+1}q_i}(c)=f^{\circ q_{i+1}}(c),c,f^{-q_{i+1}}(c),f^{\circ q_i}(c).$$
When $\rho(f)$ is irrational, \'Swia\c\negthinspace tek-Herman real a priori bounds imply that
for every $N>0$ there exists a universal constant $K_N$ such that the following holds.
For all sufficiently large $i$, the arcs
$[f^{\circ q_{i-1}+(j-1)q_i}(c),f^{\circ q_{i-1}+jq_i}(c)]$,
$[f^{-(j-1)q_i}(c),f^{-j q_i}(c)]$ and $[c,f^{\circ q_{i-1}}(c)]$ are $K_N$-commensurable, for $1\leq j\leq a_{i+1}-1$ with 
$\min(j,a_{i+1}-j)<N$.

\section{The Blaschke Model For Petersen's Theorem}
\label{sec:pet}

As a motivation for further discussion,  we present with slight modifications the construction of a model Blaschke product
for a Siegel quadratic polynomial used by Petersen in \cite{Petersen}. 
Much of the tools developed in this section will carry over to the 
Blaschke product model for mating introduced in \secref{sec:model}. It is somewhat easier, however,
to discuss them in this context.
Let us define 
\begin{equation}
\label{eqn:peter}
Q^t:z\mapsto e^{2 \pi i t} z^2 \left( \frac{z-3}{1-3z} \right).
\end{equation}
As we have seen in the previous section, the restriction $Q^t|_{\TT}$
is a critical circle map with critical value $t\in \TT$. The standard monotonicity considerations
imply that for each irrational number $0< \theta < 1$ there exists a unique value $t(\theta)$
for which the rotation number $\rho(Q^{t(\theta)}|_{\TT})=\theta$. Let us set $\QT = Q^{t(\theta)}$.
\comm{
\begin{equation}
\label{eqn:peter}
Q=\QT:z\mapsto e^{2 \pi i t(\theta)} z^2 \left( \frac{z-3}{1-3z} \right),
\end{equation}
where $0<t(\theta)<1$ is the unique angle for which the rotation number of $Q$ restricted to the unit circle $\Bbb T$ is a given irrational number $0<\theta <1$. This Blaschke product was introduced by Douady, Ghys, Herman and 
Shishikura as a model for the quadratic polynomial
$$\FT:z\mapsto e^{2 \pi i \theta}z + z^2,$$
when $\theta$ is an irrational number of bounded type \cite{Douady}. It has been used by Petersen to prove that the Julia set of $\FT$ is locally-connected and has Lebesgue measure zero \cite{Petersen} (see also \cite{Yampolsky} for a shorter proof based on the complex a priori bounds for critical circle maps). Following Petersen, we first give a symbolic description of various ``drops'' associated with $Q$ when $\theta$ is any irrational number. This turns out to be useful later when we construct topological semiconjugacies
from the Julia set $J(\FT)$ into the Julia set of the quadratic rational map $\FTN$ (see ??). We define wakes, limbs and drop-chains for $Q$. Much of these constructions will be carried over to the new Blaschke product model for mating we introduce in section ??, but we mention them in this section as they are easier to construct for $Q$. Finally, we sketch how to perform surgery on $\QT$ to obtain the quadratic polynomial $\FT$ when $\theta$ is of bounded type.
}
\subsection{Elementary properties.}
\label{subsec:Elementary properties}
For the moment, let us work with a fixed irrational $\theta$ and abbreviate $Q=Q_\theta$.
As seen from (\ref{eqn:peter}), $Q$ has superattracting fixed points at $0$ and $\infty$ and a double critical point at $z=1$.
 The immediate basin of attraction of infinity, which we denote by $A(\infty)$, is a simply-connected region on which $Q$ acts as a degree 2 branched covering. $Q$ commutes with the reflection ${\cal T}:z\mapsto 1/\ov{z}$ through $\TT$, so we have a similar description for $A(0)={\cal T}(A(\infty))$, the immediate basin of attraction of the origin.

Just as in the polynomial case, there exists a unique conformal isomorphism $\varphi: A(\infty) \iso \ov{\CC} \sm \ov{\DD}$ with $\varphi(\infty)=\infty$ and 
$ \varphi'(\infty)=1$, which conjugates $\varphi$ on $A(\infty)$ to the squaring map $z\mapsto z^2$ on $\ov{\CC} \sm \ov{\DD}$. We may use it to define the {\it external rays} $R^e(t)=\varphi^{-1} \{ re^{2 \pi i t}: r>1 \}$ 
for $t\in \TT$, and the {\it equipotentials} $E_r=\varphi^{-1} \{ re^{2 \pi i t}: t\in \TT \}$ for $r>1$. 
The ray $R^e(t)$ {\it lands} at $p$ if $\lim_{r \rightarrow 1} \varphi^{-1}(re^{2 \pi i t})=p$.  

\begin{proposition}
\label{basin}
$A(\infty)=\ov{\CC} \sm \ov{\bigcup_{n \geq 0}Q^{-n}(\DD)}$.
\end{proposition}

\begin{pf}
Let us put $U=\ov{\CC} \sm \ov{\bigcup_{n \geq 0}Q^{-n}(\DD)}$. Clearly $A(\infty) \subset U$ and $f(U)\subset U$. Since $\ov{\bigcup_{n \geq 0}Q^{-n}(\TT)} =J(Q)$, $U$ is a subset of the Fatou set of $Q$. Assume by way of contradiction that $A(\infty) \neq U$. Then there must be a connected component of $U$ other than $A(\infty)$ which eventually maps to a periodic Fatou component $V$ by Sullivan's No Wandering Theorem. We have $V\neq A(\infty)$, since otherwise $Q$ would have to have a pole $\neq \infty$ in $U$. According to Fatou-Sullivan, $V$ is either the attracting basin of an attracting or parabolic periodic point, or a Siegel disk or a Herman ring. In the first two cases, there must be a critical point in $V$ which converges to the periodic orbit. But $V\subset \CC \sm \ov{\DD}$ and there is no critical point of $Q$ in $\CC \sm \ov{\DD}$. In the last two cases, some critical point in $J(Q)$ must accumulate on the boundary of the Siegel disk or Herman ring. The only critical point in $J(Q)$ is $z=1$ whose forward orbit is dense on the unit circle $\TT$. It follows that $\TT$ must be the boundary of the Siegel disk or a component of the boundary of the Herman ring. Evidently this is impossible since $\TT$ is accumulated from both sides by points in $J(Q)$ near the critical point $z=1$.
\end{pf} 
\realfig{juliaqt}{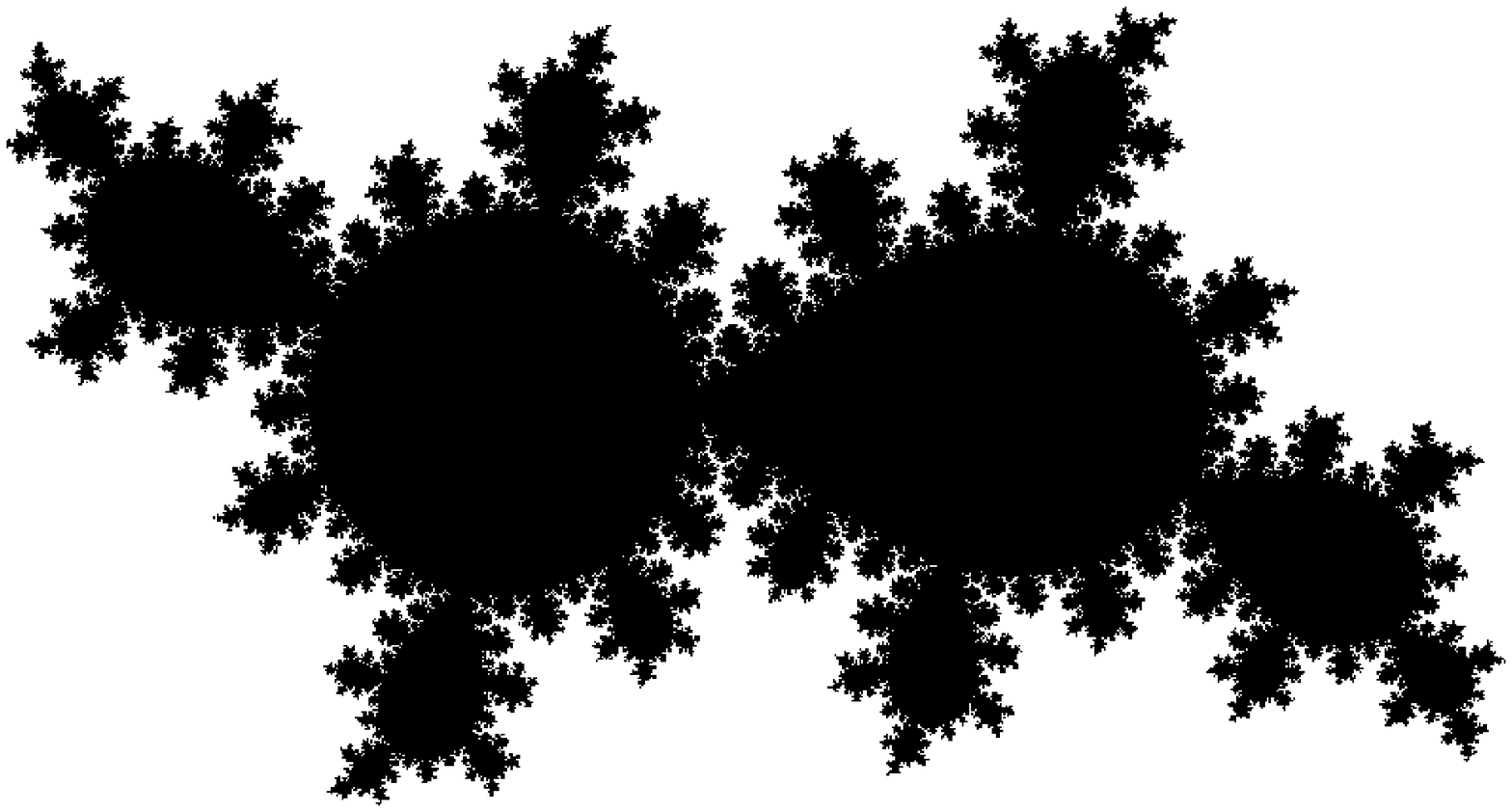}{{\sl ``Filled Julia set'' $K(\TQT)$ 
for $\theta=(\sqrt{5}-1)/2$.}}{8cm}

By the theorem of Yoccoz (see subsection \ref{subsection:circle maps}), there exists a unique homeomorphism $h:\TT \rightarrow \TT$ with $h(1)=1$ 
such that $h\circ Q|_\TT=\varrho_{\theta}\circ h$, where $\varrho_{\theta}: z\mapsto e^{2\pi i \theta}z$ is the rigid rotation 
by angle $\theta$. Let $H:\DD \rightarrow \DD$ be a homeomorphic extension of $h$ to the unit disk. To have 
a canonical homeomorphism at hand, we assume that $H$ is given by the Douady-Earle extension of circle homeomorphisms
 \cite{Douady-Earle}. Define a modified Blaschke product
\begin{equation}
\label{eqn:modify}
\tilde{Q}(z)=\TQT(z)= \left \{
\begin{array}{ll}
Q(z) & |z|\geq 1 \\
(H^{-1}\circ \varrho_{\theta} \circ H)(z) & |z| \leq 1
\end{array}
\right.
\end{equation}
where the two definitions match along the boundary of $\DD$.
Evidently, $\tilde{Q}$ is a degree $2$ branched covering of the sphere which is holomorphic outside of the unit 
disk and is topologically conjugate to a rigid rotation on the unit disk. 
Imitating the polynomial case, we define the ``filled Julia set'' of $\tilde{Q}$ by 
$$K(\tilde{Q})= \{ z \in \CC : \mbox{The orbit $\{ \tilde{Q}^{\circ n}(z) \}_{n \geq 0}$ is bounded} \} $$
and the ``Julia set'' of $\tilde{Q}$ as the topological boundary of $K(\tilde{Q})$:
$$J(\tilde{Q})=\partial  K(\tilde{Q}).$$
By \propref{basin}, we have 
$$K(\tilde{Q})=\ov{\CC} \sm A(\infty), \hspace{0.1in} J(\tilde{Q})=\partial A(\infty).$$  
In particular, $K(\tl Q)$ is full. \figref{juliaqt} shows the set $K(\tl Q)$ for the golden mean
$\theta=(\sqrt{5}-1)/2$; In this case, $t(\theta)=0.613648\ldots$.

\subsection{Drops and their addresses.}
\label{subsec:Drops and their addresses}
In what follows we collect basic facts about the ``drops'' associated with $\tilde Q$ and their addresses (see \cite{Petersen}, and
 compare \cite{Zakeri2} for a more general notion of a drop in a similar family of degree $5$ Blaschke products). By definition,
 the unit disk $\DD$ is called the $0$-{\it drop} of $\tilde Q$. For $n\geq 1$, any component $U$ of $\tilde{Q}^{-n}(\DD) \sm \DD$ is a 
Jordan domain called an $n$-{\it drop}, with $n$ being the {\it depth} of $U$. The map
 $\tilde{Q}^{\circ n}=Q^{\circ n}:U \rightarrow \DD$ is a conformal isomorphism. The unique point $z=z(U) \in U$ with the 
property $\tilde{Q}^{\circ n}(z)=H^{-1}(0)$ is called the {\it center} of $U$. This is the point in $U$ which eventually maps
 to the fixed point of the topological rotation on $\tilde{Q}:\DD\rightarrow \DD$. The unique point
 $\tilde{Q}^{-n}(1)\cap \partial U$ is called the {\it root} of $U$ and is denoted by $x(U)$. The 
boundary $\partial U$ is a real-analytic Jordan curve except at the root where it has a definite angle $\pi / 3$. 
We simply refer to $U$ as a drop when the depth is not important for us. Note that there is a unique $1$-drop $U_1$ 
which is the large Jordan domain attached to the unit disk at its root $x=1$ (see \figref{juliaqt}). 

Let $U$ and $V$ be two drops of depths $m$ and $n$ respectively. Then either 
$\ov{U}\cap \ov{V}=\emptyset$, or else $\ov{U}$ and $\ov{V}$ intersect at a unique point, in which case we necessarily have $m\neq n$. If we assume for example that $m<n$, then it is easy to check that $\ov{U}\cap \ov{V}=x(V)$. When this is the case, we call $U$ the {\it parent} of $V$, or $V$ a {\it child} of $U$. It is not hard to check that every $n$-drop with $n\geq 1$ has a unique parent which is an $m$-drop with $0 \leq m <n$. In particular the root of this $n$-drop belongs to the boundary of its parent. 

By definition, $\DD$ is said to be of {\it generation} $0$. Any child of $\DD$ is of generation $1$. In general, a drop is of generation $k$ if and only if its parent is of generation $k-1$. 

\begin{lemma}[Roots determine children]
\label{roots}
Given a point $p\in \bigcup_{n\geq 0} \tilde{Q}^{-n}(1) \sm \DD$, there exists a unique drop $U$ with $x(U)=p$. In particular, two distinct children of a parent have distinct roots.
\end{lemma}

\begin{pf}
It suffices to show that $U_1$ is the only child of $\DD$ whose root is $z=1$. Suppose that $U\neq U_1$ is an $n$-drop with $x(U)=1$. Then $\tilde{Q}^{\circ n-1}(U)=U_1$ implies $\tilde{Q}^{\circ n-1}(x(U))=x(U_1)$, or $\tilde{Q}^{\circ n-1}(1)=1$. Since $n> 1$ by the assumption, this contradicts the fact that the rotation number of $\tilde{Q}|_\TT=Q|_\TT$ is irrational. 
\end{pf} 

We give a symbolic description of various drops by assigning an address to every drop. 
This is a slightly modified version of Petersen's approach, based on a suggestion of J. Milnor. 
Set $U_0=\DD$. For $n\geq 1$, let $x_n=\tilde{Q}^{-n+1}(1) \cap \TT$ and $U_n$ be the $n$-drop with root $x_n$, 
which is well-defined by \lemref{roots}. Now let $\io=\io_1 \io_2 \cdots \io_k$ be any multi-index of length $k$, 
where each $\io_j$ is a positive integer. 
We inductively define the $(\io_1+\io_2+\cdots +\io_k)$-drop $U_{\io_1 \io_2 \cdots \io_k}$ of generation $k$ with root 
\begin{equation}
\label{eqn:root1}
x(U_{\io_1 \io_2 \cdots \io_k})=x_{\io_1 \io_2 \cdots \io_k}
\end{equation}
as follows. We have already defined these for $k=1$. For the induction step, suppose that we have defined $x_{\io_1 \io_2 \cdots \io_{k-1}}$ for all multi-indices $\io_1 \io_2 \cdots \io_{k-1}$ of length $k-1$. Then, we define
\begin{equation}
\label{eqn:rootx}
x_{\io_1 \io_2 \cdots \io_k}= \left \{
\begin{array}{ll}
\tilde{Q}^{-1}(x_{(\io_1 -1) \io_2 \cdots \io_k}) \cap \partial U_{\io_1 \io_2 \cdots \io_{k-1}} & \mbox{if}\ \io_1>1 \\
\tilde{Q}^{-1}(x_{\io_2 \cdots \io_k}) \cap \partial U_{\io_1 \io_2 \cdots \io_{k-1}} & \mbox{if}\ \io_1=1
\end{array}
\right.
\end{equation}
The drop $U_{\io_1 \io_2 \cdots \io_k}$ will then be determined by (\ref{eqn:root1}) and \lemref{roots} (see \figref{drops}).
\realfig{drops}{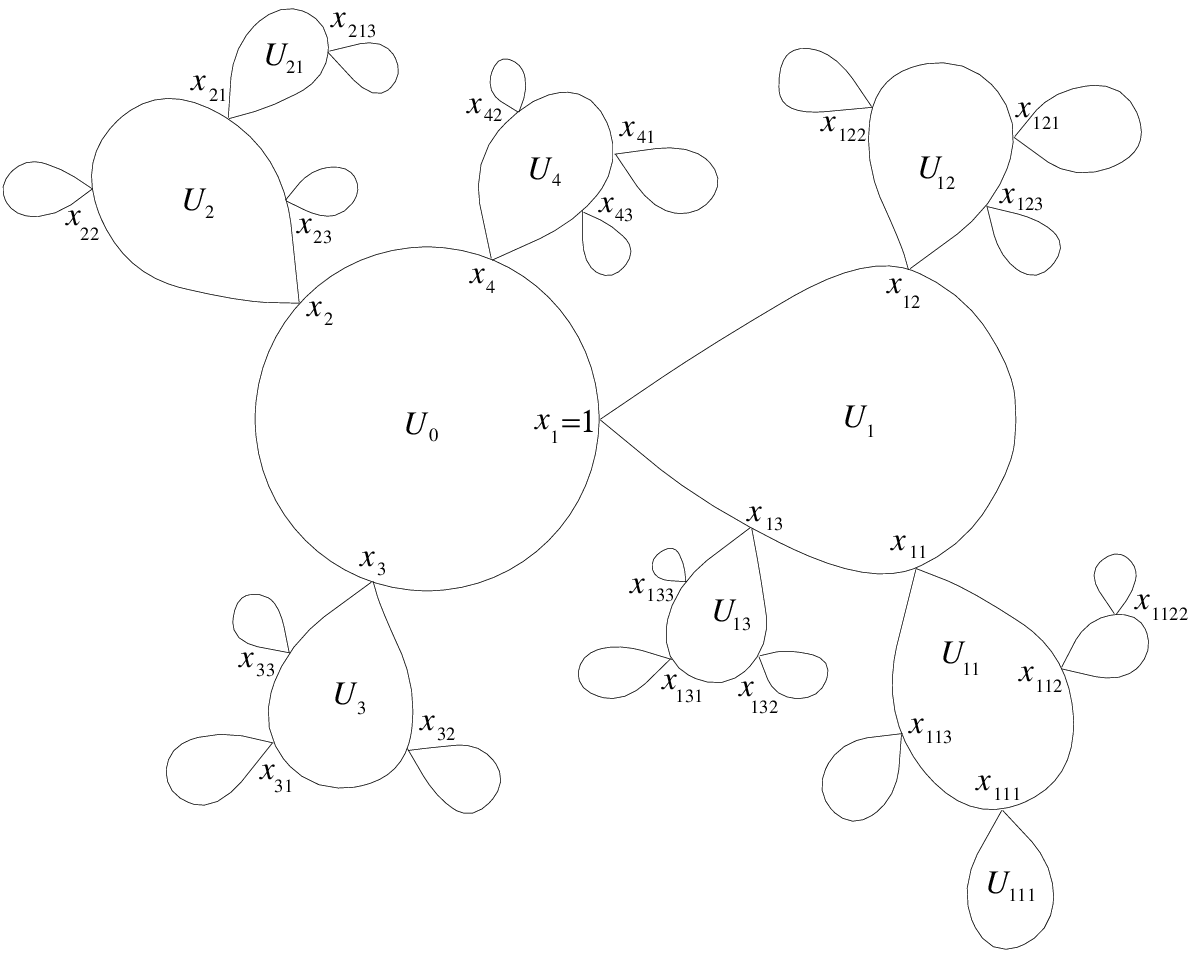}{{\sl Examples of some drops and their addresses.}}{10cm}

By the way these drops are given addresses, we have
\begin{equation}
\label{eqn:dynamics} 
\tilde{Q}(U_{\io_1 \io_2 \cdots \io_k})= \left \{
\begin{array}{ll}
U_{(\io_1 -1) \io_2 \cdots \io_k} & \mbox{if}\ \io_1>1 \\
U_{\io_2 \cdots \io_k} & \mbox{if}\ \io_1=1
\end{array}
\right.
\end{equation}

\subsection{Limbs and wakes.}
\label{subsec:Limbs and wakes}
Let us fix a drop $U_{\io_1 \cdots \io_k}$. By definition, the {\it limb} $L_{\io_1 \cdots \io_k}$ is the closure of the union of this drop and all its descendants (i.e., children and grand children etc.):
$$L_{\io_1 \cdots \io_k}=\ov{\bigcup U_{\io_1 \cdots \io_k\cdots}}\ .$$
Note that $L_0=K(\tilde Q)$. If $\io_1 \cdots \io_k \neq 0$, we call $x_{\io_1 \cdots \io_k}$ the root of $L_{\io_1 \cdots \io_k}$.
 
It is not immediatley clear from this definition that limbs provide a useful partition of the filled Julia set
$K(\tilde{Q})$. Indeed, it may happen {\it a priori} that the boundary of a limb$\neq L_0$ is the whole $J(\tl Q)$.
This is ruled out by the  following key lemma of Petersen \cite{Petersen}:

\begin{lemma}[Only two rays]
\label{tworays}
Suppose that $0< \theta <1$ is an irrational number. Then the critical point $z=1$ of $\QT$ is the landing point of two and only two external rays $R^e(t)$ and $R^e(s)$ in $A(\infty)$. 
\end{lemma}

Let $W_1$ denote the connected component of $\CC \sm (R^e(t) \cup R^e(s) \cup \{ 1\} )$ containing the drop $U_1$.
 We call $W_1$ the {\it wake} with root $x_1$. 
Given an arbitrary multi-index $\io_1 \cdots \io_k$, we define the wake $W_{\io_1 \cdots \io_k}$ as the appropriate pull-back
 of $W_1$. More precisely, consider the two external rays landing at $x_{\io_1 \cdots \io_k}$ which map to $R^e(t)$ and 
$R^e(t)$ under $\tilde{Q}^{\circ n}$, where $n=\io_1 + \cdots +\io_k$. These rays separate the plane into two simply-connected regions.
 The wake $W_{\io_1 \cdots \io_k}$ will then be the region containing the drop $U_{\io_1 \cdots \io_k}$.
It is immediately clear that
$$L_{\io_1 \cdots \io_k}=\ov{W}_{\io_1 \cdots \io_k} \cap K(\tilde{Q})$$
(see \figref{tail}). The integers $n$ and $k$ are respectively called the depth and generation of $W_{\io_1 \cdots \io_k}$ as well as $L_{\io_1 \cdots \io_k}$.
\realfig{tail}{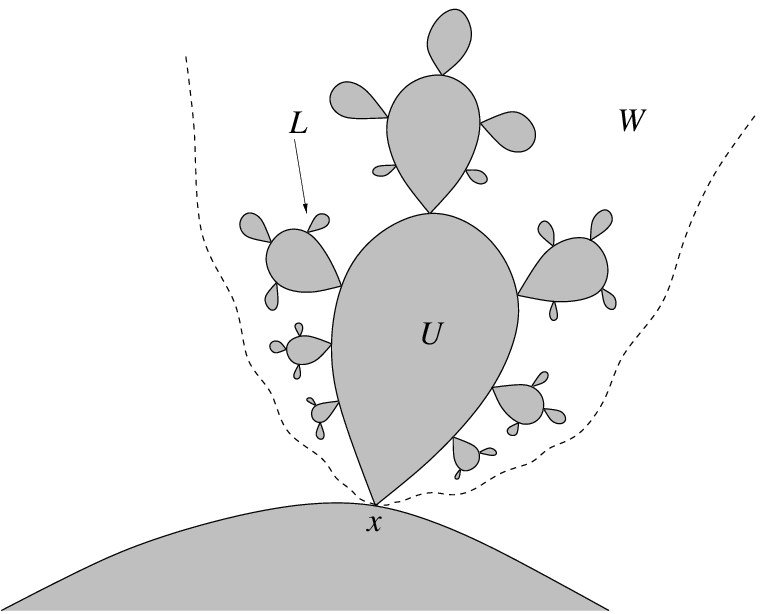}{{\sl A drop $U$ with root $x$, and the associated limb $L$ and wake $W$.}}{6cm}

\noindent
The next proposition follows directly from the above definitions:

\begin{proposition}[Properties of limbs and wakes]
\label{wake}
Consider $\TQT$ for an irrational number $0< \theta <1$. Then
\begin{enumerate}
\item[(i)]
If a drop $U$ is contained in a limb $L$, then any child of $U$ is also contained in $L$.
\item[(ii)]
Any two limbs and any two wakes are either disjoint or nested.  
\item[(iii)]
For any limb $L_{\io_1 \cdots \io_k}$, we have
$$\TQT(L_{\io_1 \cdots \io_k})= \left \{
\begin{array}{ll}
L_{(\io_1 -1) \io_2 \cdots \io_k} & \mbox{if}\ \io_1>1 \\
L_{\io_2 \cdots \io_k} & \mbox{if}\ \io_1=1
\end{array}
\right.
$$
In particular, every limb eventually maps to $L_1$ and then to the whole filled Julia set $K(\TQT)$. The same relation holds for wakes.
\end{enumerate}
\end{proposition}

The following theorem is a central result of \cite{Petersen}. 

\begin{theorem}[Local-connectivity]
\label{diam->0}
Suppose that $0< \theta <1$ is an irrational number. Then as the depth of a limb $L$ of $\TQT$ goes to infinity,
 $\diam (L) \rightarrow 0$. This implies that the Julia set $J(\QT)$, hence $J(\TQT)$, is locally-connected.
\end{theorem}

\noindent
In particular, it follows that the diameter of a drop goes to zero as the depth goes to infinity,
 simply because every drop is a subset of the limb with the same root.

One important implication of this result is the lack of the so-called ``ghost limbs'':

\begin{corollary}[No ghost limbs]
\label{ghost}
Suppose that $0< \theta <1$ is an irrational number. Then the filled Julia set  $K(\TQT)$ is the union of $\ov{\DD}$ and all
 the limbs of generation $1$:
$$K(\TQT)=\ov{\DD} \cup \bigcup_{n\geq 1}L_n.$$
\end{corollary}

\noindent
This follows from the fact that distinct $L_n$'s  are separated by their wakes and $\dia(L_n) \rightarrow 0$
 as $n \rightarrow \infty$. 

\subsection{Drop-chains.}
\label{subsec:Drop-chains}
\begin{definition} 
\label{dropchain}
Consider a sequence of drops $\{ U_0=\DD, U_{\io_1}, U_{\io_1 \io_2},  U_{\io_1 \io_2 \io_3}, \cdots \}$ where each $U_{\io_1 \cdots \io_k}$ is the parent of $U_{\io_1 \cdots \io_{k+1}}$. The closure of the union 
$$\CCC =\ov{\bigcup_k U_{\io_1 \cdots \io_k}}$$
is called a {\it drop-chain}.
\end{definition}

Since in a drop-chain $\CCC$ each parent touches the child at its root and the diameter of the subsequent children goes to zero by \thmref{diam->0}, the tail of $\CCC$ must converge to a well-defined point in the Julia set of $\tilde Q$. In other words, there exists a unique point $p=p(\CCC )$ such that in the Hausdorff topology, $\lim_{k\rightarrow \infty} \ov{U}_{\io_1 \cdots \io_k}=\{ p\}$. It follows that 
$$\CCC=\bigcup_k \ov{U}_{\io_1 \cdots \io_k} \cup \{ p \}.$$
In particular, $\CCC$ is compact, connected and locally-connected.

Another way to characterize $p(\CCC)$ is as follows: Consider the corresponding limbs 
$$K(\tilde Q)=L_0 \supset L_{\io_1} \supset L_{\io_1 \io_2} \supset L_{\io_1 \io_2 \io_3}\supset \cdots$$
which are nested by \propref{wake}. Since $\dia(L_{\io_1 \cdots \io_k})\rightarrow 0$ as $k \rightarrow \infty$ by \thmref{diam->0}, the intersection of these limbs must be a unique point, namely $p(\CCC)$:
$$p(\CCC)=\bigcap_k L_{\io_1 \cdots \io_k}.$$     

By a {\it ray} in a drop $U$ we mean a hyperbolic geodesic which connects some boundary point $p\in \partial U$ to the center $z(U)$. This ray is denoted by $[\![p, c(U)]\!]$. For two distinct points $p,q \in \partial U$, we use the notation $[\![p,q]\!]$ for the union of the rays $[\![p,c(U)]\!] \cup [\![c(U),q]\!]$. 

Given any drop-chain $\CCC$, there exists a unique ``most efficient'' path $R=R(\CCC)$ in $\CCC$ which connects $0$ to $p(\CCC)$. In fact, if $\CCC$ is of the form $\ov{\bigcup_k U_{\io_1 \cdots \io_k}}$, we define
$$R(\CCC)=[\![0,x_{\io_1}]\!] \cup \bigcup_{k\geq 2} [\![x_{\io_1 \cdots \io_k},x_{\io_1 \cdots \io_{k+1}}]\!] \cup \{ p(\CCC) \}.$$
(see \figref{rays}). It is easy to see that $R(\CCC)$ is a piecewise analytic embedded arc in the plane. We call $R(\CCC)$ the {\it drop-ray} associated with $\CCC$. We often say that $R(\CCC)$, or $\CCC$, {\it lands} at $p(\CCC)$.

\realfig{rays}{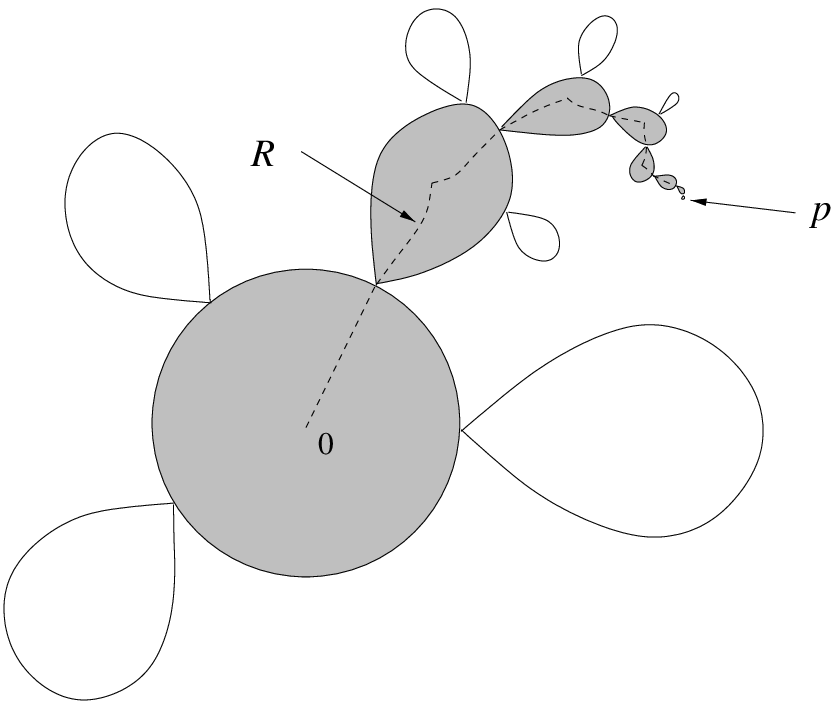}{{\sl A drop-chain 
and the drop-ray associated with it.}}{8cm}

\begin{proposition}
\label{either/or}
Every point in the filled Julia set $K(\TQT)$ either belongs to the closure of a drop or is the landing point of a unique drop-chain.
\end{proposition}

\begin{pf}
Let $p \in K(\TQT)$ and assume that $p$ does not belong to the closure of any drop. Then by \corref{ghost}, $p$ belongs to some limb $L_{\io_1}$, and inductively, it follows that it belongs to the intersection of a decreasing sequence of limbs $L_{\io_1} \supset L_{\io_1 \io_2} \supset L_{\io_1 \io_2 \io_3}\supset \cdots$. Hence $p$ is the landing point of the corresponding drop-chain $\CCC=\ov{\bigcup_k U_{\io_1 \cdots \io_k}}$. Uniqueness of this drop-chain follows from \propref{tree} below.
\end{pf}  

It follows from the next proposition that the union of drop-rays associated with all drop-chains has the structure of an infinite topological tree (a ``dendrite'') in the plane.

\begin{proposition}
\label{tree}
The assignment $\CCC \mapsto p(\CCC)$ is one-to-one. In other words, different drop-rays land at distinct points.
\end{proposition}

\begin{pf}
Suppose that $\CCC_1$ and $\CCC_2$ are two distinct drop-chains. Let $U_{\io_1 \cdots \io_k} \subset \CCC_1$ be the drop of smallest generation $k$ which is disjoint from $\CCC_2$, and similarly define $U_{\io'_1 \cdots \io'_k} \subset \CCC_2$. The limbs $L_{\io_1 \cdots \io_k}$ and $L_{\io'_1 \cdots \io'_k}$ are disjoint by \propref{wake}. Since $p(\CCC_1)\in L_{\io_1 \cdots \io_k}$ and $p(\CCC_2) \in L_{\io'_1 \cdots \io'_k}$, we will have $p(\CCC_1)\neq p(\CCC_2)$.
\end{pf}

\subsection{Surgery.}
\label{subsec:Surgery}
The modified Blaschke product $\tilde{Q}=\TQT$ as defined in (\ref{eqn:modify}) is a degree $2$ branched
 covering of the sphere.
When the rotation number $\theta$ is irrational of bounded type, the action of $\TQT$ is in fact
conjugate to that of a quadratic polynomial. This follows from a {\it quasiconformal surgery}
construction due to Douady, Ghys, Herman, and Shishikura \cite{Douady}. 

\comm{
This has been made possible by the following theorem of Swiatek and Herman \cite{Swiatek}, \cite{Herman2}:

\begin{theorem}[Linearization of Critical Circle Maps]
\label{ccmap}
Let $f:\TT \rightarrow \TT$ be a real-analytic homeomorphism with finitely many critical points and rotation number $\theta$. Then there exists a quasisymmetric homeomorphism 
$h:\TT \rightarrow \TT$ which conjugates $f$ to the rigid rotation $R_{\theta}$ if and only if $\theta$ is an irrational number of bounded type. Moreover, $h$ can be chosen to be $k(\theta)$-quasisymmetric, where $k(\theta)$ is a constant which only depends on the rotation number $\theta $ and not on the choice of $f$.
\end{theorem}

}
Let us fix an irrational number $0< \theta < 1$  of bounded type. By Herman's Theorem (see subsection \ref{subsection:circle maps})
the unique homeomorphism $h:\TT \rightarrow \TT$ with $h(1)=1$ which conjugates $Q|_{\TT}$ to $\varrho_{\theta}$ 
is  quasisymmetric. In this case, the Douady-Earle extension $H:\DD \rightarrow \DD$ of
 $h$ is a quasiconformal homeomorphism whose dilatation only depends on the dilatation of $h$ \cite{Douady-Earle}.
 The modified Blaschke product $\TQT$ of (\ref{eqn:modify}) is then a
 quasiregular branched covering of the sphere. We define a $\TQT$-invariant conformal structure $\sigma_{\theta}$ on the plane
 as follows: On $\DD$, let $\sigma_{\theta}$ be the pull-back $H^{\ast}\sigma_0$ of the standard conformal structure $\sigma_0$.
 Since $\varrho_\theta$ preserves $\sigma_0$, $\TQT$ will preserve $\sigma_{\theta}$ on $\DD$. For every $n\geq 1$, pull $\sigma_{\theta}|_{\DD}$
 back by $\tilde {\QT}^{\circ n}= \QT^{\circ n}$ on $\TQT^{-n}(\DD)\smallsetminus \DD$, which consists of all drops of $\QT$ 
of depth $n$. Since $\QT^{\circ n}$ is 
holomorphic, this does not increase the dilatation of $\sigma_{\theta}$. Finally, let 
$\sigma_{\theta}=\sigma_0$ on the rest of the plane. By  construction, $\sigma_{\theta}$
has bounded dilatation and is invariant under 
$\TQT$. Therefore, by the Measurable Riemann Mapping Theorem
(see for example \cite{Ahlfors-Bers}),
 we can find a unique quasiconformal homeomorphism $\PT:\ov{\CC} \rightarrow \ov{\CC}$,
 normalized by $\PT(\infty)= \infty $, $\PT(H^{-1}(0))=e^{2 \pi i \theta}/2$ and $\PT(1)=0$, such that $\PT^{\ast}\sigma_0=\sigma_{\theta}$. Set
\begin{equation}
\label{eqn:fmodb} 
\FT=\PT \circ \TQT \circ \PT^{-1}.
\end{equation}
Then $\FT$ is a quasiregular self-map of the sphere which preserves $\sigma_0$, hence it is holomorphic. 
Also $\FT:\CC \rightarrow \CC$ is a proper map of degree $2$ since $\TQT$ has the same properties. Therefore
$\FT$ is a quadratic polynomial.

Since the action of $\FT$ on $\PT (\DD)$ is quasiconformally conjugate to a rigid rotation, 
 $\PT (\DD)$ is contained in a Siegel disk for $\FT$ with rotation number $\theta$.
 As $\PT (1)=0$ is a critical point for $\FT$, it follows that the entire orbit 
$\{ \FT ^{\circ n}(0)\}_{n\geq 0}$ lies on the boundary of this Siegel disk.
 But $\{ \FT ^{\circ n}(0)\}_{n\geq 0}$ is dense on $\PT (\TT)$, so $\PT (\TT )$ 
is exactly the boundary of this Siegel disk, which is a quasicircle passing through the critical point $0$ of $\FT$. Up to affine conjugacy there is only one quadratic polynomial with a fixed 
Siegel disk of the given rotation number $\theta$. By the way we normalized $\PT$, we must have 
$\FT:z\mapsto z^2+c_{\theta}$ as in (\ref{eqn:ft}).

We summarize the above as follows:
\begin{theorem}[Douady, Ghys, Herman, Shishikura]
\label{dghs}
Let $f$ be a quadratic polynomial which has a fixed Siegel disk $\Delta$ of rotation number $\theta$.
 If $\theta$ is of bounded type, then $f$ is quasiconformally conjugate to $\TQT$ in $(\ref{eqn:modify})$. 
In particular, $\partial \Delta$ is a quasicircle passing through the critical point of $f$.
\end{theorem}

\noindent
In particular, this surgery procedure allows us to define drops, limbs, wakes, drop-chains and drop-rays for the quadratic polynomial $\FT$.   

\vspace{0.17in}

\section{A Blaschke Model For Mating}
\label{sec:model}

The object of this section is to construct, for a pair of numbers $0< \theta ,\nu <1$ with $\theta \neq 1-\nu$,
a Blaschke product $\BTN$. When $\theta$ and $\nu$ are irrationals of bounded type, $\BTN$
plays the role of a model for the quadratic rational map $\FTN$ of (\ref{eqn:normal form}) in the same way as $\QT$ does 
for the quadratic polynomial $\FT$. After showing the existence of such
$\BTN$, we will define drops, limbs, drop-chains and drop-rays for the ``modified'' $\TBTN$ in an analogous way.

\subsection{Existence.}
\label{subsec:Existence}
We would like to prove the following result:

\begin{theorem}[Existence of Blaschke models for mating]
\label{realize}
Let  $0\leq \theta <1$, $0\leq \nu <1$ and $\theta\ne 1-\nu$.
 Then there exists a degree $3$ Blaschke product
\begin{equation}
\label{eqn:BB}
B=\BTN: z\mapsto \frac{e^{-2\pi i\nu}}{ab}\ z\ \left ( \frac{z-a}{1-\ov{a}z} \right ) \left (
\frac{z-b}{1-\ov{b}z} \right )
\end{equation}
with the following properties:
\begin{enumerate}
\item[(i)]
$0<|a|<1$ and $|b|=|a|^{-1}>1$, with $a \ov{b} \neq 1$,
\item[(ii)]
$B$ has a double critical point at $z=1$, and
\item[(iii)]
The restriction $B|_{\Bbb T}$ is a critical circle map with rotation number $\theta$.
\end{enumerate}
\end{theorem}

The proof of this theorem will be given in the rest of this subsection. In (i) the condition $a\ov{b}\neq 1$
is necessary simply because when $a \ov{b}=1$, $B$ reduces to the linear map $z\mapsto e^{-2\pi i \nu}z$. 
  
For simplicity, let us set
\begin{equation}
\label{eqn:TZ1}
\begin{array}{l} 
\kappa =ab,\ \mbox{where $|\kappa|=1$ by (i)} \\
\zeta =a+b
\end{array}
\end{equation}
Using the equation (\ref{eqn:BB}), the condition $B'(z)=0$ may be written in the form
$$A_1 z^4+A_2 z^3+A_3 z^2+\ov{A}_2z+\ov{A}_1=0,$$ where
\begin{equation}
\label{eqn:ABC}
\begin{array}{l}
A_1 = \ov{a}\ \ov{b}=\ov{\kappa}, \\ A_2 =-2(\ov{a}+\ov{b})=-2\ov{\zeta}, \\ A_3 =2+|a+b|^2=2+|\zeta|^2.
\end{array}
\end{equation}
A brief computation  shows that the condition of $z=1$ being a double critical point of $B$ translates
into
$$ \left \{
\begin{array}{l}
4A_1 + 3 A_2 +2 A_3=-\ov{A}_2\\ 3A_1 + 2 A_2 + A_3 =\ov{A}_1
\end{array}
\right. $$
or by (\ref{eqn:ABC})
\begin{equation}
\label{eqn:TZ2}
\left \{ 
\begin{array}{l}
2\kappa -3\zeta +2+|\zeta|^2=\ov{\zeta}\\ 3\kappa -4\zeta +2+|\zeta|^2=\ov{\kappa}
\end{array}
\right.
\end{equation}
Subtracting the second equation in (\ref{eqn:TZ2}) from the first equation, we find that
$$\zeta -\kappa =\ov{\zeta}- \ov{\kappa} \Longrightarrow \zeta -\kappa\in {\Bbb R}.$$ Set $\kappa=x+iy$ and
$\zeta=u+iy$ and substitute them into the first equation in (\ref{eqn:TZ2}) to obtain
$$u^2-4u+(2x+y^2+2)=0,$$ which, by $x^2+y^2=1$, has solutions $u=x+1$ and $u=-x+3$. These correspond to
$\zeta=\kappa +1$ and $\zeta=-\ov{\kappa}+3$. By (\ref{eqn:TZ1}), the choice of $\zeta=\kappa +1$ leads to
$a=\kappa$ or $a=1$, which is not appropriate since we want $|a|<1$. Therefore, we are left with the only
possibility
\begin{equation}
\label{eqn:ZZ}
\zeta=-\ov{\kappa}+3.
\end{equation}
Let $\kappa=e^{2 \pi i t}$ with $t \in {\Bbb R}$. From (\ref{eqn:TZ1}) and (\ref{eqn:ZZ}) it follows that $a$
and $b$ are the solutions of the quadratic equation
\begin{equation}
\label{eqn:AB}
z^2+(\ov{\kappa}-3)z+\kappa=0.
\end{equation}

\begin{lemma}
\label{twobranch}
As $\kappa=e^{2 \pi i t}$ goes around the unit circle, the two solutions of the quadratic equation
(\ref{eqn:AB}) define two closed curves $t\mapsto a(t)$ and $t\mapsto b(t)$ in the complex plane with the
following properties (see \figref{szplot}):
\begin{enumerate}
\item[(i)]
$a(t+1)=a(t)$ and $b(t+1)=b(t)$,
\item[(ii)]
$0<|a(t)| \leq 1$ and hence $|b(t)|=|a(t)|^{-1} \geq 1$,
\item[(iii)] 
$|a(t)|=1$ if and only if $t\in {\Bbb Z}$, or equivalently $\kappa=1$, in which case $a(t)=b(t)=1$,
\item[(iv)]
$a(t)\ov{b(t)}\neq 1$ unless $t\in {\Bbb Z}$ so that $a(t)=b(t)=1$.
\end{enumerate}
\end{lemma}

\begin{pf}
Let us first note that the solutions $z_1,z_2$ of (\ref{eqn:AB}) lie on the unit circle ${\Bbb T}$ if and
only if $\kappa=1$ in which case there is a double root at $z_1=z_2=1$. In fact, if $|z_1|=|z_2|=1$, then
$$2= 3-|\ov{\kappa}| \leq |\ov{\kappa}-3|=|z_1+z_2| \leq |z_1|+|z_2|=2.$$ Hence $|\ov{\kappa}-3|=2$, or
equivalently, $\kappa=1$.

Now let $\kappa=e^{2 \pi i t}$ go around ${\Bbb T}$. Then the double root at $z=1$ splits into distinct roots
$a=a(t)$ and $b=b(t)$ which by inspecting the explicit formula for $a$ and $b$ are real-analytic functions
of $t$ away from integer values and are labeled so that (ii) holds. Clearly
$a$ and $b$ are ${\Bbb Z}$-periodic, so (i) holds trivially.
\realfig{szplot}{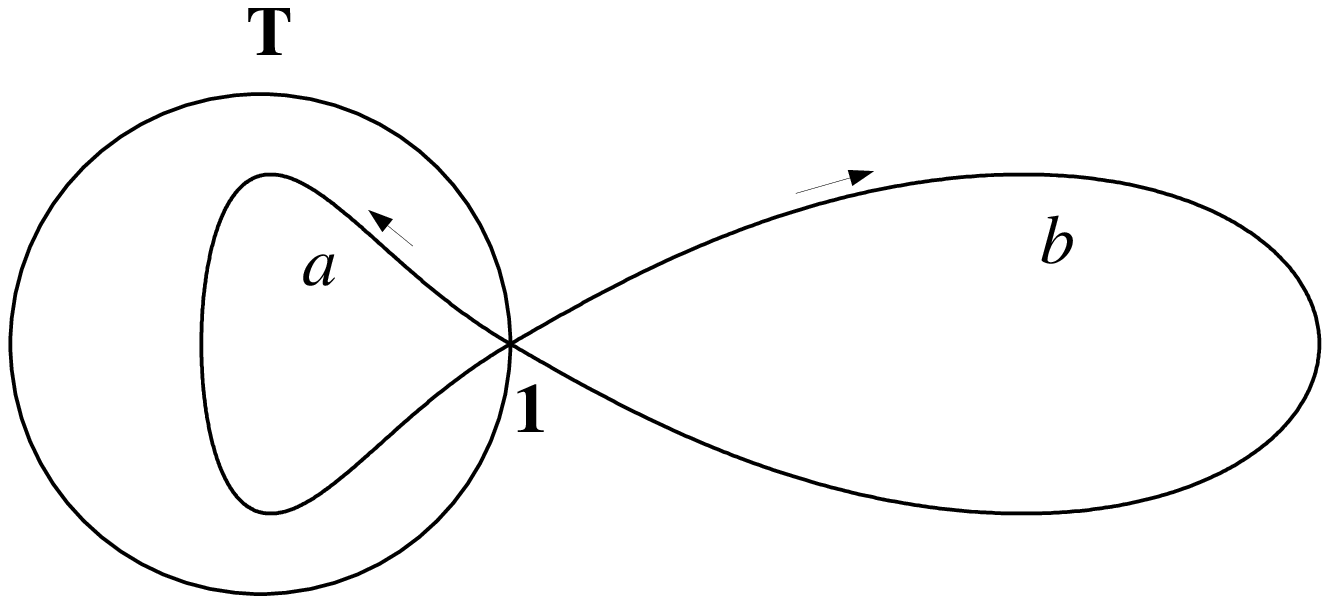}{}{7cm}

Finally, suppose that for some $t\in {\Bbb R}$, $a=a(t)$ and $b=b(t)$ satisfy $a\ov{b}=1$. Then
$a/\ov{a}=\kappa$, or $\ov{a}=a\ \ov{\kappa}$. Since $a$ is a solution of (\ref{eqn:AB}), we have
$$\ov{a}^2+(\kappa-3)\ov{a}+\ov{\kappa}=0 \Longrightarrow a^2 \ov{\kappa}^2+(\kappa-3)a\ov{\kappa}+\ov{\kappa}=0,$$ or,
after multiplying by $\kappa^2$,
\begin{equation}
\label{eqn:AT}
a^2+\kappa(\kappa-3)a+\kappa=0.
\end{equation}
Comparing (\ref{eqn:AT}) and (\ref{eqn:AB}) for $z=a$, we conclude that
$$\kappa(\kappa-3)=\ov{\kappa}-3 \Longrightarrow \kappa^2(\kappa-3)=1-3\kappa \Longrightarrow (\kappa-1)^3=0$$ which
shows $\kappa=1$.
\end{pf}

\begin{lemma}
\label{index}
For any $z\in {\Bbb T}$, the closed curve $\Gamma_z : [0,1]\rightarrow {\Bbb T}$ defined by
\begin{equation}
\label{eqn:gam} 
\Gamma_z(t)=\left ( \frac{z-a(t)}{1-\ov{a(t)}z} \right ) \left ( \frac{z-b(t)}{1-\ov{b(t)}z} \right )
\end{equation}
is null-homotopic.
\end{lemma}

Note that when $z=1$, there is no ambiguity in the definition of $\Gamma_z$. In fact, by (\ref{eqn:TZ1})
and (\ref{eqn:ZZ}),
$$\Gamma_1=\frac{1-\zeta+\kappa}{1-\ov{\zeta}+\ov{\kappa}}=\frac{-2+\kappa+\ov{\kappa}} {-2+\kappa+\ov{\kappa}}\equiv
1$$ so that $\Gamma_1$ is the constant loop $1$.

\begin{pf}
Consider the two homotopies $(t,s)\mapsto a(t,s)$ and $(t,s)\mapsto b(t,s)$ rel $\{ 1 \}$ defined by
$$a(t,s)=(1-s)a(t)+s,\ \ \ \ b(t,s)=(1-s)b(t)+s.$$ Note that $|a(t,s)|\leq 1$ and $|b(t,s)|\leq 1$, with
the equality if and only if $a(t,s)=1$ and $b(t,s)=1$. Consider the map defined by
$$H(t,s)=\left ( \frac{z-a(t,s)}{1-\ov{a(t,s)}z} \right ) \left ( \frac{z-b(t,s)}{1-\ov{b(t,s)}z} \right
)$$ A brief computation shows that when $z=1$, $H(t,s)\equiv 1$. Evidently $H$ defines a homotopy between
$H(\cdot,0)=\Gamma_z$ and the constant loop $H(\cdot,1)=1$.
\end{pf}
\vspace{0.15 in}
\noindent
{\it Proof of \thmref{realize}}. Start with the closed curves $t\mapsto a(t)$ and $t\mapsto b(t)$ of
\lemref{twobranch} and form the Blaschke product
$$B^t:z\mapsto e^{-2\pi i (\nu+t)}\ z\ \left ( \frac{z-a(t)}{1-\ov{a(t)}z} \right ) \left (
\frac{z-b(t)}{1-\ov{b(t)}z} \right ).$$ 
When $t$ is not an integer, $B^t$ has degree $3$ by
\lemref{twobranch}(iv) and satisfies conditions (i) and (ii) required by \thmref{realize}. 
Moreover, it maps the unit circle $\TT$ to itself, and has no critical points in $\TT$ other than $1$,
hence ${B^t}|_\TT$ is a critical circle map.
So to finish the
proof, it suffices to show that for some $t \notin {\Bbb Z}$, the rotation number of the restriction of $B^t$
to the  circle ${\Bbb T}$ is equal to $\theta$. To this end, consider the universal covering map ${\Bbb R}
\rightarrow {\Bbb T}$ given by $z=z(w)=e^{2 \pi i w}$. Since $B^0:z\mapsto e^{-2\pi i\nu}z$, a lifting of $B^0$ to
the real line will be the affine map $\hat{B}^0:w\mapsto -\nu+w$. The loop $\{ t\mapsto B^t
\}_{0\leq t \leq 1}$ can then be lifted to a path $\{ t\mapsto \hat{B}^t \}_{0 \leq t \leq 1}$, with
$$\hat{B}^t: w \mapsto -\nu -t +w+\frac{1}{2 \pi i} \log(\Gamma_{e^{2\pi i w}}(t)),$$ 
where $\Gamma_z$ is the closed curve defined in (\ref{eqn:gam}).
Let $\rho(t)=\lim_{n\rightarrow \infty} (\hat{B}^t)^{\circ n}(w)/n$. It is a standard fact that $\rho$ is well-defined and independent of $w$ and the map $t \mapsto \rho(t)$ is continuous (see for example \cite{deMelo-vanStrien}). The rotation number of $B^t$ is then the fractional part of $\rho(t)$. 
Evidently $\rho(0)=-\nu$. 
Since $\Gamma_z$ is null-homotopic by \lemref{index}, we simply have
$\hat{B}^1:w\mapsto -\nu -1+w$, so that $\rho(1)=-\nu-1$. It follows that for some $t$ between $0$ and
$1$, $\rho(t)\equiv \theta$ (mod 1).
Hence the rotation number of the corresponding $B^t$ is $\theta$. \hfill $\Box$

\subsection{Corollaries of the construction.}
\label{subsec:Immediate corollaries}
As we shall see below, 
the Blaschke product $\BTN$ we constructed above and the  Blaschke model $\QT$ of \secref{sec:pet}
share many common properties. This will allow us to define drops, limbs, drop-chains etc. in a 
similar fashion for $\BTN$. We will also describe a quasiconformal surgery transforming
$\BTN$ into the quadratic rational map $\FTN$.

Let $0<\theta<1$ be irrational and $0<\nu<1$ be irrational of Brjuno type, and set $B=\BTN$. By
(\ref{eqn:BB}), $B(z)= e^{-2 \pi i \nu}z+O(z^2)$ near $z=0$, so by the theorem of
Brjuno-Yoccoz \cite{Yoccoz2} the origin is the center of a Siegel disk $U^0$ for $B$. We have $U^0 \subset
\DD$ since the unit circle is a subset of the Julia set. Since $B$ commutes with the reflection ${\cal T}: z\mapsto 1/\ov{z}$, there exists a Siegel disk
$U^{\infty}={\cal T}(U^0)$ centered at infinity. In the local coordinate $w=1/z$ near infinity, the map
$w\mapsto 1/B(1/w)$ has the form $w\mapsto e^{2 \pi i \nu}w+O(w^2)$, so the rotation number of
$U^{\infty}$ is $\frac{1}{2\pi i} \log B'(\infty)=\nu$.

$B$ has zeros at $\{ 0, a, b \} $ and poles at $\{ \infty, 1/\ov{a}, 1/\ov{b} \}$. The preimage
$B^{-1}(\TT)$ consists of $\TT$ and an analytic closed curve homeomorphic to a figure eight with the double
point at $z=1$. This curve and the basic dynamics of $B$ are shown in \figref{eight}.
\realfig{eight}{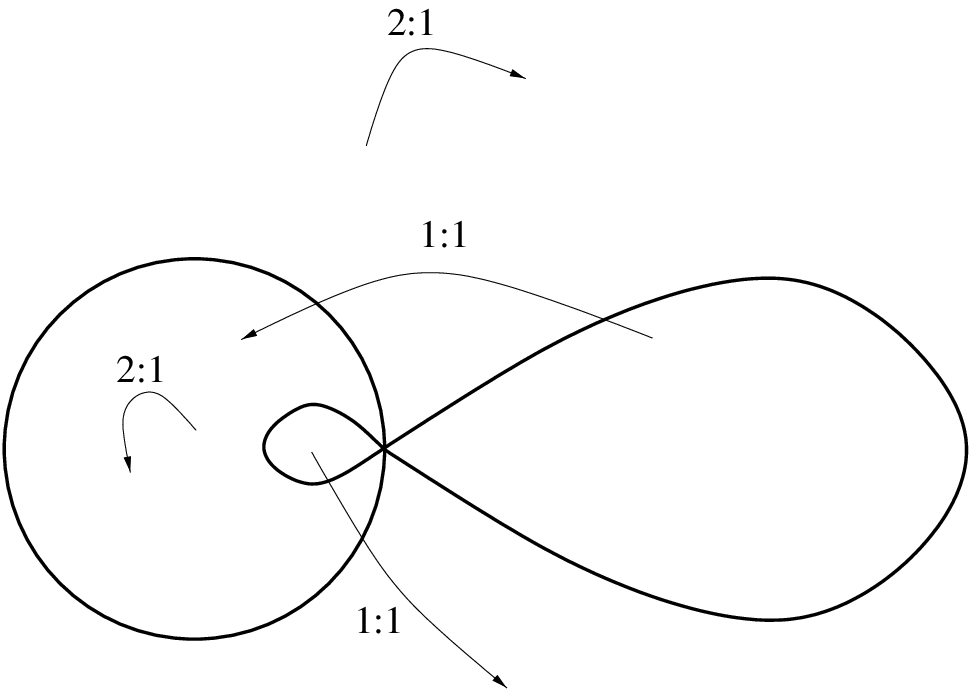}{{\sl The preimage $B^{-1}(\TT)$ and the basic dynamics of $B$.}}{7cm}
\noindent
By the theorem of Yoccoz (see subsection \ref{subsection:circle maps}), there exists a homeomorphism $h:\TT\to\TT$,
unique if we require that $h(1)=1$, such that $h\circ B|_\TT=\varrho_\theta\circ h$.
Denoting by $H:\DD\to\DD$ the Douady-Earle extension of $h$, we define the modified map $\tl B$ as
\begin{equation}
\label{eqn:modify2}
\tilde{B}(z)=\tilde{B}_{\theta, \nu}(z)= \left \{
\begin{array}{ll}
B(z) & |z|\geq 1 \\ (H^{-1}\circ {\varrho}_{\theta} \circ H)(z) & |z| \leq 1
\end{array}
\right.
\end{equation}
The map $\tilde{B}$ is a  
degree $2$ branched covering of the sphere, holomorphic outside of $\DD$. It
has a Siegel disk $U^{\infty}$ centered at $\infty$ and a ``topological Siegel disk,'' namely the unit disk $\DD$, on which its action is topologically conjugate to an irrational rotation.

The definition of drops and their addresses for the map $\tl B$
carries over word for word from subsection \ref{subsec:Drops and their addresses}. In particular, the unit disk $\DD$ is the 0-drop, and its immediate preimage $U_1=\tl B^{-1}(\DD) \sm \DD$
is the 1-drop of $\tl B$. As before, the root of the drop $U_{\io_1\io_2\ldots\io_k}$ is the point
$x_{\io_1\io_2\ldots\io_k}=\bd U_{\io_1\io_2\ldots\io_{k-1}\io_k}\cap \bd U_{\io_1\io_2\ldots\io_{k-1}}$.
%
As in subsection \ref{subsec:Drop-chains}, for each sequence of drops $\{U_0=\DD,U_{\io_1},U_{\io_1\io_2},\ldots\}$
where each $U_{\io_1 \ldots \io_k}$ is the parent of $U_{\io_1 \ldots \io_{k+1}}$, we define the drop-chain
\begin{equation}
\label{eqn:ccc}
{\cC}=\ov{\bigcup_k U_{\io_1 \ldots \io_k}},
\end{equation}
and the corresponding drop-ray $R({\cC})\subset \cC$.
We can also define the limb $L_{\io_1 \ldots \io_k}$ as the closure of the union of $U_{\io_1 \ldots \io_k}$ and all its descendants:
$$L_{\io_1\ldots\io_k}=\ov{\bigcup U_{\io_1\ldots\io_k\ldots}}\ .$$
In anticipation of the analogue of \thmref{diam->0}, let us define the {\it accumulation set} of the drop-chain $\cC$ in (\ref{eqn:ccc})
as the intersection of the decreasing sequence of limbs $L_{\io_1} \supset L_{\io_1 \io_2} \supset L_{\io_1 \io_2 \io_3}\supset \cdots$. In the case when this set is a single
point $\{ p \}$, we shall say that $R(\cC)$ or $\cC$ {\it lands at} $p$.

As an analogue to the ``filled Julia set'' $K(\tilde Q)$, we define 
$$K(\tilde{B})= K(\TBTN)= \{ z \in \CC : \mbox{The orbit $\{ \tilde{B}^{\circ n}(z) \}_{n \geq 0}$ never intersects
$U^\infty$} \} $$ 
and 
$$J(\tilde{B})=\partial K(\tilde{B}).$$ 
%
Both sets are nonempty and compact. However, $K(\tl B)$ is no longer full.
The simply-connected  basin of infinity for $\tl Q$ is replaced by the Siegel disk
$U^{\infty}$ of $\tl B$ and all its infinitely many preimages (compare \figref{goldenmate}).

\realfig{goldenmate}{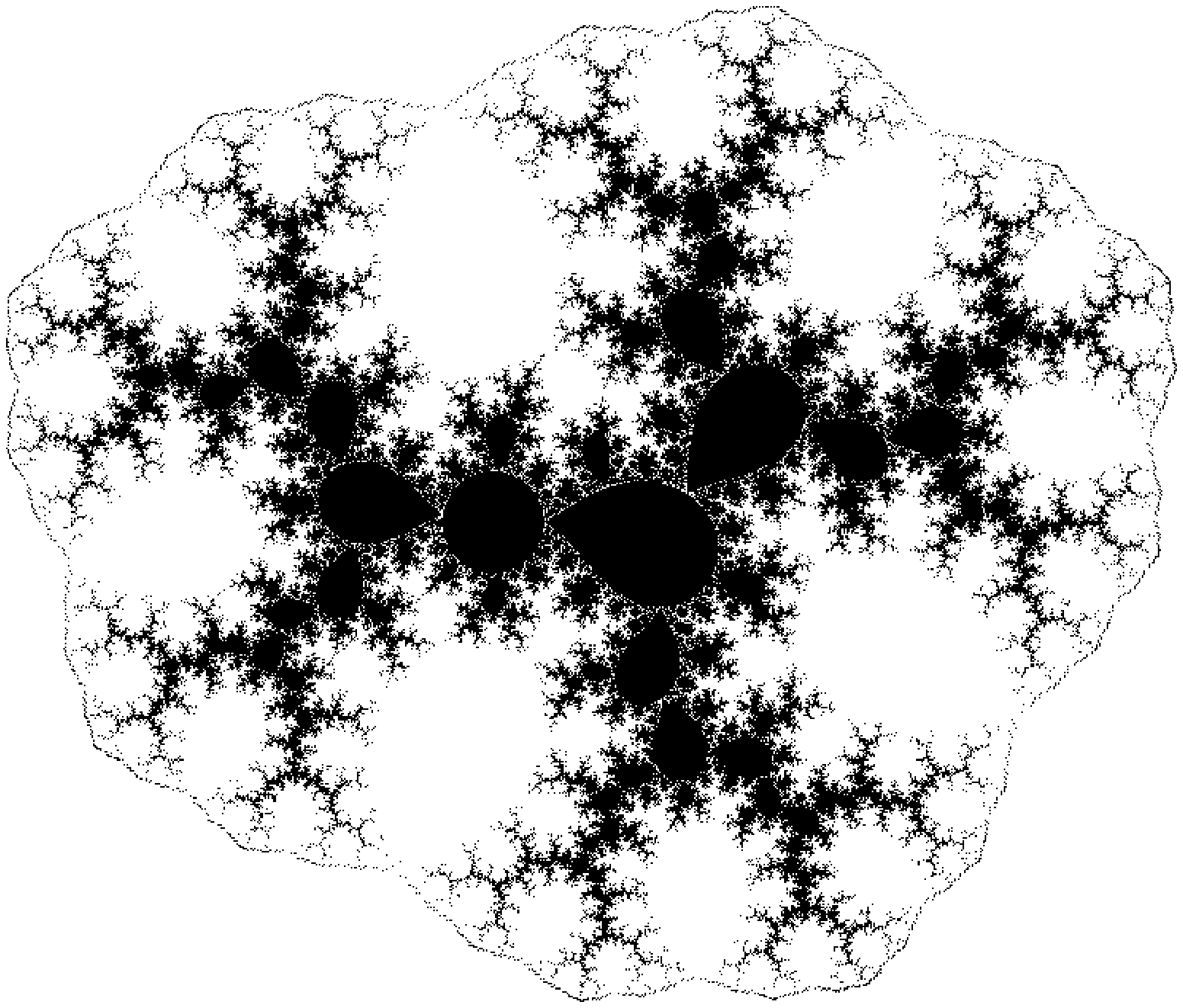}{{\sl Set $K(\tilde{B}_{\theta, \nu})$ for $\theta=\nu=(\sqrt{5}-1)/2$. Numerical experiment gives $a=-0.019048-0.298116i$, $b=3.280417-0.667122i$ for these choices of $\theta$ and $\nu$. There is a striking similarity with the corresponding picture for the quadratic rational map $F$ of \figref{blah}, up to a $90^\circ$ rotation. The reason is the existence of a quasiconformal homeomorphism conjugating $\TBTN$ to $F$ which is conformal in the white region.}}{10cm}

Finally, if $\theta$ is of bounded type, we can perform the same kind of quasiconformal surgery as in subsection \ref{subsec:Surgery}
to obtain a quadratic rational map from $\tl B$. In this case by Herman's theorem (see subsection \ref{subsection:circle maps}) the homeomorphism $h$ which linearizes $B|_{\TT}$
is quasisymmetric, therefore its Douady-Earl extension $H$ is quasiconformal.
The map $\tl B=\tl \BTN$ is a quasiregular branched covering of the Riemann sphere. We define 
a $\tl \BTN$-invariant conformal structure $\sigma_{\theta, \nu}$ on the sphere by setting it 
equal to the standard structure $\sigma_0$ on $\CC\sm K(\tl\BTN)$, to $H^*\sigma_0$ on $\DD$,
and to $(\tl\BTN^{\circ n})^*H^*\sigma_0=(\BTN^{\circ n})^*H^*\sigma_0$ on every drop of depth $n$.
The maximal dilatation of $\sigma_{\theta, \nu}$ is equal to the dilatation of $H$,
and by the Measurable Riemann Mapping Theorem, there exists a quasiconformal homeomorphism 
$\psi:\ov\CC\to\ov\CC$ with $\psi^*\sigma_0=\sigma_{\theta, \nu}$.
The conjugated map $F=\psi\circ \BTN\circ\psi^{-1}$ is a degree $2$ holomorphic branched covering of the sphere,
that is a quadratic rational map. Let us normalize $\psi$ by assuming
$\psi(\infty)=\infty$, $\psi(H^{-1}(0))=0$ and $\psi(\beta)=1$, where $\beta$ denotes the fixed point
of $\BTN$ in $\CC\sm (U^\infty\cup\DD)$. By inspection, we have $F=\FTN$ in (\ref{eqn:normal form}), so that 
$$\FTN=\psi\circ\BTN\circ\psi^{-1}.$$
Recall that $\FTN$ has two Siegel disks $\Delta^0$ and $\Delta^\infty$ centered at $0$ and $\infty$, which are the images $\Delta^0=\psi(\DD)$ and $\Delta^\infty=\psi(U^\infty)$.
As a first consequence we obtain
\begin{thm}
\label{quasicircle}
Let $0< \theta <1$ be an irrational of bounded type. Then
 the boundary of the Siegel disk $\Delta^0$ of $\FTN$ is a quasicircle passing through a single critical point of $\FTN$.
\end{thm}

Observe that there is a natural symmetry
$$\FTN={\cal I}\circ F_{\nu, \theta} \circ {\cal I},$$
where $\cal I$ is the involution $z\mapsto 1/z$.
\begin{cor}
Suppose that both $0< \theta <1$ and $0< \nu <1$ are irrationals of bounded type. Then the boundaries of the
Siegel disks $\Delta^0$ and $\Delta^\infty$ of $\FTN$ are disjoint quasicircles, each passing through a critical point of $\FTN$.
\end{cor}

\noindent
The involution $\cal I$ provides us with a quasiconformal conjugacy between $\tl \BTN$ and $\tl B_{\nu,\theta}$.
In particular, setting
$$K^\infty(\tl B_{\theta,\nu})=\ov{\CC\sm K(\tl \BTN)},$$
we have 

\begin{cor}
\label{symmetry}
There exists a quasiconformal homeomorphism of the Riemann sphere mapping the set  $K^\infty(\tl B_{\theta,\nu})$ 
to $K(\tl B_{\nu,\theta})$.
\end{cor}

\noindent
 Hence for the map $\TBTN$ we can naturally define the {\it drops growing from infinity} 
$U^\infty_{\io_1\ldots\io_k}\subset \CC\sm K(\tl\BTN)$, with $U^\infty_0=U^\infty$,
limbs growing from infinity $L^\infty_{\io_1\ldots\io_k}$, etc.

We conclude with another immediate corollary of the above construction:

\begin{corollary}
\label{common boundary}
With the above notation, $\bd K(\TBTN)= \bd K^\infty(\TBTN)$.
\end{corollary}

\begin{pf}
Under the surgery construction, both sets $\bd K(\TBTN)$ and $\bd K^\infty(\TBTN)$ correspond to the Julia set $J(\FTN)$.
\end{pf}

\vspace{0.17in}

\section{Construction of Puzzle-Pieces}
\label{sec:puzzle}
The goal of this section and the next one is to establish the following analogue of \thmref{diam->0}:
\begin{theorem}
\label{drops shrink}
Let $0< \theta, \nu <1$ be irrationals of bounded type, with $\theta \neq 1-\nu$, and consider the modified Blaschke product $\TBTN$ of (\ref{eqn:modify2}). Then as the depth of a limb $L_{\iota_1\ldots\iota_k}$ goes to infinity, $\diam (L_{\iota_1 \ldots \iota_k})$ goes to zero.
\end{theorem}
\noindent
It follows from \corref{symmetry}
that $\diam(L^\infty_{\iota_1\ldots\iota_k}) \rightarrow 0$ as $\io_1+\ldots +\io_k \rightarrow \infty$.

We start by constructing {\it puzzle-pieces}. Our construction closely parallels the one presented by Petersen in \cite{Petersen}. For simplicity, set $B= \BTN$ and $\tilde{B}=\TBTN$. Denote by $\cC$ the drop-chain 
$$\cC=\ov{U_0\cup U_1\cup U_{11}\cup U_{111}\cup \cdots}\ .$$
The following refinement of Douady-Hubbard-Sullivan Landing Theorem
can be found in \cite{Tan-Yin}:

\begin{lemma}
\label{snail}
Let $F$ be a rational map and let $\Lambda$ denote the closure of the union
of the postcritical set and possible rotation domains of $F$. Suppose
that $\gamma:(-\infty,0] \to \ov{\C} \sm \Lambda$ is a curve with
$$F^{\circ nk}(\gamma(-\infty,-k])=\gamma(-\infty, 0]$$
for all positive integers $k$. Then $\lim_{t\to -\infty}\gamma(t)$
exists and is a repelling or parabolic periodic point of $F$ whose period divides $n$.
\end{lemma}

We can apply the above lemma to the drop-chain $\cC$, setting $\gamma$ to be the drop-ray
$R(\cC)$ parameterized so that the root of the $(k+1)$-st drop corresponds to $t=-k$. We conclude that
 $R(\cC)$ lands at the unique fixed point $\beta$ of $B$ in $\ov{\CC}\sm({\D}\cup U^\infty)$.
Since $\beta$ is necessarily repelling, the size of the drops in $\cC$ decreases geometrically,
and the drop-chain $\cC$ lands at the point $\beta$. Repeating the argument,
we see that the drop-ray $R(\cD)$ associated to the drop-chain 
$$\cD=\ov{U^\infty \cup U^\infty_1\cup U^\infty_{11}\cup U^\infty_{111}\cup \cdots}$$
lands at a fixed point as well, which is necessarily $\beta$.
Let $\cC'$ be the drop-chain $\ov{U_0 \cup U_2 \cup U_{21} \cup \cdots }$  mapped to $\cC$ by $\tilde B$,
and similarly define the drop-chain $\cD'=\ov{U^\infty \cup U^\infty_2\cup U^\infty_{21} \cup \cdots}\ .$ 
Then $\cC'$ and $\cD'$ have a common landing point $\beta'\neq \beta $,
which is a preimage of $\beta$ in $\ov{\C}\sm({\D}\cup U^\infty)$.

\realfig{spine1}{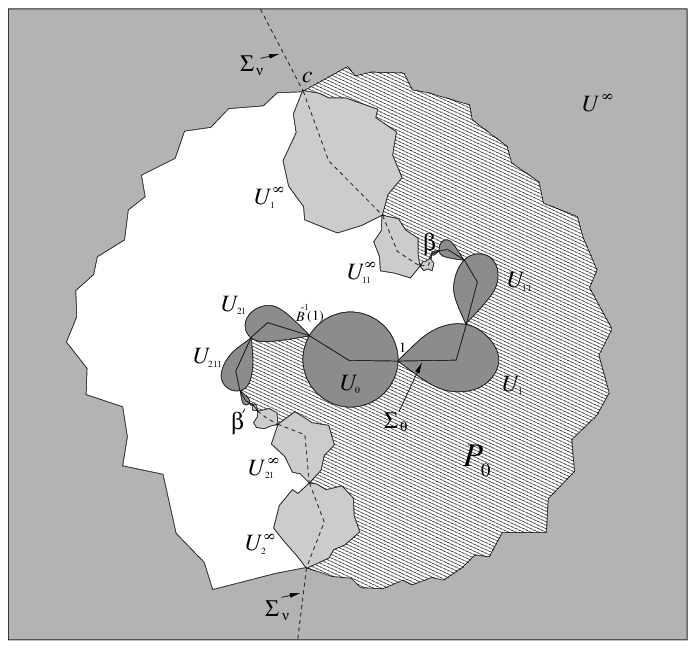}{{\sl The $0$-th critical puzzle-piece $P_0$ and the ``spines'' $\Sigma_{\theta}$ and $\Sigma_{\nu}$ (see \secref{sec:proof}).}}{11cm}

As before,
the moments of closest returns of the critical point $z=1$ are denoted by $\{ q_n \}$. Recall that these numbers appear
as the denominators of the convergents of the continued fraction of $\theta$.
We define the {\it $0$-th critical puzzle-piece} $P_0$ as the closure of
the connected component of 
$$\ov{\C}\sm ({\D}\cup U^\infty \cup \cC \cup \cC' \cup \cD \cup \cD')$$
which contains the arc $[1,B^{-1}(1)]\ni B^{\circ q_1}(1)$ in the boundary (see \figref{spine1}).  We inductively define
the $n$-th critical puzzle-piece $P_n\subset \ov{\C}\sm\D$ as the closed set which is
mapped homeomorphically onto $P_{n-1}$ by $B^{\circ q_n}$ and which contains the arc 
$[1,B^{-q_n}(1)]\subset \T$ in the boundary.
The following proposition summarizes some of the properties of critical puzzle-pieces:

\begin{proposition}[Properties of puzzle-pieces]
\label{puzzle properties}

\noindent
\begin{enumerate}
\item[(i)]
The puzzle-piece $P_n$ intersects the unit circle $\T$ along the arc 
$[1,B^{-q_n}(1)]$.
\item[(ii)]
$B^{\circ q_n}(P_n\cap \bd U_1)=[B^{\circ q_n}(1),B^{-q_{n-1}}(1)]$.
\item[(iii)]
$B^{\circ q_n+q_{n-1}+q_{n-2}}(P_n\cap \bd U_{q_n+1})=[1,B^{\circ q_{n-1}+q_{n-2}}(1)]$.
\item[(iv)]
$P_n$ contains the drop $U_{q_{n+2}+1}$.
\end{enumerate}
\end{proposition}

\begin{pf}
Observe that $B^{\circ q_n}$ is a homeomorphism $[B^{-q_n}(1),B^{-q_{n}-q_{n-1}}(1)]\iso [B^{-q_{n-1}},1]$ with one critical point at $1$. Thus the univalent inverse branch $B^{-q_n}$ sending $P_{n-1}$ to $P_n$ maps the arc $[B^{ -q_{n-1}},1]$ onto the union of $[1,B^{-q_n}(1)]$ and a subarc of $\bd U_1$.
The first three statements now follow by induction on $n$.
As seen from the combinatorics of closest returns (see \secref{subsection:circle maps})
$\bd U_{q_{n+2}+1}\cap \TT=B^{-q_{n+2}}(1)$  is contained in the arc $[1,B^{-q_n}(1)]$. Evidently, 
the drop $U_{q_{n+2}+1}$  has no intersections with $\bd P_n$, thus $U_{q_{n+2}+1}\subset P_n$. 
\end{pf}
\comm{The following simple lemma will be useful in subsequent arguments:

\begin{lemma}
\label{cut}
Let $U$ be a topological disk whose boundary is contained in a finite union of the boundary arcs of drops of $\tilde{B}$ (the preimages of the unit disk $\DD$). Then $U$ itself must be a drop.
\end{lemma}

\begin{pf}
Note that the modified Blaschke product $\tilde{B}$ satisfies the Maximum Principle in $\CC \ov U_1^{\infty}$. Since $\tilde{B}^{\circ n}(\bd U) \subset \TT$ for a large $n$, we must have $\tilde{B}^{\circ n}(U) \subset \DD$, which means $U$ itself is a drop.
\end{pf}
}
The preimages of the puzzle-piece $P_0$ have the following nesting property:
\begin{lem}
\label{nesting property}
Let $A_1$ and $A_2$ be two distinct univalent
 pull-backs of the puzzle-piece $P_0$ such that $\overset{\circ}{A_1}\cap
\overset{\circ}{A_2}\ne\emptyset$. Then either $A_1\subset A_2$ or $A_2\subset A_1$.
\end{lem}
\begin{pf}
By construction, the boundary of the puzzle-piece $P_0$ consists of an open arc $\gamma\subset \cC\cup\cC'$ 
which is made up of the boundary arcs of various drops $U_{\io_1\ldots\io_k}$, a similarly defined arc $\gamma^\infty \subset \cD \cup \cD'$ and points
$\beta$, $\beta'$ (see \figref{spine1}). Denote by $\gamma_1$, $\gamma_1^\infty$, $\beta_1$, $\beta_1'$ the corresponding parts
of $\bd A_1$, and label the boundary of $A_2$ in the same way.

Evidently $\gamma_1$ does not intersect $\gamma^\infty_2$ or the points $\beta_2$, $\beta_2'$, so it can only intersect $\gamma_2$. Similarly, $\gamma_1^{\infty}$ can only intersect $\gamma_2^{\infty}$. Suppose that $y \in \{ \beta_1, \beta_1' \} \cap \{ \beta_2, \beta_2' \}$. Then $B^{-k}(\beta)=y$ for some choice of the inverse branch. Since $\beta$ is not in the post-critical set of $B$, this branch of $B^{-k}$ has a univalent extension to a neighborhood of $\beta$ intersecting the boundary of $P_0$ along a non-empty open arc. Pulling back, it follows that for some neighborhood $D$ of $y$, $\gamma_1\cap D=\gamma_2 \cap D$ and 
$\gamma_1^\infty \cap D=\gamma_2^\infty \cap D$.

Now assume that the claim is false. Let $A_1=B^{-m}(P_0)$ and $A_2=B^{-n}(P_0)$, with $m\leq n$. Then by the above observation, either $\gamma_2$ or $\gamma_2^{\infty}$ intersects both $\stackrel{\circ}{A_1}$ and $\CC \sm A_1$. Therefore, either $B^{\circ m}(\gamma_2)$ or $B^{\circ m}(\gamma_2^{\infty})$ intersects both $\stackrel{\circ}{P_0}$ and $\CC \sm P_0$. To fix the ideas, let us assume that $B^{\circ m}(\gamma_2)$ does. Note that $B^{\circ m}(\gamma_2)\cap \bd P_0 \subset \gamma$, hence $B^{\circ m}(\gamma_2)$ must intersect the union of the drop-rays $R(\cC) \cup R(\cC')$ {\it transversally} at a root $x$ of a drop in $\cC \cup \cC'$. Now under $B^{\circ n-m}$ a small open subarc of $B^{\circ m}(\gamma_2)$ around $x$ maps homeomorphically to a subarc $\delta \subset \gamma$ around $B^{\circ n-m}(x)$. Since the orbit $x, B(x), \ldots, B^{\circ n-m}(x)$ does not contain the critical point $1$, it follows that $\delta$ also intersects $R(\cC) \cup R(\cC')$ transversally at $B^{\circ n-m}(x)$, which is impossible.
\end{pf}
\comm{
Let $A$ be an interior component of the symmetric difference
$(A_1\sm A_2)\cup (A_2\sm A_1)$. Suppose that $A$ is bounded by a finite collection of boundary arcs 
of drops $U_{\io_1\ldots\io_k}$. 
Applying the Maximum Principle, we see that $A$ itself is a preimage of the unit disk $U_0$.
Thus two distinct preimages of $U_0$ share a non-trivial boundary arc, which is impossible. Similarly, $A_3$ can not be bounded
by a finite collection of boundary arcs of drops at infinity. 
This implies that  $(A_1\sm A_2)\cup (A_2\sm A_1)$ has exactly two interior
components.
It remains to show that  $\{\beta_{A_1}, \beta_{A_1}'\}\cap \{\beta_{A_2}, \beta_{A_2}'\}=\emptyset$.
Let us consider the possibility that there are points $y_1\in\{\beta_{A_1},\beta_{A_1}'\}$,
$y_2\in \{\beta_{A_2},\beta_{A_2}'\}$ such that $y_1=y_2$. Thus $B^{-k}(\beta)=y_1=y_2$ for some choice of 
the inverse branch. Since $\beta$ is not in the post-critical set of $B$, this branch of $B^{-k}$ 
has a univalent extension to a neighborhood of $\beta$ intersecting the boundary of $P_0$ over a non-empty open arc.
Thus $A_1$ and $A_2$ share a common boundary arc which contains $y_1=y_2$ in the interior. It follows that some 
component of $(A_1\sm A_2)\cup (A_2\sm A_1)$ is bounded by a finite collection of boundary arcs of drops,
and we use the same argument as above to arrive at a contradiction.

Let us label $A_3$ the set $A_2\sm A_1$. By the above discussion the boundary of $A_3$ contains two of the four points 
$\beta_{A_1}$, $\beta_{A_1}'$, $\beta_{A_2}$, $\beta_{A_2}'$.
To fix our ideas, let us assume that these points are  $\beta_{A_1}$ and  $\beta_{A_2}'$,
the other possibilities are handled in the same way. 
Consider the image $B^{\circ n}(A_3)$. Since $B^{\circ n}(A_2)=P_0$, $B^{\circ n}(A_3)\cap \overset{\circ}{P_0}\ne \emptyset$.
Moreover, $B^{\circ n}(\beta'_{A_2})=\beta'$, hence $\bd B^{\circ n}(A_3)$ contains an arc $\delta'\subset \bd P_0$
with the point $\beta'$ in the interior. Also, $B^{\circ m}(\beta_{A_1})=\beta$, and thus $\bd B^{\circ m}(A_3)$ contains
an arc $\delta\subset \bd P_0$ having $\beta$ in the interior. By construction, $B(\delta)\supset \delta$, hence
$\delta\subset \bd B^{\circ n}(A_3)$. The properties of drops imply that $\bd B^{\circ n}(A_3)\subset \cC\cup \cC'\cup \cD\cup \cD'$.
Thus $\bd B^{\circ n}(A_3)$ partitions the sphere in two regions, one of which contains $P_0$. As 
$B^{\circ n}(A_3)\cap \overset{\circ}{P_0}\ne \emptyset$, $B^{\circ n}(A_3)\supset P_0$. This implies that $A_3=A_2$,
contradictory to the assumptions of the Lemma.
}

\begin{cor}
\label{puzzle-pieces nested}
For all $n\geq 0$ we have $P_{n+2}\subsetneqq P_n$.
\end{cor}
\begin{pf}
It is clear from the definition of critical puzzle-pieces that $\overset{\circ}{P}_{n+2}\cap \overset{\circ}{P}_n\ne \emptyset$.
By \propref{puzzle properties}(i), $P_{n+2}\cap \TT \subsetneqq P_n\cap\TT$. The claim now follows from \lemref{nesting property}.
\end{pf}

\begin{lemma}
\label{cut}
Let $U$ be a topological disk whose boundary is contained in a finite union of the boundary arcs of drops (resp. drops growing from infinity). Then $U$ itself must be a drop (resp. drop growing from infinity).
\end{lemma}

\begin{pf}
Let us consider the case of drops. The proof for the case of drops growing from infinity is similar. The modified Blaschke product $\tilde{B}$ is an open mapping, so it satisfies the Maximum Principle in $\CC \sm U_1^{\infty}$. Since $\tilde{B}^{\circ n}(\bd U) \subset \TT$ for a large $n$, we must have $\tilde{B}^{\circ n}(U) \subset \DD$, which means $U$ itself is a drop.
\end{pf}

\begin{lem}
\label{limb inside}
Let $A$ be a univalent pull-back of the puzzle-piece $P_0$. 
Suppose that a drop at infinity $U^\infty_{\io_1\ldots\io_k}$ is contained in $A$. Then $A$ contains the whole
limb $L^\infty_{\io_1\ldots\io_k}$.
\end{lem}

\begin{pf} 
Let us denote by $\gamma_A^\infty\subset \bd A$
the part of the boundary of $A$ made up of the boundary arcs of drops at infinity. Assume by way of contradiction that there is a drop at infinity $U^\infty_{\io_1\ldots\io_k\ldots\io_{k+m}}\not\subset A$.
Let $\cD$ be a drop-chain containing $U^\infty_{\io_1\ldots\io_k\ldots\io_{k+m}}$. Let $\delta \subset \bd \cD$ be an arc connecting the root of 
$U^\infty_{\io_1}$ to a point in $\bd U^\infty_{\io_1\ldots\io_k\ldots\io_{k+m}} \sm \ov{A}$. Then $\delta$ goes in and out 
of $A$, but it only intersects $\bd A$ at the points of $\gamma_A^\infty$.
Thus the curves $\delta$ and $\gamma_A^\infty$ 
 bound a topological disk $U\subset \overset{\circ}{A\ }$. By \lemref{cut}, $U$ itself is a drop growing from infinity. Since $U$ shares a non-trivial boundary arc with another drop growing from infinity, we arrive at a contradiction.
\end{pf}

\begin{lemma}
\label{commensurable disk}
The puzzle-piece $P_n$ contains a Euclidean disk $D$ centered at a point in $J({B})$ with $\diam (D) > K|[1,B^{-q_n}(1)]|$ 
for some $K$ independent of $n$.
\end{lemma}

\begin{pf}
Note first that by \propref{puzzle properties}(iv), $U_{q_{n+2}+1}\subset P_n$.
Since $B^{\circ q_{n+2}}(1)$ is a closest return of the critical point $1\in\TT$, $B^{-q_{n+2}}|_\TT$ maps
the arc  $(B^{-q_{n+2}}(1),B^{\circ q_{n+2}}(1))$ diffeomorphically onto $(B^{-2q_{n+2}}(1),1)$.
This inverse branch has a univalent extension to a neighborhood of $1$, which we denote by $\psi_n$.
By \'Swia\c\negthinspace tek-Herman real a priori bounds (see the discussion in the end of 
 \secref{subsection:circle maps}), the segments $[B^{-2q_{n+2}}(1),B^{-q_{n+2}}(1)]$,
$[B^{-q_{n+2}}(1),1]$ and $[1,B^{\circ q_{n+2}}(1)]$ are $K_1$-commensurable. Here the constant $K_1$
 becomes universal for sufficiently large $n$
and therefore can be chosen independent of $n$.  Moreover, $1/K_2\leq|\psi_n'(1)|\leq K_2$ for some 
$K_2>1$ which is  also independent of $n$.
 By Koebe Distortion Theorem we may choose a Euclidean disk $D$ around the point $1$
commensurable with $[1,B^{\circ q_{n+2}}(1)]$ such that $\psi_n$ has bounded distortion in $D$.
Now let us pull back a sub-disk $D'\subset D$ centered at a point in $\bd U_1$
to obtain a Euclidean disk $D_1\subset \CC\sm\DD$ around a point in $\bd U_{q_{n+2}+1}$
such that both $\diam D_1$ and $\dist(D_1,B^{-q_{n+2}}(1))$ are  $K_3$-commensurable with 
$[1,B^{\circ q_{n+2}}(1)]$ for some  $K_3$ independent of $n$.

Denote by $D_1'\subset \DD$ the disk symmetric to $D_1$ with respect to $\TT$. Let $D_2\subset U_1$ be given by 
$B(D_2)=D_1'$. It is clear that $D_2$ is again commensurable with $[1,B^{\circ q_{n+2}}(1)]$,
and so is $\dist(D_2,1)$. By Koebe Distortion Theorem, the image $\psi_n(D_2)\subset U_{q_{n+2}+1}\subset P_n$
contains a Euclidean disk with the desired properties.
\end{pf}

The last property of puzzle-pieces we  need is the following:

\begin{lemma}
\label{leaves delta}
There exists $N>0$ such that for all $n\geq N$ the puzzle-piece $P_n$ does not intersect $\bd U^\infty$.
\end{lemma}

\begin{pf}
Since the boundary of the Siegel disk $U^\infty$ is forward-invariant, we only
need to show the existence of one $N$ such that $P_N\cap\bd U^\infty=\emptyset$.
 Assume this is false. Let us denote by $l_n$ the  boundary arc of $P_n$
connecting $1$ to $\bd U^\infty$. By \lemref{nesting property}, the  curves in the 
 orbit 
\begin{equation}
\label{curves orbit}
l_n,B(l_n),\ldots, B^{\circ q_n-1}(l_n)
\end{equation} are disjoint. 
By the theorem of Yoccoz (see subsection \ref{subsection:circle maps}) the maps $B|_\TT$ and $B|_{\bd U^\infty}$ are topologically
conjugate to rigid rotations. Since the inverse orbit of a point under an irrational rotation is dense on the circle,
the maximum diameter of the pieces into which the curves (\ref{curves orbit}) partition the boundaries of
 $\D$ and $U^\infty$ goes to zero as $n\to\infty$. We may therefore construct an orientation-reversing 
topological conjugacy between the circle maps $B|_{\TT}$ and $B|_{\bd U^\infty}$. This contradicts the fact that $\theta\ne 1-\nu$.
\end{pf}

\vspace{0.17in}
\section{Complex Bounds}
\label{sec:bounds}

The proof of Petersen's Theorem presented in \cite{Yampolsky} is based
on a version of estimates employed in the same paper for proving
a renormalization convergence result. In renormalization theory 
it is customary to use the term {\it complex a priori bounds}
for such estimates. Our goal in this section is to adapt these bounds to the
 Blaschke product model introduced in \secref{sec:model}.

As before, let us fix irrationals $0< \theta,\nu <1$ of bounded type, with $\theta \neq 1-\nu$, and set $B=\BTN$, $\tilde{B}=\TBTN$. Recall that $B$ is a Blaschke product of the form 
$$B= z\mapsto \lambda \ z\ \left ( \frac{z-a}{1-\ov{a}z} \right ) \left ( \frac{z-b}{1-\ov{b}z} \right ),$$
where $|\lambda|=1$, $0<|a|<1$ and $|b|=|a|^{-1}$. We set 
$$B(1)=e^{2 \pi i \tau}\ \ \operatorname{with}\ 0<\tau<1.$$
The convergents of the continued fraction $\theta=[a_1, a_2, a_3,\ldots]$
will be denoted $\{ p_n/q_n \}$. 
First note that $(B(z)-B(1))/(z-1)^3$ is a bounded holomorphic function in
the domain ${\C}\sm \ov{({\D}\cup U^\infty\cup U^\infty_1)}$.
As a consequence,
\begin{equation}
\label{cubic estimate}
C^{-1}|z-1|^3<|B(z)-B(1)|<C|z-1|^3
\end{equation}
in this domain, for some positive constant $C$.

Let $S$ be the translation-invariant  infinite strip
which is mapped onto the open topological annulus $\CC\sm (\ov{U^0} \cup \ov{U^\infty})$ by the exponential map $z\mapsto e^{2\pi i z}$. Let us denote by $S_J$ the domain obtained
by removing from $S$ the points of the real line that do not belong to the interval $J\subset \RR$:
$$S_J=(S\sm \RR)\cup J.$$
Let $\hat{B}(z)$ denote the (multi-valued) meromorphic function $\frac{1}{2\pi i} \log B(e^{2\pi i z})$ on $S$. On the real line $\hat{B}$ has singularities 
at the integer points, whose images lie at the integer translates
of $0<\tau<1$. Its other singularities
lie at the boundary curves of $S$ at the points $\pm s +j,\;j\in\ZZ$,
which are mapped by the exponential map to the critical points on
the boundaries of the Siegel disks $U^0$ and $U^\infty$ of $B$.
 By the Monodromy Theorem, in the domain $S_{(\tau+i,
\tau+i+1)}$ with the critical values removed, we have well-defined branches
$\phi_{i,m}$ of the inverse $\hat{B}^{-1}$, mapping the 
open interval $(\tau+i,\tau+i+1)$
homeomorphically onto the interval between two consecutive integers $(m,m+1)$ (see \figref{strip}).
The maps $\phi_{i,m}$ range over the simply-connected regions
\begin{equation}
\label{range of phi}
S_{(m,m+1)}\sm \left [ (\frac{\pm 1}{2\pi i}\log (\ov{U}^{\infty}_1)) \cup
(\frac{\pm 1}{2\pi i}\log (\ov{U}_1)) \right ].
\end{equation}
\realfig{strip}{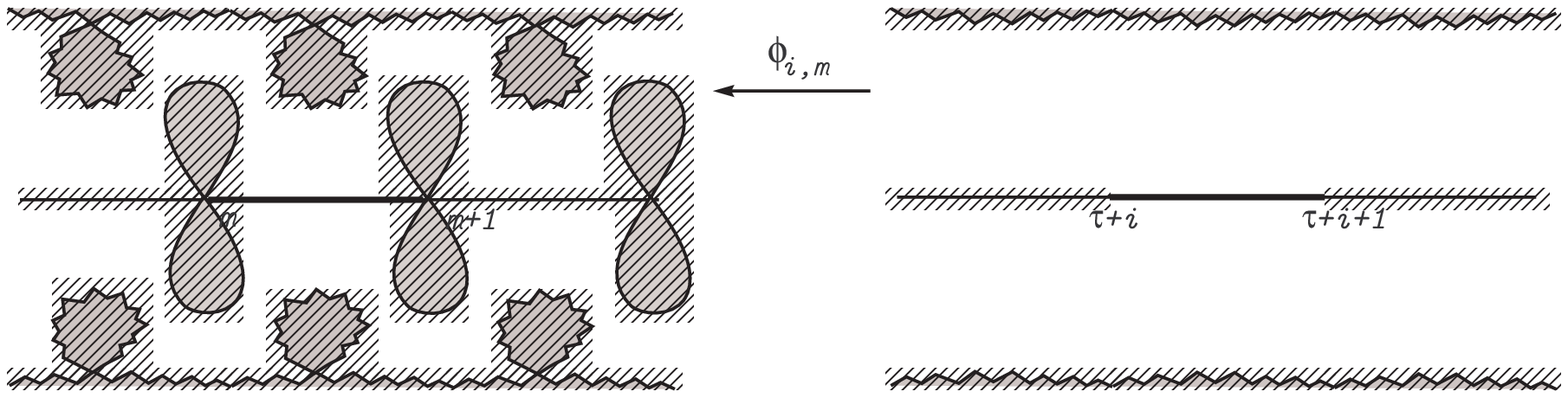}{}{\hsize}

Denote by $\Upsilon:{\TT}\sm\{B(1)\}\to I=(\tau-1,\tau)$ 
the single-valued branch of $\frac{1}{2\pi i}\log(z)$ mapping $1$ to $0$.
Define the (discontinuous) map $\phi:I\to I$ by
$$\phi(z)=\left\{ 
\begin{array}{ll}
\phi_{-1,0}(z) & \text{ for }z\in(\tau-1,\Upsilon(B^{\circ 2}(1))],\\
\phi_{-1,-1}(z) & \text{ for }z\in(\Upsilon(B^{\circ 2}(1)),\tau).
\end{array}
\right.$$

\noindent
Let us fix an  $n\geq 2$ and consider the inverse orbit
\begin{equation}
\label{circle-orbit}
(1,B^{\circ q_n}(1)),(B^{-1}(1),B^{\circ q_n-1}(1)),\ldots,(B^{ -q_n}(1),1).
\end{equation}
Set $J_{-i}=\Upsilon((B^{ -i}(1),B^{\circ q_n-i}(1)))$ and consider the $\phi$-orbit
\begin{equation}
\label{J orbit}
J_0,J_{-1},J_{-2},\ldots,J_{-q_n}.
\end{equation}
By the combinatorics of closest returns (see subsection \ref{subsection:circle maps}) the smallest value of $i>0$ for which
the arc $B^{-i}((B^{-q_n}(1),1))\subset \TT$  contains the critical point $1$ is $q_{n+1}$.
Also, the smallest $j>0$ for which $1 \in B^{\circ j}((B^{-q_n}(1),1))$ is $q_{n+1}+q_n$.
As $q_{n+1}\geq q_n+2$, the interval $(B^{-k-2}(1),B^{\circ q_n-k-2}(1))$ does not contain $1$ for $0\leq k\leq q_n-1$.
Hence, $B^{\circ 2}(1) \notin (B^{-k}(1),B^{\circ q_n-k}(1))$ for $0\leq k\leq q_n-1$. In other words,
the intervals $J_{-k}$ of the orbit (\ref{J orbit}) for $0\leq k\leq q_n-1$ do not contain the point of discontinuity
of the map $\phi$. By its definition, the map $\phi:J_{-k}\to J_{-k-1}$ for $0\leq k\leq q_n-1$ has a univalent
extension to $S_{J_{-k}}$. As seen from  (\ref{range of phi}) the range of this univalent map is a subset of $S_{J_{-k-1}}$,
hence the composition $\phi^l:J_{-i}\to J_{-i-l}$ for $0\leq i< i+l\leq q_n$ univalently extends to the entire $S_{J_{-i}}$.

Consider the univalent extensions of the iterates $\phi^k:J_0\to J_{-k}$ to the strip
$S_{J_0}$ for $1\leq k\leq q_n$. Applying these univalent branches to a point
 $z\in S_{J_0}$, we obtain the {\it inverse orbit, corresponding to the orbit (\ref{J orbit})}
\begin{equation}
\label{z orbit}
z=z_0,z_{-1},z_{-2},\ldots,z_{-q_n},\text{ where }z_{-k}=\phi^k(z_0).
\end{equation}
A corresponding inverse orbit of a subset of $S_{J_0}$ is similarly defined.

Let ${\CC}_J\supset S_J$ denote the slit plane $({\CC}\sm{\RR})\cup J$. One easily constructs
a conformal mapping of this domain to the upper half-plane to verify that
the hyperbolic neighborhood $\{z\in{\CC}_J|\dist_{{\CC}_J}(z,J)<r\}$ for $r>0$
is the union $D_\theta(J)$ of two Euclidean disks of equal radii
with common chord $J$ intersecting the real axis at an outer angle
$\theta=\theta(r)$ (see \cite{deMelo-vanStrien}). An elementary computation yields in this case
$$r=\log\tan(\pi/2-\theta/4).$$
\noindent
The standard properties of conformal maps imply that the hyperbolic 
neighborhood $\{z\in {S}_J|\dist_{{S}_J}(z,J)<r\}$
also forms angles $\theta=\theta(r)$ with $\RR$. We choose the notation $G_\theta(J)$ for this neighborhood. The Schwarz Lemma implies that $G_\theta (J)\subset D_\theta(J)$.

Let $\breve S\subset \CC$ be a horizontal strip invariant under the unit translation, which is compactly contained in $S$.
A specific choice of $\breve S$ will be made later in our arguments (see the remarks before \lemref{ind3}).
Let $I$ be a bounded interval in $\RR$. For a point $z\in S_I$ not belonging to $\RR$, denote by $0< \ang{z}{I} < \pi/2 $ the least of the outer angles the segments joining $z$ to the end-points of $I$ form with the real line.
The following adaptation of Lemma 2.1 of \cite{Yampolsky} will be used
to control the expansion of inverse branches:

\begin{lem}
\label{good angle}
Let us fix $n$ and consider the inverse orbit (\ref{z orbit}). Let $k\leq q_n-1$.
Assume that for some $i$ between $0$ and $k$, $z_{-i}\in\breve S$ and  $\ang{z_{-i}}{J_{-i}}>\eps>0$.
Then
$$\frac{\dist(z_{-k},J_{-k})}{|J_{-k}|}\leq C\ \frac{\dist(z_{-i},J_{-i})}
{|J_{-i}|}$$
for some constant $C=C(\eps,\breve S)>0$.
\end{lem}
\begin{pf}
First observe that $B^{-q_n}|_\TT$ is a diffeomorphism on the arc $[B^{\circ 2q_n}(1),B^{-q_n}(1)]\subset \TT$
which contains the arc $[B^{\circ q_n}(1),1]$ in its interior. Moreover, by 
\'Swia\c\negthinspace tek-Herman real a priori bounds (see subsection \ref{subsection:circle maps}),
the latter arc is contained well inside of the former. As seen from the combinatorics of closest returns,
the iterates $B^{-j}([B^{\circ 2q_n}(1),B^{-q_n}(1)])$ do not contain $B^{\circ 2}(1)$ for $j\leq q_n-1$.
Setting $H=\Upsilon([B^{\circ 2q_n}(1),B^{-q_n}(1)])$, we see that $J_0$ is contained well inside of $H$,
and $\phi^{j}:J_0\to J_{-j}$ univalently extends to $S_H$ for $1\leq j\leq q_n-1$. 
Set $T=\phi^i(H)\supset J_{-i}$.
By Koebe Distortion Theorem, there exists $\rho>0$ such that both 
components of $T\sm J_{-i}$ have length at least $2\rho|J_{-i}|$.
Note that the iterate $$\phi^{\circ k-i}:J_{-i}\to J_{-k}$$
has a univalent extension to $S_T$.

Let us normalize the situation by considering the orientation-preserving
affine maps
$$\alpha_1:J_{-i}\to [0,1]\text{ and }\alpha_2:J_{-k}\to [0,1].$$
The composition $\alpha_2\circ\phi^{\circ k-i}\circ \alpha_1^{-1}$ is defined in a straight horizontal strip 
$$Y=\{z\in {\CC}_{[-2\rho,1+2\rho]}:|\operatorname{Im} z|<M\}$$
 for some $M>0$ independent of $n$. The space of normalized univalent maps 
of $Y$ is compact by Koebe Theorem, thus the statement
is true if $\dist(z,J_{-i})/|J_{-i}|<\rho$.

Now assume $\dist(z,J_{-i})/|J_{-i}|>\rho$. Consider the smallest closed
hyperbolic neighborhood $\ov{G_\theta(J_{-i})}$
containing $z_{-i}$. Recall that $z_{-i}$ is contained in a  strip 
$\breve S\Subset S$. For a point $\zeta\in {\CC}_I$ with $\dist(\zeta,I)>\rho|I|$ and $\ang{\zeta}{I}>\eps$, the smallest closed neighborhood $\ov{D_\theta(I)}\ni \zeta$ satisfies $\diam D_\theta(I)\leq C(\rho,\eps)\dist(\zeta,I)$ (see \cite{Yampolsky}, Lemma 2.1). Therefore, we have $\diam G_\theta(J_{-i})\leq C(\rho,\eps,\breve S)\dist (z,J_{-i})$ and by Koebe Theorem,
$$\frac{\diam G_\theta(J_{-i})}{|J_{-i}|}\sim
\frac{\diam G_\theta(J_{-k})}{|J_{-k}|}.$$
By the Schwarz Lemma, $z_{-k}\in G_\theta(J_{-k})$ and the claim follows.
\end{pf}

Set $I_m=\Upsilon([1,B^{\circ q_m}(1)])$, and let $G_m$ denote the 
hyperbolic neighborhood 
$$G_\alpha(\Upsilon([B^{\circ q_{m+1}}(1),B^{q_m-q_{m+1}}(1)]))$$ where 
$0<\alpha<\pi/2$ will be specified later.
The following two lemmas are direct adaptations of Lemmas 4.2 and 4.4 of \cite{Yampolsky}, for which the reader is referred for a detailed discussion supplemented with figures. In both lemmas we work with the orbit (\ref{z orbit}) for some fixed value of $n$.

\begin{lem}     
\label{ind1}
Let $J$ and $J'$ be two consecutive returns of the orbit (\ref{J orbit}) of $J_0$ to $I_m$ for  $n>m>1$
 and let $\zeta$, $\zeta'$ be the corresponding points
of the inverse orbit (\ref{z orbit}). If $\zeta\in G_m$, then either
$\zeta'\in G_m$ or $\ang{\zeta'}{J'}>\eps$ and $\dist(\zeta',J')<C|I_m|$
where the constants $\eps$ and $C$ are independent of $m$.
\end{lem}

\noindent
We remark that the constants $\eps$ and $C$ will in general depend on the
choice of the Blaschke product $B$. The argument is illustrated 
in \figref{ind1fig}.
\begin{pf}
Note that $J=J_{-i}$ and $J'=J_{-i-q_{m+1}}$ for some $i<q_n-q_{m+1}$. Recall that 
$G_m=G_\alpha(\Upsilon([B^{\circ q_{m+1}}(1),B^{q_m-q_{m+1}}(1)]))$.
Let $\check G_m$ denote the pull-back of $G_m$ along the inverse orbit
$J,\ldots,J'$. Also let $ G'_m$ denote the pull-back of $G_m$ along the
piece of the orbit $J,\ldots,\phi^{\circ q_m-1}(J)$, and let $G''_m=\phi(G'_m)$.

The combinatorics of closest returns (see subsection \ref{subsection:circle maps}) implies that the restriction $B^{-q_{m-1}}|_{(B^{\circ q_{m+1}}(1),B^{q_m-q_{m+1}}(1))}$ is a diffeomorphism.
Hence the pull-back  of $G_m$ along the orbit $J,\ldots,\phi^{\circ q_m-1}(J)$
is univalent. By the Schwarz Lemma, 
$$G'_m\subset G_\alpha(\Upsilon([B^{\circ q_{m+1}-q_m+1}(1),B^{1-q_{m+1}}(1)])).$$
By \'Swia\c\negthinspace tek-Herman real a priori bounds, the 
critical value $\tau$ divides the interval $\Upsilon([B^{\circ q_{m+1}-q_m+1}(1),B^{1-q_{m+1}}(1)])$
into $K_1$-commensurable pieces,
 where $K_1$ becomes universal for large $m$, and can
therefore be chosen simultaneously for all $m$.
As the absolute value of the derivative of the exponential map is bounded
away from $0$ and $\infty$ on the strip $S$, the estimate (\ref{cubic estimate}) is 
still valid for the lifted map near the critical point. Together
with  elementary properties of the cube root map
this implies that $G''_m\subset G_\beta([\Upsilon(B^{\circ q_{m+1}-q_m}(1),1]))$
for some $\beta>0$ independent of $m$. Let $V_0 \subset S$ be the union of the connected components of $\frac{\pm 1}{2 \pi i} \log (\ov{U}_1)$ attached to $0$ (see \figref{ind1fig}). Since the boundary
of $ G''_m$ contains a segment of $\partial V_0$ which forms outer angles
$\pi/3$ with $\RR$ at $0$, we have $G''_m\subset G_\gamma([\Upsilon(
B^{\circ q_{m+1}-q_m}(1)),a_1])\cup G_\sigma([a_2,0])$ where the points 
$\Upsilon(B^{\circ q_{m+1}-q_m}(1)), a_1, a_2$, $0$ form a $K_2$-bounded 
configuration
with $K_2,\gamma>0$ and $\sigma>\pi/2>\alpha$ independent of $m$.

The pull-back of $G''_m$ to $\check G_m$ is univalent.
Applying the Schwarz Lemma we have $\check G_m\subset G_m\cup G_\gamma([0,
\Upsilon(B^{ -q_{m+1}+q_m}(a_1))])$ and the claim follows.
\end{pf}
\realfig{ind1fig}{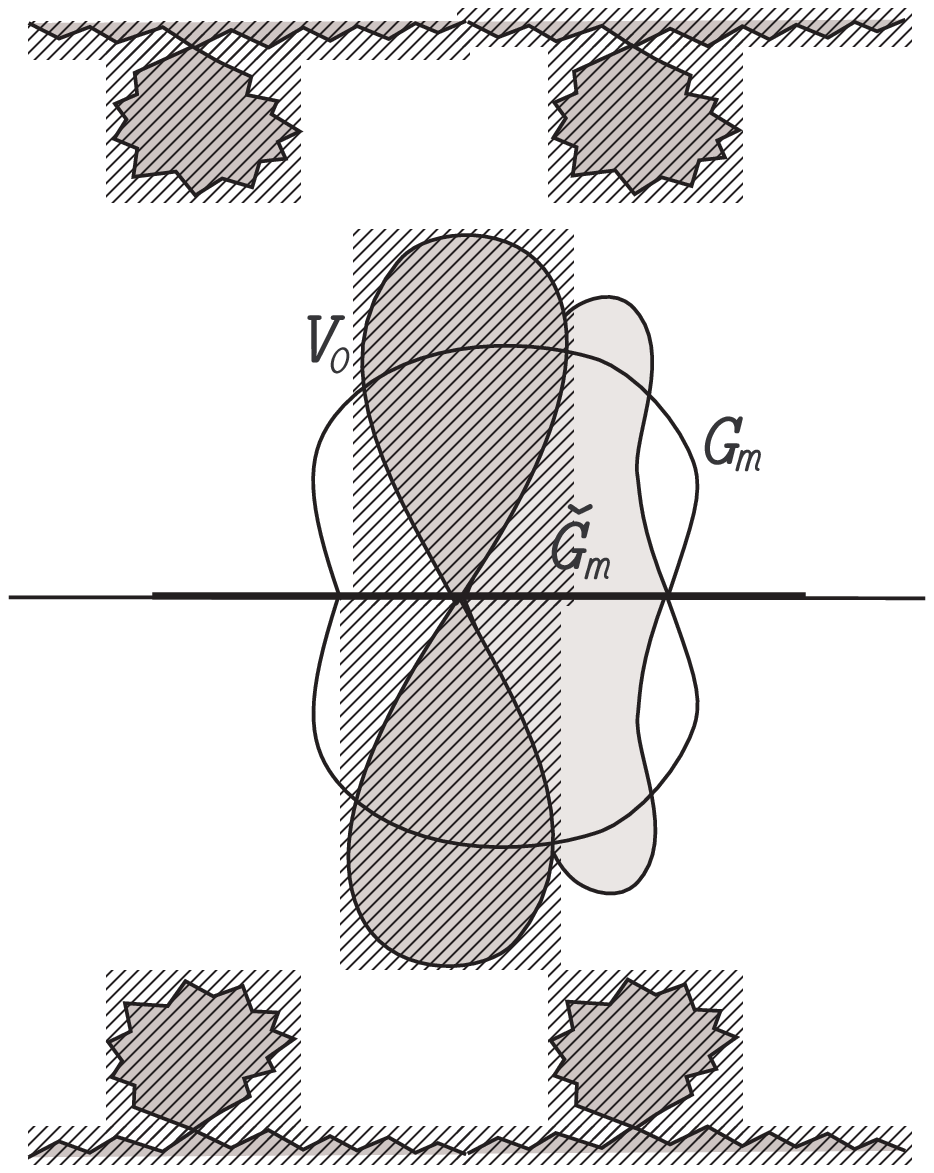}{}{6.5cm}

\begin{lem}
\label{ind2}
Let $J$ be the last return of the orbit (\ref{J orbit}) to the interval
$I_m$ preceding the first return to $I_{m+1}$ for $n-1>m>1$, and let $J'$ and $J''$
be the first two returns to $I_{m+1}$. Let $\zeta$, $\zeta'$ and $\zeta''$
be the corresponding points in the inverse orbit (\ref{z orbit}), so that
$\zeta'=\phi^{\circ q_m}(\zeta)$, $\zeta''=\phi^{\circ q_{m+2}}(\zeta')$.
Suppose that $\zeta\in G_m$. Then either $\ang{\zeta''}{I_{m+1}}>\eps=\eps(B)>0$ and 
$\dist(\zeta'',J'')<C(B)|I_{m+1}|$, or $\zeta''\in G_{m+1}$.
\end{lem}
\begin{pf}
Note that $J\subset \Upsilon([B^{\circ q_{m+1}+q_m}(1),B^{\circ q_m}(1)])$.
By the Schwarz Lemma, $$\zeta'\in G_\beta(\Upsilon([B^{\circ q_{m+1}-q_m}(1),1]))$$
for some $\beta>0$ independent of $m$.
Denote by $\hat J$ and $\check J$ the intervals of (\ref{J orbit})
such that $\phi^{\circ q_{m+1}-q_m}(J')=\hat J$ and $\phi^{\circ q_m}(\hat J)=\check J$,
and let $\hat \zeta$, $\check\zeta$ be the corresponding points of (\ref{z orbit}). We have $\hat J\subset \Upsilon([B^{\circ q_m}(1),B^{\circ q_m-q_{m+1}}(1)])$ and 
$\hat \zeta\in G_\beta(\Upsilon([B^{\circ q_m}(1),B^{\circ q_m-q_{m+1}}(1)]))$.
By the Schwarz Lemma and elementary properties of the map $B$ 
(see (\ref{cubic estimate})), there exist points $b_1$, $b_2$ in $\Upsilon([1,
B^{ -q_{m+1}}(1)])$ such that $0$, $b_1$, $b_2$, $\Upsilon(B^{ -q_{m+1}}(1))$
form a $K$-bounded configuration, and
$$\check\zeta\in G_\sigma([0,b_1])\cup G_\gamma([b_2,\Upsilon(B^{ -q_{m+1}}(1))])$$
for $\sigma$ and $\gamma$ independent of $m$ and $\sigma>\pi/2$.
The claim now follows from the Schwarz Lemma.
\end{pf}

Let us now select a strip $\breve S\Subset S$ used in \lemref{good angle}. By \lemref{leaves delta}
there exists $N>0$ such that $P_n\cap \bd U^\infty=\emptyset$ for all $n\geq N$.
Let $E$ be an annulus around the unit circle, compactly contained in the domain $\CC\sm (\ov{U}^\infty\cup \ov{U}^0)$
and such that $P_N\cup P_{N+1}\subset E$. We set $\breve S$ to be the strip $\frac{1}{2\pi i}\log (E)$.
Let $\hat P_n$ denote the component of $\frac{1}{2 \pi i}\log(P_n)$ attached to
$\Upsilon([1,B^{ -q_n}(1)])$.
Our argument culminates in the next lemma:
\begin{lem}
\label{ind3}
As before let $P_n$ denote the $n$-th critical puzzle piece and $N$ be as above.
Then for all $n\geq N+3$ we have
\begin{equation}
\label{cube root}
\diam P_{n}\leq C_1\sqrt[3]{\frac{\diam P_{n-1}}{|[B^{\circ q_{n-1}}(1),1]|}}\cdot |[1,B^{ -q_n}(1)]|+C_2
\end{equation}
for positive constants $C_1$, $C_2$ independent of $n$.
Moreover, for $z\in  \hat P_{n}$, either $z\in G_{n-1}$ or $\ang{z}{I_{n-1}}>\eps >0$, where $\eps$ is again independent of $n$.
\end{lem}

\begin{pf}
Choose $\alpha>0$ in the definition of $G_n$ so that 
 $$\hat P_{N+2}\cup \hat P_{N+3}\subset G_\alpha(\Upsilon([B^{\circ q_{N+2}}(1),B^{q_{N+1}-q_{N+2}}(1)]))=G_{N+1}.$$ 
By \corref{puzzle-pieces nested}, $P_{n+2}\subset P_n$ for all $n$, hence $\hat P_n\subset G_{N+1}$
for all $n\geq N+3$.
Fix a value of  $n>N+4$. Let
\begin{equation}
\label{pi orbit}
\Pi_0= \hat P_{n-1}, \Pi_{-1},\ldots,\Pi_{-q_n}= \hat P_{n}
\end{equation}
be the inverse orbit corresponding to the orbit (\ref{J orbit}).
We begin by establishing
\begin{equation}
\label{linear expansion}
\frac{\diam \Pi_{-(q_n-1)}}{|J_{-(q_n-1)}|}\leq K_1 \frac{\diam \hat P_{n-1}}{|J_0|}
\end{equation}
for some constant $K_1$ which does not depend on $n$.

Let $z\in \bd \hat P_{n-1}$ and consider the inverse orbit (\ref{z orbit}). Let $m\leq n$ be the largest
value for which $z\in G_m$. We will prove the estimate (\ref{linear expansion}) using an  induction on $m$. 
Let $T_{-1},\ldots,T_{-k}$
be the consecutive returns of the orbit of $J_0$ as (\ref{J orbit}) to $I_m$ until the first return
to $I_{m+1}$, and let $\zeta_{-1},\ldots,\zeta_{-k}$ be the corresponding points
in (\ref{z orbit}). 
Note that by \'Swia\c\negthinspace tek-Herman real a priori bounds, the intervals $T_{-i}$ are all $K$-commensurable
with $J_0$, for some $K$ independent of $n$.
It is easily seen from the combinatorics of the closest returns that the elements $\Pi_{-k_i}$ of the inverse orbit
(\ref{pi orbit}) corresponding to the points $\zeta_{-i}$ intersect the real axis along 
a subset of $(\hat P_N\cup \hat P_{N+1})\cap \RR$. 
By \lemref{nesting property}, $\Pi_{-k_i}\subset \hat P_N\cup \hat P_{N+1}$,
so $\zeta_{-i}\in\breve S$.
By \lemref{ind1}, either there exists a moment $i$ between $0$ and $k$
such that 
$$\ang{\zeta_{-i}}{I_m}>\eps\text{ and }\dist(\zeta_{-i},T_{-i})<C|I_m|,$$
or $\zeta_{-k}\in G_m$. In the former case we derive (\ref{linear expansion}) from
 \lemref{good angle}. In the latter case, consider the point $\zeta''$ which corresponds to 
the second return of (\ref{J orbit}) to $I_{m+1}$. By \lemref{ind2}, either $\ang{\zeta''}{I_{m+1}} > \eps$
and $\dist(\zeta'',I_{m+1})<C|I_{m+1}|$, or $\zeta''\in G_{m+1}$. 

In the first case we are done again by \lemref{good angle}. In the second case the proof of
 (\ref{linear expansion}) is completed by induction on $m$.
The same argument implies that either $\ang{z_{-q_n}}{J_{-q_n}}>\eps$, or $z_{-q_n}\in G_{n-1}$.
The estimate (\ref{cube root}) follows from (\ref{linear expansion}) and (\ref{cubic estimate}).
\end{pf}

The estimate (\ref{cube root}) implies that if $\displaystyle\frac{\diam P_{n-1}}{|[B^{\circ q_{n-1}}(1),1]|}>K$
for a large $K>0$, then 
$$1< \frac{\diam P_n}{|[1,B^{  -q_n}(1)]|}<\frac{1}{2}\cdot\frac{\diam P_{n-1}}{|[B^{\circ q_{n-1}}(1),1]|}.$$
This implies that for large $n$ the puzzle-piece $P_n$ is commensurable with its base arc
$[1,B^{  -q_n}(1)]$. In combination with the previous lemma, this shows that 
$P_n\subset G_\sigma(\Upsilon(I_{n-1}))$ for some fixed $\sigma>0$. 
Applying the Schwarz Lemma to the inverse orbit
$$P_n, B^{\circ q_{n+1}-q_n}(P_{n+1}),B^{\circ q_{n+1}-2q_n}(P_{n+1}),\ldots,B^{\circ q_{n-1}}(P_{n+1}),P_{n+1},$$
we see that

\begin{cor}
\label{in nbhd}
There exists an angle $\gamma>0$ such that for large values of $n$, $$\hat P_{n+1}\subset G_\gamma
(\Upsilon([1,B^{  -q_{n+1}}(1)])).$$
\end{cor}

Let us summarize the consequences. We first prove the following:

\begin{lem}[Only two drop-chains]
\label{biaccessible}
There are exactly two drop-chains of the form 
$\cD_1=\ov{\bigcup_k U^\infty_{\iota_1\ldots\iota_k}}$ and $\cD_2=\ov{\bigcup_k U^\infty_{\iota'_1\ldots\iota'_k}}$
accumulating at the critical point $1$. Moreover, both of these drop-chains land at $1$, and they separate $U_1$ from $\D$, in the sense that $U_1$ and $\DD$ belong to different components of $\ov{\CC} \sm (\cD_1 \cup \cD_2)$. 
\end{lem}
\begin{pf}
Let $\cD=\ov{\bigcup_k U^\infty_{\iota_1\ldots\iota_k}}$ be any drop-chain accumulating at $1$. This implies
that for an arbitrarily large $n$ there is a drop $U^\infty_{\iota_1\ldots\iota_k}\subset \cD$ which intersects
the critical puzzle-piece $P_n$. Since $U^\infty_{\iota_1\ldots\iota_k}$ cannot intersect $\bd P_n$, 
$U^\infty_{\iota_1\ldots\iota_k}\subset P_n$. By \lemref{limb inside}, the whole limb
$L^\infty_{\iota_1\ldots\iota_k}$ is contained in $P_n$.
By \corref{in nbhd}, $\diam P_n\to 0$, hence the drop-chain $\cD$ lands at $1$.

By \lemref{common boundary} every puzzle-piece $P_n$ contains a drop at infinity $U^\infty_{\iota_1\ldots\iota_k}$.
Since $P_{n+2}\subset P_n$ (\corref{puzzle-pieces nested}) and $\overset{\circ}{P}_n\cap\overset{\circ}{P}_{n+1}=\emptyset$,
there exist at least two distinct drop-chains landing at $1$ (passing through $P_n$'s with even and odd $n$'s respectively).
Clearly these drop-chains  separate $U_1$ from $\DD$.

Assume that there is a third drop-chain landing at $1$. This implies that there are two distinct drop-chains landing at the
critical value $B(1)$. Then the complement of the union of these drop-chains has a component $O$ which does not 
 contain any of the drops $U_i$. This implies that $O\subset \bigcup B^{-n}(U^\infty)$, which is a contradiction.
\end{pf}

The above lemma implies that for every $i\geq 1$ there are exactly two drop-chains  $\cD_1^{i}$, $\cD_2^{i}$
 accumulating at the
 point $x_{i}=B^{-i+1}(1)\in\TT$. These drop-chains land at $x_{i}$ and separate $U_{i}$ from $\D$. We may now define, as in subsection \ref{subsec:Limbs and wakes}, the {\it wake with root $x_{i}$} to be the 
  the connected component  $W_{i}$ of $\ov{\C}\sm (\cD_1^{i}\cup \cD_2^{i})$
containing $U_{i}$. 
For the corresponding limb we clearly have $L_i\subset \ov{W}_i$.
Due to the symmetry of the surgery (\corref{symmetry}), all the objects we have defined have their symmetric counterparts.
That is there is a sequence of critical puzzle-pieces $P_n^\infty$ converging to the critical point $c\in\bd U^\infty$,
wakes $W^\infty_i\supset U^\infty_i$ with $L^\infty_i\subset \ov{W}_i^\infty$, etc.

We now proceed to give the proof of \thmref{drops shrink}, which will occupy the rest of the section.

\smallskip
\noindent
{\it Proof of \thmref{drops shrink}.}
Let $\cD=\ov{\bigcup_k U^\infty_{\iota_1\ldots\iota_k}}$ be a drop-chain accumulating at a point $z\in J(\tilde{B})$.
We would like to show that $\diam L^\infty_{\iota_1\ldots\iota_k}\to 0$, which in turn will imply that $\cD$ lands at $z$. By symmetry of the surgery (\corref{symmetry}) this will prove the desired statement. Denote by $z_i$ the forward iterate $B^{\circ i}(z)$. Let us consider the two possibilities:

$\bullet$ {\it Case 1.} There exist $n$ and $m$ such that for $i>m$, $z_i\notin P_n\cup P_{n+1}\cup P^\infty_n\cup P^\infty_{n+1}$.
Let $\zeta$ be a limit point of the sequence $\{z_i\}$.
Since the rotation numbers $\theta$, $\nu$ are irrational, our assumption implies that $\zeta\notin\TT\cup\bd U^\infty$.
Clearly, the point $\zeta$ must be contained in a wake at infinity, which we call $W^\infty_j$. Denote by $i_k$ the moments $z_{i_k}\in L^\infty_j$,
and by $\Omega_k$ the univalent pull-back of $W^\infty_j$ along the orbit $z,z_1,\ldots,z_{i_k}$.
We refer to the following lemma to show that $\diam(\Omega_k)\to 0$ as $k \to \infty$ (see for example \cite{L3}, Prop. 1.10):
\begin{lem}[Shrinking Lemma]
Let $F$ be a rational map. Let $\{F^{  -m}_i\}$ be a family of 
univalent branches of the inverse maps in a domain $U$.
If $U\cap J(F)\neq \emptyset$, then for any $V$ such that $\ov{V}\subset U$, we have $\diam(F^{  -m}_i V)\to 0$ as $m \to \infty$.
\end{lem}

\noindent
Applying this lemma to our situation, we conclude that $\diam \Omega_k \to 0$. 
A drop $U^\infty_{\iota_1\ldots\iota_k}$ does not intersect the boundary of $\Omega_k$. Moreover, 
by the same argument as in \lemref{limb inside}, if a drop  $U^\infty_{\iota_1\ldots\iota_k}$ is contained in $\Omega_k$,
then $L^\infty_{\iota_1\ldots\iota_k}\subset \Omega_k$.
Thus the diameters of the limbs $L^\infty_{\iota_1\ldots\iota_k}$ shrink to zero,
 and hence the drop-chain $\cD$ lands at $z$.

$\bullet$ {\it Case 2}. To fix the ideas, let us assume that the critical point $1$ is a limit point of the sequence $\{z_i\}$.
Let $z_{i_n}$ be the first point in the orbit $\{ z_i \}$ contained in the puzzle-piece $P_n$.
Denote by 
\begin{equation}
\label{Q orbit}
Y_0= P_n,Y_{-1},\ldots, Y_{-{i_n}}
\end{equation}
the univalent preimages of $P_n$ along the inverse orbit $z_{i_n},\ldots,z$.

\begin{lem}
\label{at most once}
There exist at most one moment $i$ between $1$ and $i_n$ such that element $Y_{-i}$ of the inverse 
orbit (\ref{Q orbit}) hits the critical point $1$. Moreover, the pull-back (\ref{Q orbit}) decomposes
into two maps with bounded distortion and, possibly, one iterate of $B^{-1}$ near the 
critical value.
\end{lem}
\begin{pf}
Let us prove the first statement.
To be definite let us assume that $P_n$ is above  the critical point $1$.
Note that if $Y_{-i}\cap{\TT}=\emptyset$ for some $i\leq q_{n+1}$, then the
inverse orbit (\ref{Q orbit}) never hits the critical point for $1<i\leq i_n$. Otherwise denote
by $A$ and $B$ the ``above'' and ``below'' $B^{\circ q_{n+1}}$-preimages of $P_n$. 
One verifies directly, using the observations made in \lemref{puzzle properties} that $A\cap(\TT)\subsetneqq
 P_n\cap(\TT) $
 (compare \cite{Yampolsky}, Lemma 6.11). By \lemref{nesting property}, $A\subset P_n$, and thus $z_{q_{n+1}}\notin A$. 
The next possible moment when (\ref{Q orbit}) hits $1$ is $i=q_{n+1}+q_n$.
However, if $Y_{-q_{n+1}-q_n}\cap {\TT}\ne\emptyset$, then we may verify again that $Y_{-q_{n+1}-q_n}\subset P_n$,
which is not possible by our assumption.

Now let $k\leq i_n$ be the last  moment when $Y_{-k}\cap \TT\ne\emptyset$. As seen from 
the above argument, in combination with \'Swia\c\negthinspace tek-Herman real a priori bounds and \corref{in nbhd},
the pull-back $Y_0\to\cdots\to Y_{-k}$ decomposes into two maps with bounded distortion and, possibly,
one branch of $B^{-1}$ near the critical value. The combinatorics of closest returns and real a priori bounds
also imply that $\dist(Y_{-k},B(1))$ is 
greater than $K_1\diam Y_{-k}$ for some constant $K_1>0$. Hence the distance from $Y_{-k-1}$ 
to $\TT\cup \bd U^\infty$ is greater than $K_2\diam Y_{-k-1}$ for $K_2>0$, and the rest of the pull-back
$Y_{-k}\to\cdots\to Y_{-i_n}$ has bounded distortion by the Koebe Theorem.
\end{pf}

\noindent
By \lemref{commensurable disk} and \corref{in nbhd} the puzzle-piece $P_n$ contains a Euclidean disk,
whose diameter is commensurable with $\diam P_n$, centered at a point in $J(B)$.
Therefore, by \lemref{at most once}, the domain
$ Y_{-{i_n}}\ni z$ contains a Euclidean disk centered at a point of $J({B})$ whose diameter is
commensurable with $\diam  Y_{-{i_n}}$. This implies that $\diam  Y_{-{i_n}}\to 0$. 
By \lemref{limb inside}, if $U^\infty_{\iota_1\ldots\iota_k}\subset Y_{-{i_n}}$, then $L^\infty_{\io_1\ldots\io_k}\subset Y_{-{i_n}}$.
So the diameters of limbs $L^\infty_{\io_1\ldots\io_k}$ shrink to zero, and the drop-chain $\cD$ lands at $z$. \hfill $\Box$

\vspace{0.17in}

\section{The Proof}
\label{sec:proof}

Throughout this section we fix a pair of irrationals $\theta$ and $\nu$ of bounded type, with $\theta\ne 1-\nu$. 
In what follows we prove the Main Theorem, that is we show that
 the quadratic rational map $\FTN$ of (\ref{eqn:normal form}) is in fact {\it the} mating
of the quadratic polynomials $\FT$ and $\FN$ in the sense we described in the introduction.

\subsection{Spines and itineraries.}
\label{subsec:Spines and itineraries}
Let $\TQT$ be the modified Blaschke product $\TQT$ of (\ref{eqn:modify}). Consider the two drop-chains
$$\CCC=\ov{U_0 \cup U_1 \cup U_{11} \cup \cdots},\ \ \CCC'=\ov{U_0 \cup U_2 \cup U_{21} \cup \cdots}$$ with
$\TQT(\CCC')=\CCC$. 
Applying \lemref{snail} again, we see that  $\CCC$ and $\CCC'$
land respectively at the repelling fixed point $\beta$ and its preimage $\beta'$. By the {\it spine} of
$\TQT$ we mean the union of the drop-rays
$$\ST=R(\CCC) \cup R(\CCC')$$ (compare \figref{angle}, 
where the image of the spine of $\TQT$ is shown in the filled Julia set of 
the quadratic polynomial $\FT$ for $\theta=(\sqrt{5}-1)/2$). 
Every point on the spine which is not in the interior of $K(\TQT)$ is either one of the
endpoints $\beta$, $\beta'$, or a preimage of the critical point $z=1$.

By Petersen's  \thmref{diam->0} the Julia set $J(\TQT)$ is locally-connected. 
Thus the B\"ottcher map extends continuously from  the basin of infinity of $\TQT$
to its boundary.
As a consequence, there exists a Carath\'eodory loop
 $\GT: \RR/ \ZZ \to J(\TQT)$  which conjugates the angle-doubling map to $\TQT$.
 A point $z\in J(\TQT)$ is the landing point of
an external ray $R^e(t)$ if and only if $\GT(t)=z$. It is easy to see that $\GT(0)=\beta$ and $\GT(1/2)=\beta'$. 

By \lemref{tworays} the critical point $z=1$, hence every preimage of it, is {\it biaccessible}, that is it is
the landing point of exactly two external rays. For the quadratic polynomial $\FT$ the converse statement is true for an arbitrary $\theta$ of Brjuno type:  Every biaccessible point in the
Julia set $J(\FT)$ eventually hits the critical point \cite{Zakeri1}.
The two external rays landing at the critical point of $\TQT$ are both mapped to the external ray landing at the critical value $\TQT(1)$. This means that they have angles of the form  $\omega/2$ and
$(\omega+1)/2$, where $\omega=\omega(\theta)$ is a well-defined 
irrational number in the interval $(0,1)$. It can be shown that the 
function $\theta \mapsto \omega(\theta)$ is effectively computable (see \cite{Bullett} and compare with subsection \ref{subsec:Rotation sets of the doubling map}).

Consider the two connected subsets of the Julia set
\begin{equation}
\label{eqn:pieces}
\begin{array}{rl}
J_{\theta}^0= & \{ z\in J(\TQT) : z=\GT(t)\ \mbox{for some}\ 0\leq t<1/2 \},\\ J_{\theta}^1= & \{ z\in
J(\TQT) : z=\GT(t)\ \mbox{for some}\ 1/2 \leq t<1 \}.
\end{array}
\end{equation}
By local-connectivity of $J(\TQT)$ (\thmref{diam->0}), $J_{\theta}^0 \cup J_{\theta}^1= J(\TQT)$, and
evidently $J_{\theta}^0 \cap J_{\theta}^1=(\bigcup_{n=0}^{\infty} \TQT^{-n}(1))\cap \ST= \{ 1=x_1, x_{11},
x_{111}, \ldots \} \cup \{ x_2, x_{21}, x_{211}, \ldots \}$.

We proceed to define the {\it itinerary} of a point $z\in J(\TQT)$ with respect to $\ST$. This will be
a dynamically-defined infinite sequence of $0$'s and $1$'s which gives the binary expansion of the angle of
an external ray landing at $z$ (see \cite{Douady3} for a general dicscussion on how one computes angles in
similar situations). In the case where $z$ is biaccessible, we define two different itineraries corresponding to the
angles of the two rays landing at $z$. Set $z_0=z,\ z_k=\TQT(z_{k-1})$. We consider three distinct cases:

$\bullet$ {\it Case 1.} The orbit of $z$ never hits the spine. Then $z$ is not biaccessible and hence there
exists a unique angle $t$ with $z=\GT(t)$. Define the itinerary of $z$ to be  the sequence
$\varepsilon=(\varepsilon_0, \varepsilon_1, \varepsilon_2, \ldots )$, where $\varepsilon_i \in \{ 0, 1\}$
is determined by the condition
$$z_i \in J_{\theta}^{\varepsilon_i}, \ \ \ i=0,1,\ldots$$ Then it is easy to see that the angle $t$ has
binary expansion $0.\varepsilon_0 \varepsilon_1 \varepsilon_2 \cdots$.

$\bullet$ {\it Case 2.} The orbit of $z$ eventually hits the $\beta$-fixed point, i.e., there exists the
smallest integer $n\geq 0$ such that $z_n=\beta$. In this case, the angle $t$ with $z=\GT(t)$ is still
unique. The itinerary of $z$ is defined as $\varepsilon=(\varepsilon_0, \varepsilon_1, \ldots ,
\varepsilon_n, 0,0,0, \ldots )$, where $$z_i \in J_{\theta}^{\varepsilon_i}, \ \ \ i=0,1,2,\ldots, n.$$ The
binary digits of the angle $t$ are then given by the itinerary of $z$.

$\bullet$ {\it Case 3.} The orbit of $z$ eventually hits the critical point at $1$. In this case there are
exactly two angles $0<t<s<1$ with $\GT(t)=\GT(s)=z$. Let us assume that the angles corresponding to the
critical point have binary expansions $\omega/2=0.0\omega_1 \omega_2 \ldots$ and $(\omega+1)/2=0.1\omega_1
\omega_2 \ldots$. Then the critical value $v=\TQT(1)$ has a unique ray landing on it with angle
$\omega=0.\omega_1 \omega_2 \ldots$. Since $v$ can never hit the spine, by {\it Case 1} above, the binary
digits of $\omega$ are uniquely determined by the condition
$$\TQT^{\circ i}(1)\in J_{\theta}^{\omega_i}, \ \ \ i=1,2,\ldots$$ We are going to define two itineraries
for $z$. Let $n\geq 0$ be the smallest integer such that $z_n \in \ST \sm \{ \beta, \beta' \}$. Define the
{\it initial segment} $(\varepsilon_0, \ldots , \varepsilon_{n-1})$ of both itineraries of $z$ by the
condition
$$z_i \in J_{\theta}^{\varepsilon_i}, \ \ \ i=0,1,\ldots, n-1.$$ (If $n=0$, we define the initial segment
to be empty.) Let $m\geq 1$ be defined by the condition $z_{n+m}=v=\TQT(1)$. Since $z_n, \ldots, z_{n+m-1}$
all belong to the intersection $J_{\theta}^0 \cap J_{\theta}^1$, there is an ambiguity in assigning digits
to the points of this part of the orbit of $z$. So consider $z_n$ and replace it by two points $a_n\in
J_{\theta}^0$ and $b_n \in J_{\theta}^1$, both sufficiently close to $z_n$. It is easy to see that the
points of the orbits $a_n, \ldots, a_{n+m-1}$ and $b_n, \ldots, b_{n+m-1}$ have well-defined itineraries
$(\varepsilon_n, \ldots , \varepsilon_{n+m-1})$ and $(\varepsilon'_n, \ldots , \varepsilon'_{n+m-1})$
determined by the conditions
$$a_i \in J_{\theta}^{\varepsilon_i}, \ \ \ i=n, n+1,\ldots, n+m-1,$$
$$b_i \in J_{\theta}^{\varepsilon'_i}, \ \ \ i=n, n+1, \ldots, n+m-1.$$ We call these two segments {\it
ambiguous}. Note that $\varepsilon_i+\varepsilon'_i=1$ for $n \leq i \leq n+m-1$. Finally, follow these two
by the well-defined itinerary of the critical value. We thus obtain two itineraries for $z$:
$$\varepsilon=(\underbrace{\varepsilon_0, \ldots, \varepsilon_{n-1}}_{\tiny \mbox{initial segment}},\
\underbrace{\varepsilon_n, \ldots, \varepsilon_{n+m-1}}_{\tiny \mbox{ambiguous segment}},\
\underbrace{\omega_1, \omega_2, \ldots}_{\tiny \mbox{itinerary of $v$}} ),$$
$$\varepsilon'=(\underbrace{\varepsilon_0, \ldots, \varepsilon_{n-1}}_{\tiny \mbox{initial segment}},\
\underbrace{\varepsilon'_n, \ldots, \varepsilon'_{n+m-1}}_{\tiny \mbox{ambiguous segment}},\
\underbrace{\omega_1, \omega_2, \ldots}_{\tiny \mbox{itinerary of $v$}} ).$$ These two itineraries give the
binary digits of the two angles $t$ and $s$.

Since $\TQT$ and $\FT$ are quasiconformally conjugate for $\theta$ of bounded type, with the conjugacy
being conformal in the basin of infinity, we have a completely similar description for the spine and
itineraries of points in the quadratic Julia set $J(\FT)$. \figref{angle} shows the spine and selected rays for $\FT$
with $\theta=(\sqrt{5}-1)/2$.
\realfig{angle}{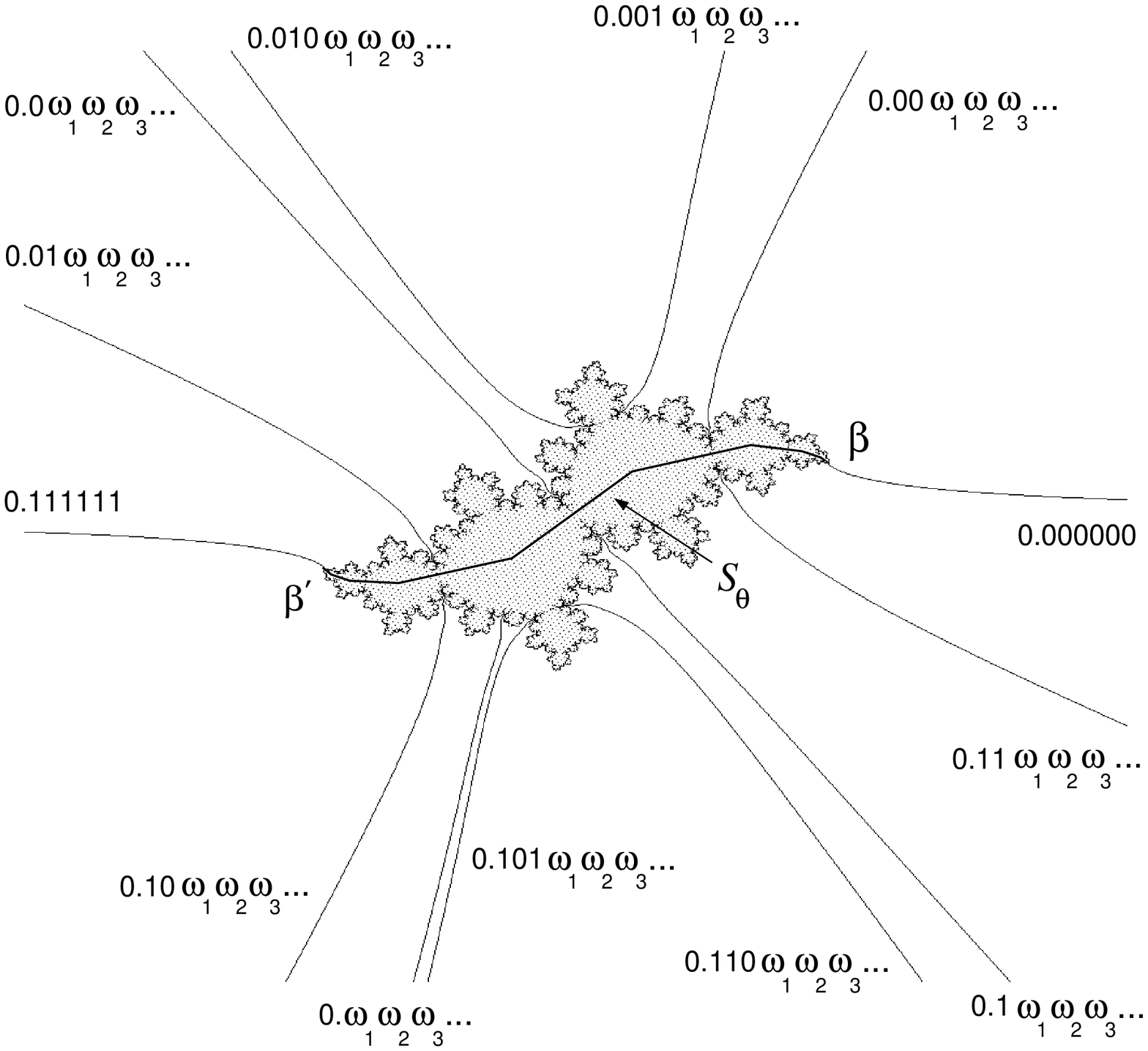}{{\sl This picture shows the filled Julia set of the quadratic polynomial $\FT$, for $\theta=(\sqrt{5}-1)/2$. The spine is shown by a thick path connecting the repelling fixed point $\beta$ to its preimage $\beta'$. Selected rays and angles in base 2 are shown. Here $\omega=0.\omega_1 \omega_2 \omega_3 \ldots$ is the unique angle corresponding to the ray which lands at the critical value. For this value of $\theta$, $\omega$ is given by the continued fraction $[1,2,2,2^2,2^3,2^5, \ldots]$, where the powers of $2$ form the Fibonacci sequence. Hence $\omega_1=1$, $\omega_2=0$, $\omega_3=1$, etc.}}{11cm}

We summarize the above discussion in the following proposition:

\begin{proposition}
\label{itin}

\noindent
\begin{enumerate}
\item[(i)]
Let $z\in J(\TQT)$. Then the angle(s) of the external ray(s) landing at $z$ is(are) determined by the
itinerary(ies) of $z$, that is by the answer to the purely topological question  of whether a point in the
forward orbit of $z$ belongs to $J_{\theta}^0$, $J_{\theta}^1$, or to which point of the spine. In
particular, two points in the Julia set having the same itinerary must coincide.
\item[(ii)]
Every infinite sequence of $0$'s and $1$'s can be realized as the itinerary of a unique point in $J(\TQT)$.
\end{enumerate} 
\end{proposition}

\subsection{Main reduction.}
\label{subsec:Main reduction}
A key ingredient in the proof of the main theorem is the following reduction step:

\begin{theorem}
\label{reduction}
Let $0< \theta, \nu <1$ be irrationals of bounded type and $\theta \neq 1-\nu$. Then there exist continuous
maps $\ZT: K(\TQT) \rightarrow \ov{\CC}$ and $\ZN: K(\TQN)\rightarrow \ov{\CC}$ such that
\begin{equation}
\label{eqn:semi}
\begin{array}{lll}
\ZT \circ \TQT & = \TBTN \circ \ZT & \mbox{on}\ K(\TQT)\\
\ZN \circ \TQN & = \TBTN \circ \ZN & \mbox{on}\ K(\TQN).  
\end{array}
\end{equation}
$\ZT$ and $\ZN$ can be chosen to be quasiconformal homeomorphisms in the interiors of $K(\TQT)$ and
$K(\TQN)$ respectively. Moreover, $\ZT(K(\TQT))\cup \ZN(K(\TQN))=\ov{\CC}$ and $\ZT(z)=\ZN(w)$ if and only
if there exists an angle $t \in \RR /\ZZ$ such that $z=\GT(t)$ and $w=\GN(-t)$.
\end{theorem} 

Before starting the proof, we fix some notation. For simplicity, we set
$K(\TQT)=\KT$, $K(\TQN)=\KN$. We also recall the definition of the compact set $K(\TBTN)=\KTN$ as the set
of all points whose forward orbits under the iteration of $\TBTN$ never hit the Siegel disk
$U^\infty$. Similarly, $\KTN^\infty=\ov{\CC\sm\KTN}$ is the set of points whose forward orbits never hit the 
 ``Siegel disk'' $U_0=\DD$. \vspace{0.12in}\\
{\it Proof of \thmref{reduction}.} We begin by constructing $\ZT$. The map $\ZN$ can be constructed in a
similar fashion. Consider the modified Blaschke products $\TQT$ of (\ref{eqn:modify}) and $\TBTN$ of
(\ref{eqn:modify2}). Since both of these are quasiconformally conjugate to the rigid rotation $z\mapsto
e^{2 \pi i \theta}z$ on the unit disk, one can define a quasiconformal conjugacy $\ZT: \DD \rightarrow \DD$
between them, which extends homeomorphically to a conjugacy $\ZT:\ov{\DD}\rightarrow \ov{\DD}$. This $\ZT$
can be extended to the union of the closures of all drops of $\TQT$ by  pulling back. To this
end, let $U_{\io_1 \ldots \io_k}$ be any drop of $\TQT$ of generation $k$ and consider the corresponding
drop $U'_{\io_1 \ldots \io_k}$ of $\TBTN$ {\it with the same address}. Let $n=\io_1+ \cdots + \io_k$ and
define $\ZT: \ov{U}_{\io_1 \ldots \io_k}\iso \ov{U'}_{\io_1 \ldots \io_k}$ by
$$\ZT= \TBTN^{-n} \circ \ZT \circ \TQT^{\circ n}.$$ An easy induction on $n$ shows that $\ZT$ defined this
way is a conjugacy between $\TQT$ and $\TBTN$ on $\bigcup_k \bigcup_{\io_1, \ldots, \io_k} \ov{U}_{\io_1
\ldots \io_k}$ which is quasiconformal on the union $\bigcup_k \bigcup_{\io_1,\ldots, \io_k} U_{\io_1
\ldots \io_k}=\mbox{int}(\KT)$.

We would like to extend $\ZT$ to a continuous semiconjugacy $\KT \rightarrow \KTN$. By
\propref{either/or}, every point in $\KT$ is either in the closure of a drop or is the landing point of a
unique drop-chain. Since $\ZT$ is already defined on $\bigcup_k \bigcup_{\io_1, \ldots, \io_k}
\ov{U}_{\io_1 \ldots \io_k}$, it suffices to define it at the landing points of drop-chains of $\TQT$. Take
a drop-chain $\CCC=\ov{\bigcup_k U_{\io_1 \ldots \io_k}}$ which lands at $p$ and consider the corresponding
drop-chain of $\TBTN$, $\CCC'=\ov{\bigcup_k U'_{\io_1 \ldots \io_k}}$, whose drops have {\it the same}
addresses. By \thmref{drops shrink}, the diameter of $U'_{\io_1 \ldots \io_k}$ goes to zero as
$k\rightarrow \infty$, hence $\CCC'$ lands at a well-defined point $p'\in \KTN$. Define $\ZT(p)=p'$.

Evidently $\ZT$ defined this way has the property that for any limb $L_{\io_1 \ldots \io_k}$ of $\TQT$, the
image $\ZT(L_{\io_1 \ldots \io_k})$ is precisely the limb $L'_{\io_1 \ldots \io_k}$ of $\TBTN$ with the
same address. We would like to show that $\ZT$ is continuous as a map from $\KT$ into $\ov{\CC}$. Take a
point $p\in \KT$ and a sequence $p_n \in \KT$ converging to $p$. When $p$ belongs to the interior of $\KT$
continuity is trivial. So let us assume that $p\in \bd \KT$. By \propref{either/or}, we have two
possibilities:

$\bullet$ {\it Case 1.} $p$ is the landing point of a drop-chain $\CCC=\ov{\bigcup_k U_{\io_1 \ldots
\io_k}}$. Fix a multi-index $\io_1 \ldots \io_k$ and observe that $p$ belongs to the wake $W_{\io_1 \ldots
\io_k}$. Therefore, for $n$ large enough, $p_n\in L_{\io_1 \ldots \io_k}$, which implies $\ZT(p_n)\in
L'_{\io_1 \ldots \io_k}$. It follows that $\dist(\ZT(p),\ZT(p_n)) < \diam (L'_{\io_1 \ldots \io_k})$. Since
$\diam (L'_{\io_1 \ldots \io_k})\rightarrow 0$ as $k\rightarrow \infty$ by \thmref{drops shrink}, we have
$\ZT(p_n)\rightarrow \ZT(p)$ as $n\rightarrow \infty$.

$\bullet$ {\it Case 2.} $p$ belongs to the boundary of a drop $U_{\io_1 \ldots \io_k}$ of $\TQT$ of {\it
smallest} possible generation. It might be the case that $p$ is the root of a child $U_{\io_1 \ldots \io_k
\io_{k+1}}$ in which case $p=\bd U_{\io_1 \ldots \io_k} \cap \bd U_{\io_1 \ldots \io_k \io_{k+1}}$. If for
all sufficiently large $n$, $p_n$ belongs to $\ov{U}_{\io_1 \ldots \io_k}$ (or to $\ov{U}_{\io_1 \ldots
\io_k} \cup \ov{U}_{\io_1 \ldots \io_k \io_{k+1}}$ if $p$ is the root of $U_{\io_1 \ldots \io_k
\io_{k+1}}$), then $\ZT(p_n)\to \ZT(p)$ is immediate. Hence it suffices to prove the convergence in the
case $p_n \notin \ov{U}_{\io_1 \ldots \io_k}$ (or $p_n \notin \ov{U}_{\io_1 \ldots \io_k} \cup
\ov{U}_{\io_1 \ldots \io_k \io_{k+1}}$ if $p$ is the root of $U_{\io_1 \ldots \io_k \io_{k+1}}$). Under
this assumption, it follows from $p_n \to p$ that $p_n$ belongs to a limb $L(n)$ with root $x(n)\in \bd
U_{\io_1 \ldots \io_k}$ (or $x(n)\in \bd U_{\io_1 \ldots \io_k} \cup \bd U_{\io_1 \ldots \io_k \io_{k+1}}$
if $p$ is the root of $U_{\io_1 \ldots \io_k \io_{k+1}}$) such that $x(n) \to p$ as $n \to \infty$. Then
$\ZT(p_n)$ belongs to the limb $ L'(n)$ of $\TBTN$ with the same address as $L(n)$ whose root
$x'(n)=\ZT(x(n))$ converges to $\ZT(p)$ as $n \to \infty$. Since $\diam (L'(n)) \to 0$ by \thmref{drops
shrink}, we must have $\ZT(p_n)\rightarrow \ZT(p)$ as $n\to \infty$ as well. This finishes the proof of
continuity.

We can define $\ZN$ and prove its continuity in a similar way. It is clear from the above construction that
the semiconjugacy relations (\ref{eqn:semi}) hold and $\ZT(\KT)=\KTN$ and similarly $\ZN(\KN)=\KTN^\infty$.

It remains to prove the last property of $\ZT$ and $\ZN$. Consider the spines $\ST$ and $\SN$ for $\TQT$
and $\TQN$ as in subsection \ref{subsec:Spines and itineraries} and map them to get simple arcs
$\Sigma_{\theta}=\ZT(\ST)$ and $\Sigma_{\nu}=\ZT(\SN)$ (compare \figref{spine1}). Set
$$\Sigma=\Sigma_{\theta}\cup \Sigma_{\nu}.$$

\begin{lemma}
\label{not both}
Two simple curves $\Sigma_{\theta}$ and $\Sigma_{\nu}$ do not intersect except at the two end-points
$\beta$ and $\beta'$. Hence $\Sigma$ is a Jordan curve on the Riemann sphere.
\end{lemma}

\begin{pf}
Clearly the intersection $\Sigma_{\theta} \cap \Sigma_{\nu}$ is a subset of $\bd \KTN \cap \Sigma$. Every
point in the latter intersection is either $\beta$ or $\beta'$, or is a preimage of $1$ or $c$, where $c$
is the critical point of $\BTN$ on the boundary of $U^\infty$. Since $1$ and $c$ have disjoint forward
orbits, the conclusion follows.
\end{pf} 

Now consider the four connected sets
$$\LT^i=\ZT(J_{\theta}^i),\ \LN^i=\ZN(J_{\nu}^i)\ \ \ i=0,1,$$ where $J_{\theta}^i$ and $J_{\nu}^i$ are the
subsets of the Julia sets $J(\TQT)$ and $J(\TQN)$ we defined in (\ref{eqn:pieces}). Let
$$X=\{ 1=x_1, x_{11}, x_{111}, \ldots, x_2, x_{21}, x_{211}, \ldots \}$$ and
$$Y=\{ c=x_1^\infty, x_{11}^\infty, x_{111}^\infty, \ldots, x_2^\infty, x_{21}^\infty, x_{211}^\infty,
\ldots \}$$ be the preimages of the critical points $1$ and $c$ along $\Sigma$.  It is clear from the
definition that
$$X \subset \LT^0 \cap \LT^1 \subset X \cup Y,$$
$$Y \subset \LN^0 \cap \LN^1 \subset X \cup Y.$$ But in fact we have the following much sharper statement:

\begin{lemma}
\label{intersection}
With the above notation, we have
$$\LT^0 \cap \LT^1=\LN^0 \cap \LN^1= X \cup Y. $$
\end{lemma}

\begin{pf}
Take a point $y\in Y$ and assume that $\TBTN^{\circ n}(y)=c$. By \lemref{biaccessible}, there are exactly
two drop-chains which land at the critical point $c$ from different sides of $\Sigma_{\nu}$. Then the
pull-backs of these drop-chains along the orbit $y, \TBTN(y), \ldots, \TBTN^{\circ n}(y)=c$ give two
drop-chains which land at $y$ from different sides of $\Sigma_{\nu}$. These drop-chains are clearly subsets
of the compact set $\KTN$. The fact that they land at $y$ from different sides of $\Sigma_{\nu}$ implies
$y\in \LT^0 \cap \LT^1$. The proof of the other equality is similar.
\end{pf}

\begin{corollary}
\label{equal}
With the above notation, we have
$$\LT^0=\LN^1\ \ \ \mbox{and} \ \ \ \LT^1=\LN^0.$$
\end{corollary}

\begin{pf}
Let $\ov{\CC} \sm \Sigma= O_1 \cup O_2$, where $O_i$ are disjoint topological disks with $\LT^0 \subset
\ov{O}_1$ and $\LT^1 \subset \ov{O}_2$. Taking the orientations on the sphere into account, we have $\LN^1
\subset \ov{O}_1$ and $\LN^0 \subset \ov{O}_2$. Since $\LT^0 \cup \LT^1= \LN^0 \cup \LN^1= \bd
\KTN$ and $\LT^0 \cap \LT^1=\LN^0 \cap \LN^1$ by \lemref{intersection}, it follows that $\LT^0=\LN^1$ and
$\LT^1=\LN^0$.
\end{pf} 

We can now define the itinerary of a point $p\in \bd \KTN$ with respect to $\Sigma_{\theta}$ by looking at
the points in the forward orbit of $p$ and deciding whether they belong to $\LT^0$, $\LT^1$, or to $\LT^0
\cap \LT^1$. As in the discussion of itineraries for the points in the Julia set $J(\TQT)$ (see subsection
\ref{subsec:Spines and itineraries}), we may face an ambiguity in defining the digits when some forward
iterate of $p$, say $p_n$, belongs to $\LT^0 \cap \LT^1$. In this case, we perturb $p_n$ to obtain a pair
of nearby points in $\LT^0$ and $\LT^1$ and keep iterating the two points to decide to which piece of the
Julia set they belong. After a finite number of iterations, we are off the spine $\Sigma$ and the rest of
the itinerary can be defined in an unambiguous way. Since $\bd \KTN= \bd \KTN^\infty$ by \corref{common boundary},
a similar procedure can be used to define the itinerary or two itineraries of $p$ with respect to
$\Sigma_{\nu}$. In short,

\begin{proposition}[Two or four itineraries]
\label{two itin}
Let $p\in \bd \KTN$. Then, either $p$ is not a preimage of $1$ or $c$ in which case it has unique
itineraries $\varepsilon_{\theta}$ with respect to $\Sigma_{\theta}$ and $\varepsilon_{\nu}$ with respect
to $\Sigma_{\nu}$, or $p$ is a preimage of $1$ or $c$ in which case it has two different itineraries
$\varepsilon_{\theta}, \varepsilon'_{\theta}$ with respect to $\Sigma_{\theta}$ and $\varepsilon_{\nu},
\varepsilon'_{\nu}$ with respect to $\Sigma_{\nu}$.
\end{proposition}     

Since the $i$-th digit of the itinerary or itineraries of a point $p$ with respect to $\Sigma_{\theta}$ is
determined by the condition $\TBTN^{\circ i}(p)\in \LT^0$ or $\LT^1$, and similarly for the itineraries
with respect to $\Sigma_{\nu}$, we have the following consequence of \corref{equal}:

\begin{proposition}[$\Sigma_{\theta}$- and $\Sigma_{\nu}$-itineraries have opposite digits]
\label{01-10}
Let $p\in \bd \KTN$ have itinerary $\varepsilon_{\theta}(p)=(\varepsilon_0, \varepsilon_1, \varepsilon_2,
\ldots)$ with respect to $\Sigma_{\theta}$. Then the itinerary of $p$ with respect to $\Sigma_{\nu}$ is
obtained by converting all $0$'s to $1$'s and all $1$'s to $0$'s in $\varepsilon_{\theta}$. In other words,
$\varepsilon_{\nu}(p)=\ov{\varepsilon}_{\theta}(p)=(\ov{\varepsilon}_0, \ov{\varepsilon}_1,
\ov{\varepsilon}_2, \ldots)$, where $\ov{\varepsilon}_i=1-\varepsilon_i$. In the case where $p$ has two itineraries,
we have $\varepsilon_{\nu}(p)=\ov{\varepsilon}_{\theta}(p)$ and
$\varepsilon'_{\nu}(p)=\ov{\varepsilon'}_{\theta}(p)$.
\end{proposition}

The following lemma is a straightforward consequence of the above construction:

\begin{lemma}[Itineraries match]
\label{same}
Let $z\in \KT$ and $p=\ZT(z)\in \KTN$.
\begin{enumerate}
\item[(i)]
Suppose that $z$ is not a preimage of the critical point $1$. Then the unique itinerary of $z$ with respect
to $\ST$ coincides with $\varepsilon_{\theta}(p)$ when $p$ is not a preimage of $c$, and it coincides
with one of the two itineraries $\varepsilon_{\theta}(p)$ or $\varepsilon'_{\theta}(p)$ when $p$ is a
preimage of $c$.
\item[(ii)]
Suppose that $z$ is a preimage of $1$. Then the two itineraries of $z$ with respect to $\ST$ coincide with
the two itineraries $\varepsilon_{\theta}(p)$ and $\varepsilon'_{\theta}(p)$.
\end{enumerate} 
\end{lemma} 

\begin{corollary}[Itineraries determine points]
\label{determ}
Two points in $\bd \KTN$ with the same itinerary with respect to $\Sigma_{\theta}$ or $\Sigma_{\nu}$ must
coincide.
\end{corollary}

\begin{pf}
Let $p,q \in \bd \KTN$ and assume for example that $\varepsilon_{\theta}(p)=\varepsilon_{\theta}(q)$. When $p$ (hence $q$) is a preimage of $1$ or $c$, it is easy to see that identical $\Sigma_{\theta}$-itineraries implies $p=q$. So let us assume that $p$ and $q$ are not preimages of $1$ or $c$. Since
$\ZT: \KT \to \KTN$ is surjective, we have $p=\ZT(u)$ and $q=\ZT(v)$ for some $u,v \in \bd \KT=
J(\TQT)$. By the above lemma, $u$ and $v$ have the same itineraries with respect to $\ST$. By \propref{itin}(i), $u=v$. Hence $p=q$.
\end{pf} 

Now consider two points $z\in \KT$ and $w\in \KN$ such that  $z=\GT(t)$ and $w=\GN(-t)$ for some $t\in \TT$. Set
$p=\ZT(z)$ and $q=\ZN(w)$. 
The binary digits  $(\varepsilon_0, \varepsilon_1, \varepsilon_2, \ldots)$ of the angle $t$ form 
an  itinerary of $z$ with respect to $\ST$. Since $t=0.\varepsilon_0 \varepsilon_1
\varepsilon_2 \ldots $ in base $2$, $-t$ has  the binary expansion $0.\ov{\varepsilon}_0 \ov{\varepsilon}_1
\ov{\varepsilon}_2 \ldots $. Hence $(\ov{\varepsilon}_0, \ov{\varepsilon}_1, \ov{\varepsilon}_2,
\ldots)$ is an itinerary of $w$ with respect to $\SN$. Thus by \lemref{same},
$\varepsilon_{\theta}(p)=(\varepsilon_0, \varepsilon_1, \varepsilon_2, \ldots)$ and
$\varepsilon_{\nu}(q)=(\ov{\varepsilon}_0, \ov{\varepsilon}_1, \ov{\varepsilon}_2, \ldots)$. By
\propref{01-10}, $\varepsilon_{\theta}(q)=(\varepsilon_0, \varepsilon_1, \varepsilon_2, \ldots)$, which
means $p$ and $q$ have the same itinerary with respect to $\Sigma_{\theta}$. This, by \corref{determ},
implies $p=q$.

Conversley, assume that $\ZT(z)=\ZN(w)=p$. We consider two cases: First assume that $p$ is not a preimage of $1$ or $c$. Then it follows from \propref{01-10} that $\varepsilon_{\theta}(p)=\ov{\varepsilon}_{\nu}(p)=(\varepsilon_0, \varepsilon_1, \varepsilon_2, \ldots)$ and these itineraries are unique. By \lemref{same},  $(\varepsilon_0,
\varepsilon_1, \varepsilon_2, \ldots)$ is the $\ST$-itinerary of $z$ and $(\ov{\varepsilon}_0, \ov{\varepsilon}_1, \ov{\varepsilon}_2, \ldots)$ is the $\SN$-itinerary of $w$. Setting  $t=0.\varepsilon_0 \varepsilon_1 \varepsilon_2
\ldots $ in base $2$, we have $z=\GT(t)$ and $w=\GN(-t)$ and we are done. Next, assume that $p$ is a preimage of, say, $1$. Then, as $1$ and $c$ have disjoint orbits under $\TBTN$, $p$ cannot be a preimage of $c$. This implies that $z$ is a preimage of the critical point $1$ of $\TQT$ and therefore has two itineraries, and $w$ is not a preimage of the critical point $1$ of $\TQN$ and so has a unique itinerary. Let $w=\GN(-t)$, where the unique $t\in \TT$ has binary expansion $t=0.\varepsilon_0 \varepsilon_1 \varepsilon_2
\ldots $. By \lemref{same}, ${\varepsilon}_{\nu}(p)=(\ov{\varepsilon}_0, \ov{\varepsilon}_1, \ov{\varepsilon}_2, \ldots)$ is one of the $\Sigma_{\nu}$-itineraries of $p$. Hence by \propref{01-10}, $\varepsilon_{\theta}(p)=(\varepsilon_0, \varepsilon_1, \varepsilon_2, \ldots)$ is one of the $\Sigma_{\theta}$-itineraries of $p$. Therefore, by another application of \lemref{same}, $(\varepsilon_0, \varepsilon_1, \varepsilon_2, \ldots)$ is one of the two $\ST$-itineraries of $z$, implying $z=\GT(t)$. 

This covers all the cases and completes the proof of \thmref{reduction}. \hfill $\Box$ \vspace{0.12in}

We conclude with the following:

\begin{corollary}[At most three points]
\label{<=3}
Let $p\in \bd \KTN$. Then $\ZT^{-1}(p)\cup \ZN^{-1}(p)$ contains at most 3 points.
\end{corollary}

\begin{pf}
Since $p$ has at most two itineraries with respect to $\Sigma_{\theta}$ and two itineraries with respect to $\Sigma_{\nu}$, \lemref{same} and \propref{itin} imply that $\ZT^{-1}(p)$ and $\ZN^{-1}(p)$ each contain at most two points. So to prove the corollary, we assume by way of contradiction that there are four distinct point $z_1, z_2\in \KT$ and $z_3, z_4\in \KN$ such that $\ZT(z_1)=\ZT(z_2)=\ZN(z_3)=\ZN(z_4)=p$. By \thmref{reduction}, all four points have to be biaccessible. Pick, for example, $z_1$ and $z_3$ and note that they eventually map to the critical points of $\TQN$ and $\TQT$ \cite{Zakeri1}. Hence $p=\ZT(z_1)$ eventually maps to the critical point $1$ of $\TBTN$ and $p=\ZN(z_3)$ also maps to the critical point $c$ of $\TBTN$. This is clearly impossible since $c$ and $1$ have disjoint orbits.
\end{pf}

\subsection{End of the proof.} 
\label{subsec:End of the proof}
We can now prove the main theorem of this paper:

\begin{theorem}[Bounded type Siegel quadratics are mateable]
\label{main}
Let $0< \theta, \nu <1$ be two irrationals of bounded type and $\theta \ne 1-\nu$. Then the quadratic
polynomials $\FT$ and $\FN$ are topologically mateable. Moreover, there exists a quadratic rational map $F$
such that $F=\FT\mate \FN$. Any two such rational maps are conjugate by a M\"obius transformation.
\end{theorem} 

\begin{pf}
The last assertion is immediate since every quadratic rational map with two fixed Siegel disks of rotation
numbers $\theta$ and $\nu$ is holomorphically conjugate to the normalized map $\FTN$ defined in
(\ref{eqn:normal form}). By Definition IIa of the introduction, it suffices to construct continuous maps
$\VT: K(\FT) \to \ov{\CC}$ and $\VN: K(\FN) \to \ov{\CC}$ with the following properties:
\begin{enumerate}
\item[(a)]
$\VT \circ \FT =\FTN \circ \VT$ and $\VN \circ \FN =\FTN \circ \VN$.
\item[(b)]
$\VT(K(\FT)) \cup \VN(K(\FN))= \ov{\CC}$.
\item[(c)]
$\VT$ and $\VN$ are conformal in the interiors of $K(\FT)$ and $K(\FN)$.
\item[(d)]
$\VT(z)=\VN(w)$ if and only if $z$ and $w$ are ray equivalent.
\end{enumerate}

It is  clear from the preceding discussion  what these maps should be. 
By the surgery construction of subsections \ref{subsec:Surgery}
and \ref{subsec:Immediate corollaries}, there exist quasiconformal homeomorphisms $\PT, \PN, \psi :
\ov{\CC} \to \ov{\CC}$ such that
\begin{equation}
\label{eqn:threemap}
\begin{array}{rcl}
\PT \circ \TQT & = & \FT \circ \PT,\\
\PN \circ \TQN & = & \FN \circ \PN,\\
\psi \circ \TBTN & = & \FTN \circ \psi.
\end{array}
\end{equation} 
Consider the semiconjugacies $\ZT$ and $\ZN$ of \thmref{reduction}
and define
$$\VT=\psi \circ \ZT \circ \PT^{-1},$$
$$\VN=\psi \circ \ZN \circ \PN^{-1}.$$ Properties (a) and (b) above are immediate consequences of the
corresponding properties of $\ZT$ and $\ZN$ stated in \thmref{reduction}. So to finish the proof, we must
show (c) and (d).

To show (c), recall the surgery construction of subsection \ref{subsec:Surgery}. Consider the Douady-Earle
extension $H_{\theta}$ used in defining the modified Blaschke product $\TQT$ in (\ref{eqn:modify}). The
invariant conformal structure $\sigma_{\theta}$ on the unit disk $\DD$ is given by the pull-back of the
standard conformal structure $\sigma_0$ under $H_{\theta}$. Similarly, we have the Douady-Earle extension
$H_{\theta, \nu}$ for the linearizing homeomorphism of $\BTN:\TT \to \TT$ used in defining the modified
Blaschke product $\TBTN$ in (\ref{eqn:modify2}), and the invariant conformal structure $\sigma_{\theta, \nu}$
on $\DD$ as the pull-back of $\sigma_0$ under $H_{\theta, \nu}$. Both $H_{\theta}$ and $H_{\theta, \nu}$
conjugate $\TQT$ and $\TBTN$ to the rigid rotation $z\mapsto e^{2 \pi i \theta}z$. By definition of $\ZT$,
we have $\ZT=H_{\theta, \nu}^{-1} \circ H_{\theta}$ on $\DD$. This means that $\ZT$ pulls $\sigma_{\theta,
\nu}$ back to $\sigma_{\theta}$ on the unit disk. It follows that the composition $\VT=\psi \circ \ZT \circ
\PT^{-1}$ on $\DD$ pulls $\sigma_0$ back to $\sigma_0$, hence it is conformal there. Then (a) and the fact
that $\FT$ and $\FTN$ are holomorphic show that $\ZT$ is conformal in the interior of $K(\FT)$. A similar
argument applies to $\ZN$.

To show (d), we note that the quasiconformal conjugacies $\PT$ and $\PN$ are conformal outside the filled
Julia sets, so they preserve the external angles. Therefore $\gamma_{\theta}=\PT \circ \GT$ and $\gamma_{\nu}=\PN \circ \GN$, where $\gamma_{\theta}$ and $\gamma_{\nu}$ are the Carath\'eodory loops of $J(\FT)$ and $J(\FN)$. By \thmref{reduction}, $\VT(z)=\VN(w)$ implies that
$z=\gamma_{\theta}(t)$ and $w=\gamma_{\nu}(-t)$ for some $t\in \TT$, which means $z$ and $w$ are ray
equivalent. The converse statement is almost immediate because if $z\in K(\FT)$ is ray equivalent to $w\in
K(\FN)$, the same is true for $\PT^{-1}(z)$ and $\PN^{-1}(w)$. Since every pair of ray equivalent points of
the form $(\GT(t), \GN(-t))$ is mapped to the same point under $(\ZT, \ZN)$, the same must be true for
arbitrary pairs of ray equivalent points. Hence $\ZT(\PT^{-1}(z))=\ZN(\PN^{-1}(w))$, or
$\VT(z)=\VN(w)$. This proves (d), and finishes the proof of the Main \thmref{main}.
\end{pf}

\vspace{0.17in}

\section{Concluding Remarks}
\label{sec:remarks}

In this section, we discuss some corollaries of \thmref{main}. In particular, we describe the nature of the pinch points
already observed in \figref{blah}. Then we prove a number-theoretic corollary of the topological mateability part of
\thmref{main} which is related to the rotation sets of the angle-doubling map on the circle. Finally, we conclude with a
discussion of the special case of a self-mating $\FT \mate \FT$ and mating $\FT$ with the Chebyshev polynomial $z\mapsto
z^2-2$.

\subsection{Ray equivalence classes and pinch points.}
\label{subsec:Ray equivalence classes and pinchings}

Consider two irrationals $\theta$ and $\nu$ of bounded type, with $\theta \neq 1-\nu$, and the quadratic polynomials $\FT$
and $\FN$ and the rational map $\FTN$. Let
$$K(\FTN)=\{ z \in \CC : \mbox{The orbit $\{ \FTN^{\circ n}(z) \}_{n \geq 0}$ never intersects $\Delta^\infty$} \}, $$ and
similarly
$$K^\infty(\FTN)=\{ z \in \CC : \mbox{The orbit $\{ \FTN^{\circ n}(z) \}_{n \geq 0}$ never intersects $\Delta^0$} \}. $$
(In \figref{blah} these two sets are the compact sets in black and gray respectively.) As we have already noted in the
introduction, $K(\FTN)$ is not a full set. In fact, it is evident from \figref{blah} that there are infinitely many
identifications between pairs of landing points of drop-chains in $K(\FTN)$ which correspond to the pinch points of
$K^\infty(\FTN)$, that is the preimages of the critical point $c\in \bd \Delta^\infty$. Similar fact holds for drop-chains of
$K^\infty(\FTN)$ and the pinch points of $K(\FTN)$. We gave a precise version of this statement in
\lemref{biaccessible}. 
It follows that every precritical point in the Julia set of  $\FT$ (resp. $\FN$) is
identified with the landing points of two distinct drop-chains of $\FN$ (resp. $\FT$). \thmref{main} allows us to
determine exactly which two drop-chains correspond to the given pinch point. Throughout the following discussion we
continue using notations from \secref{sec:proof}.

Recall that the quasiconformal conjugacies $\PT$ (between $\TQT$ and $\FT$) and $\PN$ (between $\TQN$ and $\FN$) in
(\ref{eqn:threemap}) are conformal in the basins of infinity, so they preserve the ray equivalence classes. From this fact
and \corref{<=3}, it follows that for the formal mating of $\FT$ and $\FN$, every ray equivalence class intersects $K(\FT)
\cup K(\FN)$ in at most three points. Let $E$ denote the intersection of a ray equivalence class with the union $K(\FT)
\cup K(\FN)$. We only have three possibilities for $E$:

$\bullet$ {\it Case 1.} $E=\{ z, w \}$, where $z\in K(\FT)$ and $w\in K(\FN)$ are both the landing points of unique rays,
hence $z=\gamma_{\theta}(t)$ and $w=\gamma_{\nu}(-t)$ for a unique $t\in \TT$.

$\bullet$ {\it Case 2.} $E=\{ z, z', w \}$, where $z,z' \in K(\FT)$ are both the landing points of unique rays and $w\in
K(\FN)$ is biaccessible, hence a preimage of the critical point of $\FN$. In this case, there exist $s,t \in \TT$ such that
$z=\gamma_{\theta}(s)$, $z'=\gamma_{\theta}(t)$, and $w=\gamma_{\nu}(-s)=\gamma_{\nu}(-t)$.

$\bullet$ {\it Case 3.} $E=\{ z, w, w' \}$, where $z \in K(\FT)$ is biaccessible, and $w,w' \in K(\FN)$ are both the
landing points of unique rays. In this case, there exist $s,t \in \TT$ such that $z=\gamma_{\theta}(s)=\gamma_{\theta}(t)$,
$w=\gamma_{\nu}(-t)$, $w'=\gamma_{\nu}(-s)$.

\begin{corollary}[Pinch points in $K(\FTN)$]
\label{pinch}
The compact set $K(\FTN)$ is homeomorphic to the quotient of the filled Julia set $K(\FT)$ by an equivalence relation
$\sim$ defined as follows. Two points $z\neq z'$ in $K(\FT)$ satisfy $z \sim z'$ if and only if they are the landing points
of unique rays at angles $s,t \in \TT$, $z=\gamma_{\theta}(s)$, $z'=\gamma_{\theta}(t)$, such that
$\gamma_{\nu}(-s)=\gamma_{\nu}(-t)$. Every non-trivial equivalence class of $\sim$ contains exactly two points which are
necessarily the landing points of two distinct drop-chains of $\FT$.
\end{corollary}

\begin{pf}
Since $\VT: K(\FT) \to K(\FTN)$ is a surjective map, $K(\FTN)$ is homeomorphic to $K(\FT)/\sim$, where $z \sim z'$ if and
only if $z$ and $z'$ belong to the same fiber of $\VT$. By the above discussion ({\it Case 2}), for distinct points $z\neq
z'$, we have $\VT(z)=\VT(z')$ if and only if there exist $w\in K(\FN)$ and distinct angles $s,t \in \TT$ such that
$z=\gamma_{\theta}(s)$, $z'=\gamma_{\theta}(t)$, and $w=\gamma_{\nu}(-s)=\gamma_{\nu}(-t)$. In this case $w$ is a preimage
of the critical point of $\FN$. Both $z$ and $z'$ are landing points of distinct drop-chains of $\FT$, for otherwise $z$ or
$z'$ would belong to the closure of a drop (\propref{either/or}), hence $\VT(z)=\VT(z')$ would eventually map to the
boundary of the Siegel disk $\Delta^0$ of $\FTN$. On the other hand, $\VT(z)=\VN(w)$ eventually maps to the critical point
of $\FTN$ on the boundary of $\Delta^\infty$. This would contradict $\bd \Delta^0 \cap \bd \Delta^\infty = \es$.
\end{pf}

This completely describes which identifications are made in $K(\FT)$ in order to obtain $K(\FTN)$: Take any precritical point in the Julia set of $\FN$ and calculate the angles $s, t$ of the two external rays landing on it. Then find the landing points of the external rays at angles $-s$ and $-t$ for $\FT$, which are ends of distinct drop-chains, and identify them in $K(\FT)$. This creates a ``pinch point.'' After all such possible identifications are made, we obtain a homeomorphic copy of $K(\FTN)$. Note that not all the landing points of drop-chains of $\FT$ undergo this identification, simply because there are uncountably many drop-chains and only countably many pinch points.

\subsection{Rotation sets of the doubling map.}
\label{subsec:Rotation sets of the doubling map}

The angle $\omega=\omega(\theta)$ of the external ray landing at the critical value of the quadratic polynomial $\FT$ may be described in terms of the rotation sets of the angle-doubling map on $\TT$ defined by
$m_2: x\mapsto 2x$ (mod 1). A subset $E\subset \TT$ is called a {\it rotation set} if the restriction of $m_2$ to $E$ is
order-preserving, with $m_2(E)\subset E$. It is easy to see that in this case $E$ must be contained in a closed
semicircle. Hence the restriction $m_2|_E$ can be extended to a degree $1$ monotone map of the circle, which has a
well-defined rotation number, denoted by $\rho(E)\in [0,1)$. The following theorem can be found in \cite{Bullett}:

\goodbreak
\begin{theorem}[Rotation sets of the doubling map]
\label{rotation set}

\noindent
\begin{enumerate}
\item[(i)]
For any $0\leq \theta <1$ there exists a unique compact rotation set $E_{\theta}\subset \TT$ with
$\rho(E_{\theta})=\theta$. When $\theta$ is rational $E_{\theta}$ is a single periodic orbit of $m_2$. On the other hand,
when $\theta$ is irrational, $E_{\theta}$ is a Cantor set contained in a well-defined semicircle $[\omega/2,
(\omega+1)/2]$, with $\{ \omega/2, (\omega+1)/2 \} \subset E_{\theta}$, and the action of $m_2$ on $E_{\theta}$ is
minimal. In this case the angle $\omega$ can be computed in terms of $\theta$ as
\begin{equation}
\label{eqn:OM}
\omega=\sum_{0 < p/q < \theta} 2^{-q},
\end{equation}
where the sum is taken over all (not necessarily reduced) fractions $p/q$.
\item[(ii)]
For every $0< \omega <1 $, the semicircle $[\omega/2, (\omega+1)/2]$ contains a unique compact minimal rotation set
$E^\omega$. The graph of $\omega \mapsto \rho(E^\omega)$ is a devil's staircase.
\end{enumerate}  
\end{theorem}

The mapping $\omega \mapsto \rho(E^\omega)$ is intimately connected with the parameter rays defining the limbs of the
Mandelbrot set \cite{Bullett}.

Now consider the quadratic polynomial $\FT$ for an irrational $\theta$ of bounded type. Then the Julia set $J(\FT)$ is
locally-connected, and the boundary of the Siegel disk $\Delta$ of $\FT$ is a quasicircle passing through the critical
point $0$ (compare \thmref{diam->0} and \thmref{dghs}). We know that $0$ is the landing point of exactly two external rays
at angles $\omega/2$ and $(\omega+1)/2$, where $0<\omega <1$. Define
$$E=\{ t \in \TT: \gamma_{\theta}(t) \in \bd \Delta \}.$$ It is easy to see that $E$ is compact and contained in the
semicircle $[\omega/2, (\omega+1)/2]$, hence by the above theorem, $E=E_\omega$. On the other hand, the order of the points
in the orbit $\{ \FT^{\circ n}(0) \}_{n \geq 0}$ on the boundary $\bd \Delta$ determines the rotation number $\theta$
uniquely \cite{deMelo-vanStrien}. At the same time this order coincides with the order of the orbit of $\omega$ under
$m_2$ on the circle. It follows that $\rho(E_\omega)=\theta$.

\begin{corollary}
When $0< \theta <1$ is an irrational of bounded type, the angle $0< \omega(\theta) <1$ of the external ray landing at the critical
value of the quadratic polynomial $\FT$ is given by (\ref{eqn:OM}).
\end{corollary}

It is interesting to investigate number-theoretic properties of the numbers $\omega(\theta)$ when $\theta$ is
irrational. For example, it follows from the above discussion that for irrational $0< \theta <1$, $\omega(\theta)$ is also
irrational. When $\theta$ is of bounded type, we have the much sharper statement that $\omega(\theta)$ is not
$(2+(\sqrt{5}-1)/2 -\delta)$-Diophantine for any $\delta>0$ \cite{Bullett}. In particular, by Roth's theorem,
$\omega(\theta)$ is transcendental over $\Bbb Q$. The topological mateability part of \thmref{main} allows us to draw a
further conclusion:

\begin{theorem}
\label{number theory}
Suppose that $0< \theta, \nu <1$ are irrationals of bounded type, with $\theta \neq 1-\nu$, and consider the angles
$\omega(\theta)$ and $\omega(\nu)$. Then the equation
\begin{equation}
\label{eqn:relation} 
2^n \omega(\theta)+2^m \omega(\nu) \equiv 0\ \operatorname{(mod\ 1)}
\end{equation}
does not have any solution in non-negative integers $n,m$.
\end{theorem}

\noindent
Note that the condition $\theta \neq 1-\nu$ is necessary because $\omega(\theta)+\omega(1-\theta)=1$. Also, when
$\theta=\nu$ this theorem is only saying that $\omega(\theta)$ is irrational, a fact that is clear from \thmref{rotation
set}.

\begin{pf}
Suppose that (\ref{eqn:relation}) holds for some $n,m$. Set $t=\omega(\theta)/2^m$, so that $-2^{n+m}t \equiv 2^m
\omega(\nu)$ (mod 1). Let $z=\gamma_{\theta}(t) \in J(\FT)$ and $w=\gamma_{\nu}(-t)\in J(\FN)$. Then $\FT^{\circ
m}(z)=c_{\theta}$ is the critical value of $\FT$ and $\FN^{\circ n+m}(w)=\FN^{\circ m}(c_{\nu})$ belongs to the forward
orbit of the critical point of $\FN$. By \thmref{main}, $\FTN= \FT \mate \FN$, so $\VT(z)\in J(\FTN)$ and $\VN(w)\in
J(\FTN)$ eventually hit $\bd \Delta^0$ and $\bd \Delta^\infty$ respectively. But $z$ and $w$ are ray equivalent, so
$\VT(z)=\VN(w)$ by \thmref{main}. This contradicts $\bd \Delta^0 \cap \bd \Delta^\infty =\es$.
\end{pf}

\subsection{Mating with Chebyshev quadratic polynomial.}
\label{subsec:Mating with the Chebyshev quadratic polynomial}

When $\theta=\nu$, the self-mating $F=F_{\theta, \theta}=\FT \mate \FT$ given by \thmref{main} has a natural symmetry,
i.e., it commutes with the involution ${\cal I}:z\mapsto 1/z$ of the sphere. 
As was apparently first observed by C. Petersen, if we destroy this symmetry by passing to the quotient space, we can create new examples of mating.

Consider the quotient of the Riemann sphere by the action of $\cal I$. The resulting space is again a Riemann surface
conformally isomorphic to the sphere $\ov \CC$. Since $F \circ {\cal I}={\cal I}\circ F$, there is a well-defined rational
map $G$ which makes the following diagram commute:
$$\begin{CD}
\ov \CC @> F >> \ov  \CC\\
@VV \pi V@VV \pi V\\
\ov \CC @> G >> \ov  \CC
\end{CD}$$  
Here $\pi: \ov{\CC} \to \ov{\CC}/{\cal I} \simeq \ov{\CC}$ is the degree $2$ natural projection. Chasing around this
diagram shows that $G$ is a quadratic rational map which clearly has one Siegel disk of rotation number $\theta$. Therefore
this way of collapsing the sphere identifies the two critical points of $F$ but creates a new critical point of its own. It
is not hard to check that $G$ is M\"{o}bius conjugate to the map
\begin{equation}
\label{eqn:G}
z\mapsto \frac{4z}{((1+z)+e^{2 \pi i \theta}(1-z))^2},
\end{equation}
with a fixed Siegel disk centered at $1$. The critical point $c_1=(e^{2 \pi i \theta}+1)/(e^{2 \pi i \theta}-1)$ of this
map has the finite orbit $c_1\mapsto \infty \mapsto 0$. The second critical point $c_2=-c_1$ belongs to the boundary of the
Siegel disk (compare \figref{self1}).
\realfig{self1}{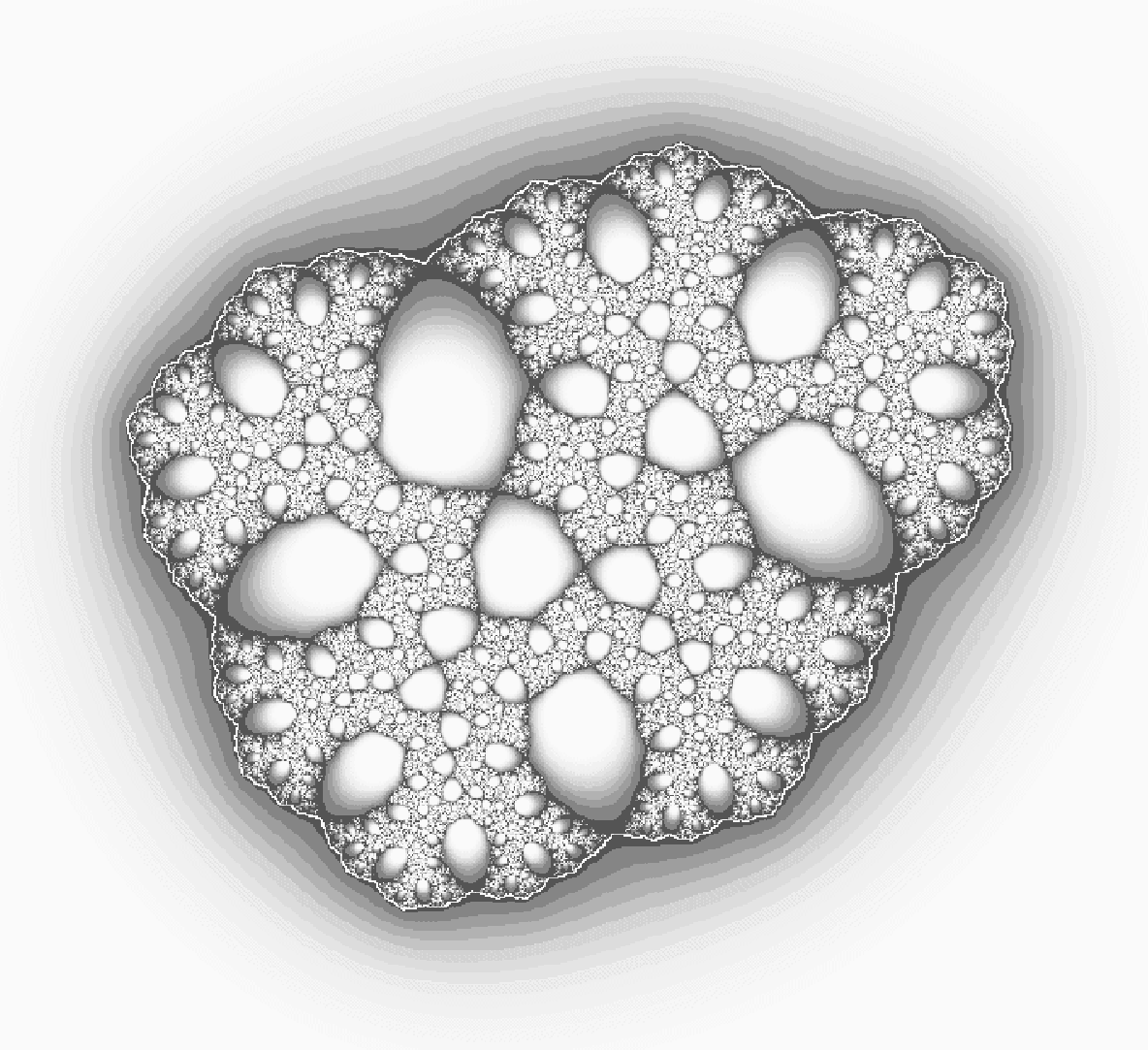}{{\sl The Julia set of the mating $\FT \mate \cheb$, where $\theta=(\sqrt{5}-1)/2$. To get a better picture we have conjugated the map in (\ref{eqn:G}) by $w=1/(z-1)$ so as to put the center of the Siegel disk at infinity and the finite critical orbit at $(e^{2 \pi i \theta}+1)/2 \mapsto 0 \mapsto -1$.}}{10cm}
 
Recall that the {\it Chebyshev} quadratic polynomial is $\cheb :z\mapsto z^2-2$. It is easy to see that the filled Julia set
$K(\cheb)=J(\cheb)$ is the closed interval $[-2,2]$. Its Carath\'eodory loop $\gcheb:\TT \to J(\cheb)$ is simply given by
$\gcheb(t)=2 \cos t$, hence $\gcheb(t)=\gcheb(s)$ if and only if $t=-s$.

We would like to show that $G$ is the mating of $\FT$ with $\cheb$. Recall that $\gamma_{\theta}$ is the Carath\'eodory
loop of $J(\FT)$ and $\VT : K(\FT) \to \ov{\CC}$ is the semiconjugacy between $\FT$ and $F$ given by \thmref{main}. Denote
by $\varphi_1$ the composition $\pi \circ \VT:K(\FT)\to \ov{\CC}$, which conjugates $\FT$ to the quadratic rational map
$G$. It is clear from the symmetry of the construction that
$$\VT( \gamma_{\theta}(-t))={\cal I} (\VT( \gamma_{\theta}(t)))$$ for all $t\in \TT$. It follows that the composition $\VT
\circ \gamma_{\theta}$ conjugates the map $t\mapsto -t$ on $\TT$ to the involution $\cal I$. Hence it descends to a map
$\varphi_2: K(\cheb)\to \ov{\CC}$ which conjugates $\cheb$ to $G$. It is easy to check that the pair $(\varphi_1,
\varphi_2)$ satisfies the conditions of Definition IIa of the introduction. Hence,

\begin{theorem}[Mating with the Chebyshev map]
\label{cheby}
Let $0< \theta <1$ be any irrational of bounded type. Then there exists a quadratic rational map $G$ such that
$$G= \FT \mate \cheb.$$ Moreover, $G$ is unique up to conjugation with a M\"obius transformation.
\end{theorem}


\vspace{0.17in}

\begin{thebibliography}{*****}

\bibitem [\bf{AB}]{Ahlfors-Bers} L. Ahlfors and L. Bers, {\it Riemann mapping's theorem for variable metrics}, Annals of Math., {\bf 72} (1960) 385-404.
\bibitem [\bf{BS}]{Bullett} S. Bullett and P. Sentenac, {\it Ordered orbits of the shift, square roots, and the devil's staircase}, Math. Proc. Camb. Phil. Soc., {\bf 115} (1994) 451-481.
\bibitem [\bf{CG}]{Carleson-Gamelin} L. Carleson and T. Gamelin, {\it Complex Dynamics}, Springer-Verlag, 1993.
\bibitem[\bf{dFdM}]{dFdM} E. de~Faria and W. de~Melo, {\it Rigidity of critical
circle mappings I}, SUNY at Stony Brook IMS Preprint, 1997/16.

\bibitem[\bf{Do1}]{Douady3} A. Douady, {\it Algorithms for computing angles in the Mandelbrot set}, in ``Chaotic Dynamics and Fractals,'' ed. Barnsley and Demko, Academic Press (1986) 155-168.
\bibitem[\bf{Do2}]{Do1} A. Douady, {\it Syst\'ems dynamiques holomorphes}, Seminar Bourbaki, Ast\'erisque, {\bf 105-106} (1983) 39-64.
\bibitem [\bf{Do3}]{Douady} A. Douady, {\it Disques de Siegel at aneaux de Herman}, Seminar Bourbaki, Ast\'erisque, {\bf 152-153} (1987) 151-172.
\bibitem [\bf{DE}]{Douady-Earle} A. Douady and C. Earle, {\it Conformally natural extension of homeomorphisms of the circle}, Acta Math., {\bf 157} (1986) 23-48.
\bibitem [\bf{DH}]{Douady-Hubbard} A. Douady and J. Hubbard, {\it A proof of Thurston's topological characterization
of rational functions}, Acta Math., {\bf 171} (1993) 263-297.
\bibitem[\bf{Ep}]{Epstein} A. Epstein, {\it Counterexamples to the quadratic mating conjecture}, Manuscript in preparation. 
\bibitem [\bf{He}]{Herman2} M. Herman, {\it Conjugaison quasisymetrique des homeomorphismes analytique des cercle a des rotations}, Manuscript.
\bibitem[\bf{Luo}]{Luo} Jiaqi Luo, {\it Combinatorics and holomorphic dynamics: Captures, matings, Newton's method},
Thesis, Cornell University, 1995.
\bibitem[\bf{Lyu}]{L3} M.Yu. Lyubich, {\it The dynamics of rational transforms: The topological picture}, Russian Math. Surveys {\bf 41} (1986) 43-117.
\bibitem [\bf{Mc}]{McMullenbook} C. McMullen, {\it Complex Dynamics and Renormalization}, Annals of Math Studies, vol. 135, 1994.
\bibitem[\bf{dMvS}]{deMelo-vanStrien} W. de Melo, S. van Strien, {\it One-dimensional dynamics}, Springer-Verlag, 1993.

\bibitem [\bf{Mi1}]{Milnor1} J. Milnor, {\it Dynamics in One Complex Variable: Introductory Lectures}, SUNY at Stony Brook IMS preprint, 1990/5.
\bibitem[\bf{Mi2}]{Milnor-remarks} J. Milnor, {\it Geometry and dynamics of quadratic rational maps}, Experim. Math., {\bf 2} (1993) 37-83.
\bibitem [\bf{Mi3}]{Milnor2} J. Milnor, {\it Periodic orbits, external rays, and the Mandelbrot set: An expository account}, Available at http://www.math.sunysb.edu/ $\tilde{}$ jack.
\bibitem [\bf{Mi4}]{Milnor4} J. Milnor, {\it Pasting together Julia sets - a worked out example of mating}, To appear.
\bibitem [\bf{Mo}]{Moore} R.L. Moore, {\it Concerning upper semi-continuous collection of continua}, Trans. Amer. Math. Soc., {\bf 27} (1925) 416-428.
\bibitem [\bf{Pe}]{Petersen} C. Petersen, {\it Local connectivity of some Julia sets containing a circle with an irrational rotation}, Acta Math., {\bf 177} (1996) 163-224.
\bibitem [\bf{Re1}]{Rees1} M. Rees, {\it Realization of matings of polynomials as rational maps of degree two}, Manuscript, 1986.
\bibitem[\bf{Re2}]{Rees2} M. Rees, {\it A partial description of parameter space of rational maps of degree two:
part I}, Acta Math., {\bf 168} (1992) 11-87.
\bibitem[\bf{Sh}]{Shishikura} M. Shishikura, {\it On a theorem of M. Rees for matings of polynomials},
Preprint, IHES, 1990.
\bibitem[\bf{ST}]{ST} M. Shishikura and  L. Tan, {\it A family of cubic rational maps and matings of cubic polynomials},
Preprint, Max Planck Inst., 1988.
\bibitem [\bf{Sw}]{Swiatek} G. Swiatek, {\it Rational rotation numbers for maps of the circle}, Comm. Math. Phys., {\bf 119} (1988) 109-128.
\bibitem[\bf{Tan}]{Tan} L. Tan, {\it Matings of quadratic polynomials}, Erg. Th. and Dyn. Sys. {\bf 12} (1992) 589-620.
\bibitem[\bf{TY}]{Tan-Yin} L. Tan and Y. Yin, {\it Local connectivity
of the Julia set for geometrically finite rational maps}, Preprint
\'Ec. Norm. Sup. de Lyon, UMPA-94-${\sl N}^{\underline{\circ}}$121, 1994.
\bibitem [\bf{Ya}]{Yampolsky} M. Yampolsky, {\it Complex bounds for renormalization of critical circle maps}, To appear in Erg. Th. and Dyn. Sys. (early version available as SUNY at Stony Brook IMS preprint, 1995/12).
\bibitem [\bf{Yo1}]{Yoccoz1} J.C. Yoccoz, {\it Il n'y a pas de contre-example de Denjoy analytique}, C. R. Acad. Sci. Paris, {\bf 298} (1984) 141-144.
\bibitem [\bf{Yo2}]{Yoccoz2} J.C. Yoccoz, {\it Petits Diviseurs en Dimension 1}, Ast\'erisque {\bf 231}, 1995.
\bibitem [\bf{Za1}]{Zakeri1} S. Zakeri, {\it Biaccessibility in quadratic Julia sets I-II}, SUNY at Stony Brook IMS preprint, 1998/1.
\bibitem [\bf{Za2}]{Zakeri2} S. Zakeri, {\it On dynamics of cubic Siegel polynomials}, SUNY at Stony Brook IMS preprint, 1998/4.

\end{thebibliography}
\end{document}